\documentclass[leqno,a4paper,11pt]{article}

\usepackage{hyperref}
\usepackage{amsthm}
\usepackage{amsmath}
\usepackage{amssymb}
\usepackage{graphicx}
\usepackage{centernot}
\usepackage[xcolor]{changebar}
\usepackage{enumitem}

\newcommand{\Con}{\ensuremath{\mathcal{C}}}

\newcommand{\loc}{\ensuremath{\text{loc}}}

\newcommand{\mb}[1]{\ensuremath{\mathbb{#1}}}
\newcommand{\N}{\mb{N}}

\newcommand{\R}{\mb{R}}

\newcommand{\grad}{\ensuremath{\mbox{\rm grad}\,}}

\newfont{\bl}{msbm10 scaled \magstep2}

\newcommand{\beq}{\begin{equation}}
\newcommand{\eeq}{\end{equation}}
\newcommand{\notmid}{\mid\kern-0.5em\not\kern0.5em}

\newcommand{\eps}{\varepsilon}

\newcommand{\supp}{\mathop{\mathrm{supp}}}

\theoremstyle{theorem}
\newtheorem{thm}{Theorem}[section]
\newtheorem{lem}[thm]{Lemma}
\newtheorem{prop}[thm]{Proposition}
\newtheorem{cor}[thm]{Corollary}
\newtheorem{defi}[thm]{Definition}

\theoremstyle{definition}
\newtheorem{rem}[thm]{Remark}
\newenvironment{pr}{\begin{proof}[\textbf{Proof:}] \ }{\end{proof}}
\newtheorem{ex}[thm]{Example}

\newcommand{\LLS}{Lorentzian length space }
\newcommand{\LLSn}{Lorentzian length space}
\newcommand{\LpLS}{Lorentzian pre-length space }
\newcommand{\LpLSn}{Lorentzian pre-length space}
\newcommand{\Xll}{$(X,d,\ll,\leq,\tau)$ }
\newcommand{\Xlln}{$(X,d,\ll,\leq,\tau)$}
\newcommand{\nll}{\centernot\ll}
\newcommand{\Tau}{\mathcal{T}}
\newcommand{\Lip}{\mathrm{Lip}}

\title{Lorentzian length spaces}
 \author{Michael Kunzinger\thanks{{\tt michael.kunzinger@univie.ac.at}, Faculty of Mathematics, University of Vienna,
 Austria.},
Clemens S\"amann\thanks{{\tt clemens.saemann@univie.ac.at}, Faculty of Mathematics, University of Vienna,
 Austria.}}

\begin{document}
\date{}
 \maketitle

 \begin{abstract}
We introduce an analogue of the theory of length spaces into the setting of Lorentzian geometry
and causality theory. The r\^ole of the metric is taken over by the time separation function,
in terms of which all basic notions are formulated. In this way we recover many fundamental
results in greater generality, while at the same time clarifying the minimal requirements for
and the interdependence of the basic building blocks of the theory. A main focus of this work
is the introduction of synthetic curvature bounds, akin to the theory of Alexandrov and CAT$(k)$-spaces,
based on triangle comparison. Applications include Lorentzian manifolds with metrics of low regularity,
closed cone structures, and certain approaches to quantum gravity. 
\vskip 1em

\medskip

 \end{abstract}

 \tableofcontents

%\nochangebars %uncomment for turning off the changebars

\section{Introduction}\label{sec-intro}

Metric geometry, and in particular the theory of length spaces, is a vast and very active field
of research that has found applications in diverse mathematical disciplines, such as differential 
geometry, group theory, dynamical systems and partial differential equations. It has led
to identifying the `metric core' of many results in differential geometry, to clarifying the interdependence
of various concepts, and to generalizations of central notions in the field to low regularity situations.
In particular, the synthetic approach to curvature bounds in the theory of Alexandrov spaces
and CAT$(k)$-spaces has turned out to be of fundamental importance (cf., e.g., \cite{BH:99,BBI:01,Pap:14}). 

The purpose of this work is to lay the foundations for a synthetic approach to Lorentzian 
geometry that is similar in spirit to the theory of length spaces and that, in particular,
allows one to introduce curvature bounds in this general setting. The motivation for our 
approach is very similar to that of metric geometry: ideally, it should  identify
the minimal requirements for and the logical dependence among fundamental results of Lorentzian geometry
and causality theory, and in this way separate the essential concepts from various
derived notions. Based on this, one may extend the validity of these results to their
natural maximal range, in particular to settings
where the Lorentzian metric is no longer differentiable, or even to situations where there is
only a causal structure not necessarily induced by a metric. Again in parallel to the
case of metric geometry, appropriate notions of synthetic (timelike or causal) 
curvature bounds based on triangle comparison occupy a central place in this development, 
leading to a minimal framework for Lorentzian comparison geometry. 
In the smooth case, triangle comparison methods were
pioneered by Harris in \cite{Har:82} for the case of timelike triangles in Lorentzian
manifolds and by Alexander and Bishop in \cite{AB:08} for the general semi-Riemannian
case and triangles of arbitrary causal character. The notions 
introduced in this paper are designed to be compatible with these works, while at
the same time introducing curvature bounds even to settings where there is no
curvature tensor available (due to low differentiability of the metric or even the
absence thereof). On the one hand, this makes it possible to establish well-known results
from metric geometry also in the Lorentzian setting (e.g.\ non-branching of maximal
curves in spaces with timelike curvature bounded below, cf.\ Theorem \ref{nonbranch_th}). 
On the other hand, it provides
a new perspective on genuinely Lorentzian phenomena, like the push-up principle
for causal curves, which in the present context appears as a consequence
of upper causal curvature bounds (see Section \ref{sec-causal-curv-bounds}).

The r\^ole of the metric of a length space in the current framework will be played by the 
time-separation function $\tau$, which will therefore be our main object of study. It is
closely linked to the causal structure of Lorentzian manifolds, and in fact for
strongly causal spacetimes it completely determines the metric by a classical result
of Hawking, cf.\ the beginning of Subsection \ref{subsec-lor-pre-lenght-spaces} 
below. We will  express all
fundamental notions of Lorentzian (pre-)length spaces (such as length and maximality of curves, 
curvature bounds, etc.) in terms of $\tau$. It turns out that this provides a satisfactory 
description of much of standard causality theory and recovers many fundamental results
in greater generality.

Apart from the intrinsic interest in a Lorentzian analogue of metric geometry, a main motivation
for this work is the necessity to consider Lorentzian metrics of low regularity. This need
is apparent both from the PDE point of view on general relativity, i.e., the initial-value-problem for the Einstein 
equations, and from physically relevant models. In fact, the standard local existence result for the vacuum Einstein 
equations (\cite{Ren:05}) assumes the metric to be of Sobolev-regularity $H^s_\loc$ (with $s>\frac{5}{2}$).  
Recently, the regularity of the metric has been lowered even further (e.g.\ \cite{KRS:15}). In many cases, spacetimes describing 
physically relevant situations require certain restrictions on the regularity of the metric. In particular, modeling 
different types of matter in a spacetime may lead to a discontinuous energy-momentum tensor, and thereby via the Einstein 
equations to metrics of regularity below $\Con^2$ (\cite{Lic:55, MS:93}). Prominent examples are spacetimes that model the 
inside and outside of a star or shock waves. Physically relevant models 
of even lower regularity include  spacetimes with 
conical singularities and cosmic strings (e.g.\ \cite{Vic:90, VW:00}), (impulsive) gravitational waves (see e.g.\ \cite{Pen:72a}, 
\cite[Ch.\ 20]{GP:09} and \cite{PSSS:14, PSSS:16, SSS:16} for more recent works), and ultrarelativistic black holes (e.g.\ 
\cite{AS:71}).

There has in fact recently been a surge in activity in the field of low-regularity Lorentzian
geometry and mathematical relativity. Some main trends in this branch of research are the studies in 
geometry and causality theory for Lorentzian metrics of low regularity (see \cite{CG:12, Sae:16}
for results on continuous metrics and \cite{Min:15,KSS:14,KSSV:14} for the $\Con^{1,1}$-setting).
Cone structures on differentiable manifolds are another natural generalization of smooth Lorentzian
geometry, and several recent fundamental works have led to far-reaching extensions of causality
theory, see \cite{FS:12,BS:18,Min:17}. Another line of research concerns the extension of the
classical singularity theorems of Hawking and Penrose to minimal regularity (\cite{KSSV:15, KSV:15, GGKS:18}).
In fact, both spacetimes with metrics of low regularity and closed cone structures
of certain types provide large natural classes of examples within the theory of Lorentzian length spaces
(cf.\ Section \ref{sec-app}).
Finally, we mention the recent breakthroughs concerning the $\Con^0$-inextendibility of 
spacetimes, pioneered by J.\ Sbierski in \cite{Sbi:18} in the case of the Schwarzschild solution
to the Einstein equations, followed by further investigations
by various authors, cf.\ \cite{GLS:18, GL:17, DL:17}. As we shall point out repeatedly throughout 
this work, there are close ties between these works and the theory of Lorentzian length spaces.
 More specifically, the follow-up work \cite{GKS:18} 
\begin{itemize}
\item introduces a notion of extendibility for Lorentzian (pre-)length spaces
that reduces to isometric embeddability in the particular case of spacetimes,
\item gives a characterization of timelike completeness in terms of the time separation function,
\item shows that timelike completeness in this sense implies inextendibility even in this general setting, and
\item for the first time, relates inextendibility to the occurrence of (synthetic) curvature
singularities.
\end{itemize}
This complements and extends the works cited above in several directions. In particular, both the
original spacetime and the extension are now allowed to be of low regularity and indeed also
non-manifold extensions can be considered.

%In particular, the methods introduced here can 
%aid in the understanding of low-regularity extendibility of spacetimes, cf.\ \cite{GKS:18}.

There have been several approaches to give a synthetic or axiomatic description of (parts of) Lorentzian geometry and 
causality. However, except the work of Busemann \cite{Bus:67} they were not in the spirit of metric geometry and length 
spaces. Moreover, triangle comparison was never used in these settings. In particular, Busemann --- a pioneer of length 
spaces --- introduced so-called \emph{timelike spaces} in \cite{Bus:67}, but his approach was too restrictive to even capture 
all smooth (globally hyperbolic) Lorentzian manifolds. Another closely related work is due to Kronheimer and Penrose (\cite{KP:67}), 
who studied the properties and 
interdependences of the causal relations in complete generality. Similar in spirit, Sorkin and Woolgar 
\cite{SW:96} established that using order-theoretic and topological methods one can extend specific results from causality 
theory to spacetimes with continuous metrics. On the other hand, Martin and Panangaden showed in \cite{MP:06} how one recovers a 
spacetime from just the causality relations in an order-theoretic manner, thereby indicating applications to quantum gravity. 
For this reason Lorentzian (pre-)length spaces might provide a fundamental framework to approaches to quantum gravity, as 
outlined in Subsection \ref{subsec-app-qg}. In a series of works Borchers and Sen \cite{BS:90, BS:99, BS:06} 
describe an approach to causality via an axiomatic notion of \emph{light rays} corresponding to null geodesics. Finally, and very relevant to the goals of this paper, Sormani and Vega \cite{SV:16} have recently introduced a metric structure on (smooth) spacetimes via what they call a null distance function. It is defined as the infimum over \emph{null lengths}, which in turn are total variations of a time function along concatenations of (future or past directed) causal curves. The null distance provides a conformally invariant pseudo-metric, and under some natural assumptions on the spacetime even a definite distance function inducing the manifold topology. This leads to an alternative starting point for studying metric geometry in the Lorentzian setting.

The plan of the paper is as follows: In the remainder of this section we
fix some basic notions used throughout this work. Section \ref{sec-bas} introduces the
fundamental causal and chronological relations, the definition of Lorentzian pre-length spaces
in terms of a lower semi-continuous time-separation function, and the associated topological notions.
It also includes the fundamentals of causal curves. To obtain a satisfying causality theory close
to that of smooth spacetimes, more structure is required. This leads, in Section \ref{sec-lorentzian-length-spaces},
to the notion of Lorentzian length space. We show that in this setting, most of standard 
causality theory remains intact, including limit curve theorems, push-up principles for causal curves, (a significant part of)
the causal ladder, and the Avez-Seifert theorem for globally hyperbolic spaces.
We also introduce a synthetic notion of regularity for \LpLSn s.  Section \ref{sec-curvature}
is devoted to synthetic curvature bounds in terms of triangle comparison. Time separation in timelike and
causal geodesic triangles is compared with corresponding quantities in Lorentzian model spaces of constant
curvature. In the smooth case, this is compatible with the Toponogov-type results of Alexander and Bishop
in \cite{AB:08}. Applications include a non-branching theorem for timelike curves in spaces with 
timelike curvature bounded below (Theorem \ref{nonbranch_th}) and a new interpretation of length-increasing
push-up in terms of upper bounds on causal curvature (Theorem \ref{thm-length-increasing-push-up_curvature}). 
We then go on to defining synthetic curvature singularities and show that the central singularity of the interior 
Schwarzschild solution can be detected via timelike triangle comparison. The final Section \ref{sec-app}
provides an in-depth study of probably the most important class of examples, namely continuous Lorentzian metrics.
Apart from identifying Lorentzian length spaces among continuous spacetimes we also derive
several new results on continuous causality theory in this section.
In addition, we consider closed cone structures following the recent fundamental work of Minguzzi
(\cite{Min:17}) and give a brief outlook on potential applications
in certain approaches to quantum gravity, namely causal Fermion systems and the theory of causal sets. 

To conclude this introduction we introduce some basic notions and fix notations. For terminology
from or motivated by causality theory we will follow the standard texts \cite{ONe:83,BEE:96,MS:08,Chr:11},
see also \cite{CG:12} for the case of continuous Lorentzian metrics. 
We will usually only formulate the future-directed case of our results, with the understanding that
there always is an analogous past-directed statement.

Among the main applications of the theory developed here will be spacetimes $(M,g)$, 
where $M$ is a differentiable manifold and $g$ is a continuous or smooth 
Lorentzian metric. We always assume that $(M,g)$ is time-oriented (i.e., there is a continuous timelike vector field on $M$) and 
we fix a smooth (without loss of generality) complete Riemannian metric $h$ on $M$, denoting the induced metric by $d^h$. 

Recall that for (continuous) Lorentzian metrics $g,\ \hat{g}$ on $M$,
$g\prec\hat{g}$ is defined as the property that the lightcones of $\hat{g}$ are strictly wider than those of $g$, i.e., if 
a non-zero vector is causal for $g$ then it is timelike for $\hat{g}$. Additionally, we also use the non-strict version, 
i.e., we mean by $g\preceq \hat{g}$ that every $g$-causal vector is causal for $\hat{g}$. 

Causal curves in spacetimes are locally Lipschitz continuous curves $\gamma\colon I\rightarrow M$ with $g(\dot\gamma,\dot\gamma)\leq 0$ 
almost everywhere. Analogously for timelike and null curves and their time-orientation.

\section{Basics}\label{sec-bas}

\subsection{Causal relations}\label{subsec-cr}

We start our analysis by introducing a slightly more general notion of causal spaces, as compared to \cite{KP:67}.

\begin{defi}\label{def-causal-space}
 Let $X$ be a set with a reflexive and transitive relation $\leq$ (a \emph{preorder}) and $\ll$ a transitive relation 
contained in $\leq$ (i.e., $\ll\,\subseteq\,\leq$, or more explicitly: if $x\ll y$ then $x\leq y$). Then $(X,\ll,\leq)$ is 
called a \emph{causal space}. % (cf.\ \cite{KP:67}, where a more specific notion of causal space is introduced). 
We write $x<y$ if $x\leq y$ and $x\neq y$.
\end{defi}

\begin{ex}\label{ex-cau-spa-con}
 Any spacetime with a continuous metric and the usual causal relations (e.g.\ $x\ll y$ if there is a future-directed 
timelike curve from $x$ to $y$) is a causal space. Contrary to \cite{KP:67} we do not require any causality 
condition to hold. Thus also the Lorentz cylinder $M=S^1_1\times\R$ (\cite[Example 5.35]{ONe:83}) is an example, 
and in this case $\ll\ =\ \leq\ = M\times M$. We will study spacetimes with continuous metrics in 
greater detail in subsection \ref{subsec-cont_lor}.
\end{ex}

\begin{defi}(Futures and pasts)
 For $x\in X$ define
 \begin{enumerate}[label=(\roman*)]
  \item $I^+(x):=\{y\in X: x\ll y\}$ and $I^-(x):=\{y\in X: y\ll x\}$, 
  \item $J^+(x):=\{y\in X: x\leq y\}$ and $J^-(x):=\{y\in X: y\leq x\}$.
 \end{enumerate}
\end{defi}

\subsection{Topologies on causal spaces}\label{subsec-top}
On a causal space one can define two natural topologies as follows. Let $(X,\ll,\leq)$ be a causal space. For $x,y\in X$ 
define $I(x,y):=I^+(x)\cap I^-(y)\subseteq X$.

\begin{defi} $($Topologies on $(X,\ll,\leq))$
Let $(X,\ll,\leq)$ be a causal space.
\begin{enumerate}[label=(\roman*)]
 \item  Define a topology $\mathcal{A}$ on $X$ by using $S:= \{I(x,y): 
x,y\in X\}$ as a subbase. We call this topology the \emph{Alexandrov topology} on $X$ with respect to $\ll$.
 \item Define a topology $\mathcal{I}$ on $X$ by using $P:=\{I^\pm(x):x\in X\}$ as a subbase. We call this topology the 
\emph{chronological topology} on $X$.
\end{enumerate}
\end{defi}

\begin{ex}\label{ex-graph} ($S$ and $P$ are not bases for topologies)
In general the sets $S$ and $P$ do not form bases for topologies as this simple example shows. Let $X=\{1,\ldots,7\}$ and let
$\ll$ be given via the graph below (e.g., $1\ll7$ etc.). This is a transitive and irreflexive relation.
\begin{figure}[h!]
  \begin{center}
\includegraphics[width=72mm, height= 48mm]{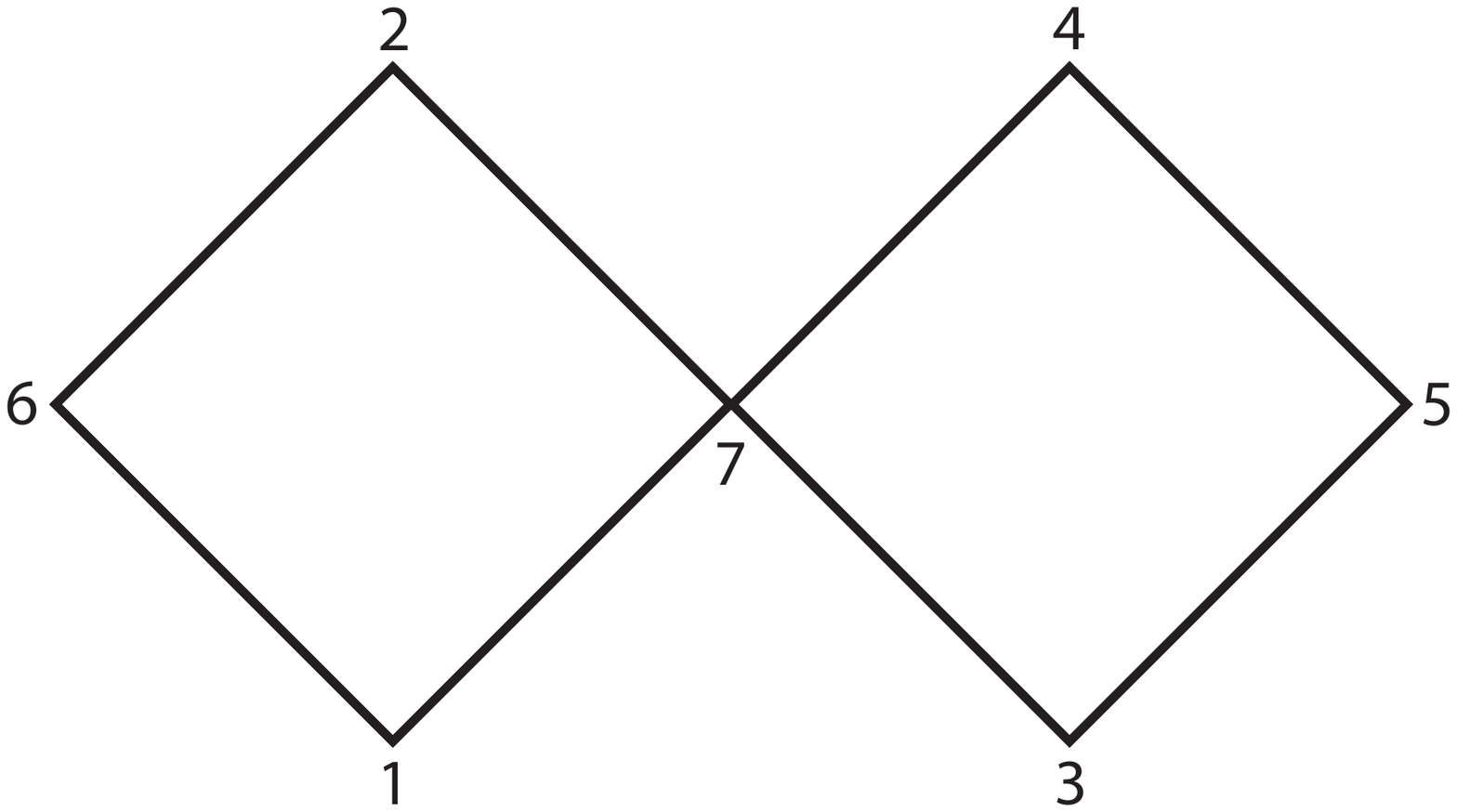}
\end{center}
\caption{Graph for Example \ref{ex-graph}}
\end{figure}

Then:
\begin{enumerate}[label=(\roman*)]
 \item The point $1$ is not covered by any $I(x,y)$ for $x,y\in X$, thus $S$ does not cover $X$.
 \item Although $P$ does cover $X$, it does not have the second property required 
for a basis: $I^+(1)\cap I^-(2) = \{6,7\}$ but there is no $x\in X$ such that $7\in 
I^\pm(x)\subseteq \{6,7\}$.
\end{enumerate}

\end{ex}

\begin{rem}\label{rem-Ip-open}
 Note that in general the future and past $I^\pm(x)$ is not open with respect to $\mathcal{A}$. As an example, take $X=\{1,2\}$ with 
$\ll:=\{(1,2)\}$, then $S=\{\emptyset\}$ hence $\mathcal{A}=\{\emptyset, X\}$ is the trivial topology and $I^+(1)=\{2\}$.
\end{rem}

\begin{prop}\label{prop-I-finer-A}
 The chronological topology $\mathcal{I}$ is finer than the Alexandrov topology $\mathcal{A}$. Consequently a map from $X$ 
into a topological space $Y$ that is $\mathcal{A}$-continuous is also $\mathcal{I}$-continuous. In particular, for 
$Y=[0,\infty]$ the same holds for semicontinuous maps.
\end{prop}
\begin{pr}
 The Alexandrov topology $\mathcal{A}$ is the coarsest topology containing $S$ and thus since $\mathcal{I}$ contains $S$ we 
have $\mathcal{A}\subseteq\mathcal{I}$.
\end{pr}
The example in \ref{rem-Ip-open} shows that in general $\mathcal{I}$ is strictly finer than $\mathcal{A}$. In that case 
$\mathcal{A}=\{\emptyset, X\}$ and $\mathcal{I}$ is the discrete topology.
\bigskip

In what follows we will not work directly with one of these topologies, but mainly use them for
comparison. Instead, 
we will assume that the causal space $X$ also comes with the structure of a metric space.

\subsection{Lorentzian pre-length spaces}\label{subsec-lor-pre-lenght-spaces}
We now introduce the central object of our study, namely the time separation function $\tau$,
in terms of which all subsequent concepts will be formulated. While it is evident that
in the smooth setting $\tau$ is closely linked to the causal structure, 
it is worth noting that by a classical result due to Hawking, King and McCarthy \cite{HKM:76}
for strongly causal spacetimes it in fact completely determines the
metric,  (cf.\ \cite[Prop.\ 3.34]{MS:08} or \cite[Thm.\ 4.17]{BEE:96}).

\begin{defi}\label{def-LpLS}
 Let $(X,\ll,\leq)$ be a causal space and $d$ a metric on $X$. Let $\tau\colon X\times X\rightarrow[0,\infty]$ be 
a lower semicontinuous map (with respect to the metric topology induced by $d$) that satisfies
\begin{equation}\label{eq-rev-tri}
 \tau(x,z)\geq \tau(x,y) + \tau(y,z)\,,
\end{equation}
for all $x,y,z\in X$ with $x\leq y\leq z$. Moreover, suppose that $\tau(x,y)=0$ if $x\nleq y$ and $\tau(x,y)>0 \Leftrightarrow x\ll y$.
Then \Xll is called a \emph{Lorentzian pre-length space} and $\tau$ is called 
the \emph{time separation function}.
\end{defi}
Since we now use the metric $d$ on $X$, all topological notions refer to the metric topology $\mathcal{D}$ induced by $d$, unless specified otherwise.

\begin{rem} It would be logically possible to introduce \LpLSn s based only on a set
endowed with a pre-order $\le$ and then \emph{define} the timelike relation $\ll$ via 
$x\ll y :\Leftrightarrow \tau(x,y)>0$. However, we prefer to view this condition
as a form of compatibility between the time separation function and the causal 
space. For an example where this compatibility is violated see Example \ref{bubble-ex}.
\end{rem}

\begin{lem}(Push-up)\label{lem-pup}
Let \Xll be a \LpLS and let $x,y,z\in X$ with $x\leq y\ll z$ or $x\ll y\leq z$. Then $x\ll y$.
\end{lem}
\begin{pr}
 Let $x\leq y\ll z$ or $x\ll y\leq z$. Then $\tau(x,z)\geq \tau(x,y)+\tau(y,z)>0$, 
%since both terms are non-negative and one of them is positive, 
which implies $x\ll z$.
\end{pr}

\begin{ex}\label{ex-lpls-smo}
 Let $(M,g)$ be a smooth spacetime with its canonical causal relations $\ll$ and $\leq$. Then by Example \ref{ex-cau-spa-con} 
$(M,\ll,\leq)$ is a causal space. The (classical) time separation function $\tau$ is lower semicontinuous with 
respect to the manifold topology (\cite[Lemma 14.17]{ONe:83}), which is induced by any Riemannian metric 
$h$ on $M$ and its associated metric $d^h$. Moreover, $\tau$ satisfies the reverse triangle inequality \eqref{eq-rev-tri} 
and $\tau(p,q)>0 \Leftrightarrow p\ll q$, thus $(M,d^h, \ll,\leq,\tau)$ is a Lorentzian pre-length space. 
Note that, in general, for a spacetime with a continuous metric the time separation function $\tau$ need not be lower semicontinuous, 
cf.\ Example \ref{bubble-ex} below.
\end{ex}

We will establish some basic facts about \LpLSn s below.

\begin{lem}\label{lem-Ip-open}
 Let \Xll be a \LpLS and $x\in X$. Then $I^\pm(x)$ is open.
\end{lem}
\begin{pr}
We establish this fact only for $I^+(x)$, the case for $I^-(x)$ works in complete analogy. Let $y\in I^+(x)$, 
so $x\ll y$ and thus $\tau(x,y)>0$. By the lower semicontinuity of $\tau(x,.)$ % (using $\eps:=\frac{\tau(x,y)}{2}>0$)
there is a neighborhood $U$ of $y$ such that $\tau(x,z)>\frac{\tau(x,y)}{2}>0$ for all $z\in U$. Consequently $x\ll z$ and 
thus $y\in U\subseteq I^+(x)$.
\end{pr}

\begin{prop}
 Let \Xll be a Lorentzian pre-length space. Then the relation $\ll$ is an open subset of $X\times X$.
\end{prop}
\begin{pr}
 Clearly, $\ll = \{(x,y)\in X\times X: \tau(x,y)>0\}$, which is open by the lower semicontinuity of $\tau$.
\end{pr}

\begin{prop}\label{prop-tau-0-inf}
 Let \Xll be a Lorentzian pre-length space. Then for every $x\in X$, either $\tau(x,x)=0$ or $\tau(x,x)=\infty$. Moreover, if $\tau(x,y)\in (0,\infty)$, then $\tau(y,x)=0$.
\end{prop}
\begin{pr}
 Let $x\in X$ with $\tau(x,x)<\infty$. By the reflexivity of $\leq$ we have $x\leq x \leq x$ and thus by the reverse 
triangle inequality we have 
$\tau(x,x) + \tau(x,x) \leq \tau(x,x)$, which implies that $\tau(x,x)\leq 0$ and so $\tau(x,x)=0$.
Finally, let $0<\tau(x,y)<\infty$ and suppose that $\tau(y,x)>0$. Then $x\ll y\ll x$,
which implies $x\ll x$ and by the above $\tau(x,x)=\infty$. But then $\tau(x,x)+\tau(x,y)
\le \tau(x,y)<\infty$ gives a contradiction.
\end{pr}

\begin{rem}
The above proposition shows that a set $X$ with a transitive relation $\ll$ that is not irreflexive
(i.e., $x\ll x$ for some $x$) cannot have a finite 
valued time separation function with respect to $\ll$ (and any transitive and reflexive relation $\leq$ on $X$ containing 
$\ll$). 
In fact, $x\ll x$ implies $\tau(x,x)>0$, but by the above $\tau(x,x)\leq 0$, if $\tau$ were finite valued.
\end{rem}

\begin{ex}\label{ex-lpls-gra}
A finite directed graph $(V,E)$ can be given the structure of a \LpLSn. Here $X$ is any finite set, $V\subseteq X$ is the 
set of vertices and $E\subseteq V\times V$ is the (directed) set of edges, i.e., $(x,y)\in E$ if and only if there is an 
edge from $x$ to $y$. 
Now define $x\ll y$ if $(x,y)\in E$ and $x\leq y$ if $x\ll y$ or $x=y$. This gives a causal space $(V,\ll,\leq)$, where 
$\ll$ is irreflexive if and only if $(V,E)$ is a directed acyclic graph. Define $\tau(x,y):= \sup \{ |C|: C $ a directed 
subgraph from $x$ to $y \}$, if such a subgraph exists and $|C|$ denotes its (finite) cardinality; otherwise set 
$\tau(x,y):=0$. Since the only metrizable topology on a finite set is the discrete 
topology, we let $d$ be the discrete metric on $V$. Thus $\tau$ is continuous and it satisfies the reverse triangle 
inequality \eqref{eq-rev-tri}. This yields that $(V,d,\ll,\leq,\tau)$ is a \LpLSn.

Note that since the topology is discrete the space is totally disconnected and hence there are no causal curves (cf.\ 
Definition \ref{def-cau-cur}, below). Causal curves will be essential in the development of the theory of \LLSn s and it will 
be a requirement on such spaces that causally related points can be connected by a (non-constant) continuous causal curve (cf.\ 
Definition \ref{def-cau-pc}). This rules out finite and, in fact, countable Lorentzian length spaces, 
as the only metrizable topology in these
cases is totally disconnected. Thus the situation is the same as for length spaces --- they are path connected and hence 
cannot be countable.

This example is closely related to the theory of \emph{causal sets}, an approach to quantum gravity, see Subsection \ref{subsec-app-qg}.
\end{ex}

\subsection{Topologies on Lorentzian pre-length spaces}

We want to relate the two natural topologies $\mathcal{I}$ and $\mathcal{A}$ to the given metric topology $\mathcal{D}$.

One can see from the proof of Lemma \ref{lem-Ip-open} that the only property used for the topology 
is that $\tau$ be semicontinuous with respect to it.
Thus if one has a function $\tau\colon X\times X\rightarrow [0,\infty]$ that satisfies the properties of a time separation 
function (as in Definition \ref{def-LpLS}) except that it is lower semicontinuous with respect to $\mathcal{A}$ (and not necessarily with respect to $d$), then 
$I^\pm(x)$ is $\mathcal{A}$-open for every $x\in X$. Consequently, $\mathcal{A}$ is finer than $\mathcal{I}$ and so by 
Proposition \ref{prop-I-finer-A} we have $\mathcal{I}=\mathcal{A}$. On the other hand the metric topology $\mathcal{D}$ is 
finer than both these topologies since by the extension of Lemma \ref{lem-Ip-open} 
mentioned above $I^\pm(x)$ and $I^+(x)\cap I^-(y)$ are $\mathcal{D}$-open 
for every $x,y\in X$. This yields that
\begin{enumerate}[label=(\roman*)]
 \item in an ``$\mathcal{A}$-Lorentzian pre-length space'' (i.e., $\tau$ lower semicontinuous with respect to $\mathcal{A}$) 
we have $\mathcal{A}=\mathcal{I}$ and thus
 \item every ``$\mathcal{A}$-Lorentzian pre-length space'' is also an ``$\mathcal{I}$-Lorentzian 
pre-length space'', which is also a Lorentzian pre-length space (since if $\tau$ is lower semicontinuous with 
respect to $\mathcal{I}$, then it is lower semicontinuous with respect to $\mathcal{D}$).
\end{enumerate}
Summing up, by using the additional freedom via the metric $d$ we (possibly) enlarge the number of Lorentzian pre-length 
spaces, as can be seen from the example below.

\begin{ex}
 In light of Remark \ref{rem-Ip-open} and the above we know that there cannot be a time separation function $\tau$ on 
$(X=\{1,2\},\ll=\{(1,2)\},\leq=X\times X)$ that is lower semicontinuous with respect to $\mathcal{A}$. This can be also seen 
directly from the fact that a prospective time separation function $\tau$ would need to satisfy $\tau(1,1)=\tau(2,2) = 
\tau(2,1)=0$ and $\tau(1,2)>0$. However, the only $\mathcal{A}$-open sets in $X\times X$ are $\emptyset$ and $X\times X$ 
and thus $\{(x,y)\in X\times X: \tau(x,y)>0\} = \{(1,2)\}$ is not $\mathcal{A}$-open (which would be required if $\tau$ were 
$\mathcal{A}$-lower semicontinuous).
\end{ex}

\subsection{Causal curves}
At this point we introduce \emph{timelike}, \emph{causal} and \emph{null curves}, which will be defined via the corresponding relations. 
One has to note that even in the case of a smooth spacetime the class of timelike or causal curves obtained in this way 
differs in general from the class of timelike or causal curves defined in the usual way (via the causal character of the tangent 
vector). If one assumes additionally that a smooth spacetime is strongly causal then the classes of causal curves are the 
same. On the other hand, the classes of timelike curves are still different in general. This will be discussed in more detail 
in Example \ref{ex-rel-cau-cur} below. It will, however, not be an issue since we are mainly interested in the length of 
curves and strong causality ensures that the lengths are unchanged, see Proposition \ref{prop-smo-stc-len}. A causal
curve is non-constant by definition but it could be constant on some interval, contrary to causal curves in a spacetime. 
However, even such curves can be parametrized with respect to $d$-arclength, see \cite[Prop.\ 1.2.2, Cor.\ 1.2.6]{Pap:14}.

\begin{defi}\label{def-cau-cur}
Let $I\subseteq \R$ be an interval. A non-constant curve $\gamma\colon I\rightarrow X$ is called 
\emph{future-directed causal (timelike)} if $\gamma$ is locally Lipschitz continuous (with respect to $d$) and for all 
$t_1,t_2\in I$, $t_1<t_2$ we have that $\gamma(t_1)\leq\gamma(t_2)$ ($\gamma(t_1)\ll\gamma(t_2)$). Furthermore, a future 
directed causal curve is called \emph{null} if no two points on the curve are related with respect to $\ll$. Analogous 
notions are introduced for \emph{past-directed} curves.
\end{defi}
Note that if $d$ is the discrete metric on $X$, then there are no causal curves, since in this case $(X,d)$ is totally disconnected 
and any continuous curve is constant.

\begin{lem}
 Let $\gamma\colon [a,b]\rightarrow X$ be a causal curve, then $\gamma$ is Lipschitz continuous and 
$d$-rectifiable.
\end{lem}
\begin{pr}
Since the domain of $\gamma$ is compact, local Lipschitz continuity implies Lipschitz continuity, which in turn 
implies $d$-rectifiability.
\end{pr}

Let us now investigate the relationship of the different notions of causal and timelike curves for the case of 
continuous or smooth spacetimes.

\begin{ex}\label{ex-rel-cau-cur}
Let $(M,g)$ be the Lorentz cylinder (cf.\ Example \ref{ex-cau-spa-con}), then since $\ll\ =\ \leq\ = M\times M$, every 
non-constant, locally Lipschitz continuous curve is timelike and causal in the sense of Definition 
\ref{def-cau-cur}. Consequently, there are no null curves (again in the sense of Definition \ref{def-cau-cur}). Thus this 
class of timelike and causal curves is larger than the usual class of timelike and causal curves. We will clarify the precise 
relationship in Lemma \ref{lem-rel-cau-cur} below.
\end{ex}

In the following result (and thereafter), when comparing the notions of causal curves in the 
present setting with the standard definition in spacetimes, it will always be understood 
that parametrizations are chosen in which the respective curves are never locally constant
(cf.\ \cite[Ex.\ 2.5.3]{BBI:01}).
\begin{lem}\label{lem-rel-cau-cur}
Let $(M,g)$ be a spacetime with a continuous metric. For clarity we call timelike and causal curves in the sense of 
Definition \ref{def-cau-cur} \emph{R-timelike} and \emph{R-causal}, respectively.
\begin{enumerate}[label=(\roman*)]
 \item \label{lem-rel-cau-cur-cR} A causal/timelike curve is an R-causal/R-timelike curve.
 \item \label{lem-rel-cau-cur-strc-Rc} If $(M,g)$ is a smooth, strongly causal (\cite[Def.\ 14.11]{ONe:83}) 
spacetime, then the notions of 
causal, R-causal and continuous causal curves (\cite[p.\ 192f.]{Wal:84}) coincide.
 \item \label{lem-rel-cau-cur-tim} If $(M,g)$ is a smooth, strongly causal spacetime, then a locally Lipschitz continuous 
curve is an R-timelike curve if and only if it is a continuous timelike curve (\cite[p.\ 192f.]{Wal:84}). However, in general 
the tangent vector of such a curve is only causal almost everywhere, as can be seen from Example \ref{ex-rel-cau-cur-nul}.
\end{enumerate}
\end{lem}
\begin{pr}
Let $\gamma\colon I\rightarrow M$ be a locally Lipschitz continuous curve.
 \begin{enumerate}[label=(\roman*)]
  \item Without loss of generality let $\gamma$ be a future-directed causal/timelike curve. Let $t_1,t_2\in I$ 
with $t_1<t_2$, then $\gamma$ is a future-directed causal/time\-like curve from $\gamma(t_1)$ to $\gamma(t_2)$, hence $\gamma$ 
is R-causal/R-timelike.

 \item Without loss of generality let $\gamma$ be a future-directed R-causal curve and let $t_0\in I$. Let $U$ be a convex 
neighborhood of $\gamma(t_0)$ and let $V\subseteq U$ be a neighborhood of $\gamma(t_0)$ such that every causal curve with 
endpoints in $V$ is contained in $U$. Let $t_1,t_2\in I$ such that $\gamma(t_1),\gamma(t_2)\in V$, hence 
$\gamma([t_1,t_2])\subseteq U$ and consequently $\gamma(t_2)\in J^+(\gamma(t_1),U)$, where $J^+(\gamma(t_1), U)$ denotes the 
set of all points which can be reached from $\gamma(t_1)$ by a Lipschitz continuous future-directed causal curve contained 
in $U$. This is the same set as constructed via piecewise smooth causal curves
(cf.\ \cite[Prop.\ 2.3]{Ler:72}, \cite[Fig.\ 44]{Pen:72b}, and \cite[Cor.\ 2.4.11]{Chr:11}
for the Lipschitz case). This establishes that 
$\gamma$ is a continuous causal curve. Moreover, by \cite[Lemma 14.2(1)]{ONe:83}, 
$\overrightarrow{\gamma(t_1)\gamma(t_2)}\,=\,\mathrm{exp}_{\gamma(t_1)}^{-1}(\gamma(t_2))$ is causal.
Consequently, at any $t_0$ where $\gamma$ is differentiable, 
\[
\dot\gamma(t_0) = \lim_{h\to 0}\frac{1}{h} \overrightarrow{\gamma(t_0)\gamma(t_0+h)}
\]
is causal, so $\gamma$ is a future-directed causal curve.

\item That any R-timelike $\gamma$ is a future-directed continuous timelike curve again follows as in point
\ref{lem-rel-cau-cur-cR} by applying \cite[Cor.\ 2.4.11]{Chr:11}. The converse is clear.
\end{enumerate}
\end{pr}

In fact, as we shall see in Proposition \ref{prop-strongly-causal} below,
the conclusion of Lemma \ref{lem-rel-cau-cur} \ref{lem-rel-cau-cur-strc-Rc} can be extended to continuous spacetimes.

\begin{ex}\label{ex-rel-cau-cur-nul}
Denote by $\R^n_1$ the $n$-dimensional Minkowski spacetime and define $\gamma\colon\R\rightarrow\R^3_1$ by 
$\gamma(t)=(t,\cos(t),\sin(t))$ for $t\in\R$. Then $\gamma$ is a future-directed null curve, but every two points on the 
curve can be joined by a timelike curve, given by the straight line $s\mapsto \gamma(t_1)+ s v$ $(s\in[0,1])$, where 
$t_1<t_2$ and $v:=\gamma(t_2)-\gamma(t_1)$ is future-directed timelike. This shows that $\gamma$ is timelike in the sense of 
Definition \ref{def-cau-cur} and hence by Lemma \ref{lem-rel-cau-cur} \ref{lem-rel-cau-cur-tim} it is also a continuous 
timelike curve. This may be viewed as a caveat concerning the notion of continuous timelike curves as introduced in 
\cite[p.\ 192f.]{Wal:84}.
\end{ex}

\begin{rem}\label{rem-isochronal} The above considerations exemplify the fact that our notion of timelike curves corresponds to maps from intervals in $\R$ into $X$ that preserve the chronal order, hence could be called \emph{isochronal}, while null curves correspond to \emph{achronal} curves. However, we feel that the danger of confusion is rather low and have therefore opted for the above, more intuitive, definitions.
\end{rem}

We now introduce the \emph{length} of a causal curve via the time separation function $\tau$.

\begin{defi}\label{def-cc-len}
 Let $\gamma\colon[a,b]\rightarrow X$ be a future-directed causal curve, then we define its \emph{$\tau$-length} by
\begin{equation}
 L_\tau(\gamma):=\inf\{\sum_{i=0}^{N-1} \tau(\gamma(t_i),\gamma(t_{i+1})): N\in \N,\  a=t_0<t_1<\ldots<t_N=b\}\,,
\end{equation}
and analogously if $\gamma$ is a past-directed causal curve.
\end{defi}

\begin{lem}\label{lem-L-tau-add}
 Let $\gamma\colon[a,b]\rightarrow X$ be a future-directed causal curve and $c\in(a,b)$. Then
\begin{equation}
 L_\tau(\gamma) = L_\tau(\gamma\rvert_{[a,c]}) + L_\tau(\gamma\rvert_{[c,b]})\,.
\end{equation}
\end{lem}
\begin{pr}
 A partition of $[a,c]$ and a partition of $[c,b]$ give naturally a partition of $[a,b]$, hence $L_\tau(\gamma) \leq  
L_\tau(\gamma\rvert_{[a,c]}) + L_\tau(\gamma\rvert_{[c,b]})$. On the other hand, given a partition $a=t_0<t_1<\ldots <t_N=b$ 
of $[a,b]$, we have to consider two cases. The first case is when there is a $k\in\{1,\ldots,N\}$ such that $t_k = c$. 
Consequently, $(t_i)_{i=0}^k$ is a partition of $[a,c]$, and $(t_i)_{i=k}^N$ is a partition of $[c,b]$. Thus
\begin{equation}
 L_\tau(\gamma\rvert_{[a,c]}) + L_\tau(\gamma\rvert_{[c,b]})\leq \sum_{i=0}^{N-1} \tau(\gamma(t_i),\gamma(t_{i+1}))\,.
\end{equation}
If there is no such $k$, then define $j:=\max\{1\leq i \leq N: t_i< c\}$. Then $(t_i)_{i=0}^j \cup \{c\}$ is a partition of 
$[a,c]$ and $\{c\}\cup (t_i)_{i=j}^N$ is a partition of $[c,b]$. Hence
\begin{align*}
 L_\tau(\gamma\rvert_{[a,c]}) + L_\tau(\gamma\rvert_{[c,b]}) &\leq 
\sum_{i=0}^{j-1}\Big(\tau(\gamma(t_i),\gamma(t_{i+1}))\Big) + 
\tau(\gamma(t_j),\gamma(c))\\
&+ \tau(\gamma(c),\gamma(t_{j+1}))+ \sum_{i=j+1}^{N-1} \Big(\tau(\gamma(t_i),\gamma(t_{i+1}))\Big)\\
 &\leq\sum_{i=0}^{N-1} \tau(\gamma(t_i),\gamma(t_{i+1}))\,,
\end{align*}
where in the last inequality we used the reverse triangle inequality. Taking now the infimum over all partitions of $[a,b]$, 
we obtain $L_\tau(\gamma\rvert_{[a,c]}) + L_\tau(\gamma\rvert_{[c,b]})\leq 
L_\tau(\gamma)$.
\end{pr}

\begin{defi}
 By a \emph{reparametrization} of a future-directed causal curve $\gamma\colon[a,b]\rightarrow X$ we mean a 
future-directed causal curve $\lambda\colon[c,d]\rightarrow X$ with $\gamma = \lambda \circ \phi$, where 
$\phi\colon[a,b]\rightarrow [c,d]$ is continuous and strictly monotonically increasing.
\end{defi}
Note that implicit in this definition is the assumption that $\lambda\circ\phi$ is Lipschitz continuous and observe that 
the inverse of such a $\phi$ is also strictly monotonically increasing and continuous.

\begin{lem}
A reparametrization does not change the causality, i.e., the causal character is the same (timelike, null, causal).
\end{lem}
\begin{pr}
 Let $\gamma = \lambda \circ \phi\colon[a,b]\rightarrow X$ be a (without loss of generality) future-directed causal 
(timelike or null) curve and its reparametrization given via $\phi\colon[a,b]\rightarrow [c,d]$. Then $\lambda$ is causal 
(timelike or null). To see this, let $c\leq s_1 < s_2 \leq d$. Since $\phi^{-1}$ is strictly monotonically increasing we have 
that $t_1:= \phi^{-1}(s_1) < \phi^{-1}(s_2)=:t_2$ and thus $\lambda(s_1) = \gamma(t_1)\leq \gamma(t_2) = \lambda(s_2)$ 
($\lambda(s_1) \ll \lambda(s_2)$ or $\lambda(s_1) \nll \lambda(s_2)$).
\end{pr}

\begin{lem}\label{lem-L-tau-inv}
 The $\tau$-length is reparametrization invariant.
\end{lem}
\begin{pr}
 Let $\gamma\colon[a,b]\rightarrow X$ be a future-directed causal curve and $\lambda\colon[c,d]\rightarrow X$ a 
reparametrization of $\gamma$ given via $\phi\colon[a,b]\rightarrow [c,d]$ (i.e., $\gamma = \lambda \circ \phi$). Let 
$a=t_0<t_1<\ldots <t_N = b$ be a partition of $[a,b]$. This yields a partition $c = t_0'< t_1'<\ldots <t_N' = d$ via 
$t_i':=\phi(t_i)$. Consequently, we have
\begin{align*}
 L_\tau(\lambda) &\leq \sum_{i=0}^{N-1}\tau(\lambda(t_i'),\lambda(t_{i+1}')) = 
\sum_{i=0}^{N-1}\tau(\lambda(\phi(t_i)),\lambda(\phi(t_{i+1})))\\
 &= \sum_{i=0}^{N-1}\tau(\gamma(t_i),\gamma(t_{i+1}))\,.
\end{align*}
Now taking the infimum over all partitions of $[a,b]$ yields $L_\tau(\lambda)\leq L_\tau(\gamma)$ and thus by symmetry also 
$L_\tau(\gamma)\leq L_\tau(\lambda)$.
\end{pr}

\begin{defi}
 A future-directed causal curve $\gamma\colon[a,b]\rightarrow X$ is called \emph{rectifiable} if 
$L_\tau(\gamma\rvert_{[t_1,t_2]})>0$ for all $a\leq t_1 < t_2\leq b$.
\end{defi}

\begin{lem}\label{lem-rect-is-timelike}
 A rectifiable causal curve is timelike.
\end{lem}
\begin{pr}
 Let $\gamma\colon[a,b]\rightarrow X$ be a future-directed rectifiable causal curve. Then for $a\leq t_1<t_2\leq b$ we have 
$0<L_\tau(\gamma\rvert_{[t_1,t_2]}) \leq \tau(\gamma(t_1),\gamma(t_2))$. Thus $\gamma(t_1)\ll\gamma(t_2)$ and $\gamma$ is 
timelike.
\end{pr}

\begin{ex}(A timelike curve with $\tau$-length zero, hence non-rectifiable)
Let $\gamma\colon\R\rightarrow\R^3_1$, $\gamma(t)=(t,\cos(t),\sin(t))$ ($t\in\R$) be the timelike curve given in 
Example \ref{ex-rel-cau-cur-nul}. Let $t_1<t_2$ and let $t_1=s_0< s_1 < \ldots < s_k = t_2$ be a partition of $[t_1,t_2]$. 
Then $L_\tau(\gamma\rvert_{[t_1,t_2]}) \leq \sum_{i=0}^{k-1} \tau(\gamma(s_i),\gamma(s_{i+1})) = \sum_{i=0}^{k-1} 
(-\eta(v_i,v_i))$, where $v_i:={\gamma(s_{i+1})-\gamma(s_i)}$ and $\eta$ is the Minkowski metric. It is not hard to 
see that $0<-\eta(v_i,v_i) = (s_{i+1}-s_i)^2 - 2 (1-\cos(s_{i+1}-s_i)) \le \frac{(s_{i+1}-s_i)^4}{12} + 
\frac{2(s_{i+1}-s_i)^6}{6!}\to 0$, for $s_{i+1}-s_i\to 0$. Consequently, for $\eps>0$ one can choose a partition of 
$[t_1,t_2]$ with mesh-size $\delta$ sufficiently small such that $\sum_{i=0}^{k-1}(-\eta(v_i,v_i))\leq k (\frac{\delta^4}{12} 
+ \frac{2\delta^6}{6!}) <\eps$, which shows that $L_\tau(\gamma)=0$. This is a direct proof for this specific curve of the 
general fact that for smooth and strongly causal spacetimes the $\tau$-length agrees with the length given by the 
norm, which we establish in the following proposition.
\end{ex}

\begin{prop}\label{prop-smo-stc-len}
 Let $(M,d^h,\ll,\leq,\tau)$ be the \LpLS induced by the smooth spacetime $(M,g)$, see Example \ref{ex-lpls-smo}. Moreover, 
denote by $L_g(\gamma)=\int_a^b{\sqrt{-g(\dot\gamma,\dot\gamma)}}\,dt$ the length of a causal curve $\gamma\colon[a,b]\rightarrow 
M$ with respect to $g$. If $(M,g)$ is strongly causal, then $L_\tau(\gamma)=L_g(\gamma)$. 
\end{prop}
\begin{pr}
By Lemma \ref{lem-rel-cau-cur} \ref{lem-rel-cau-cur-cR} we know that the notions of causal curves coincide. Let 
$\gamma\colon[a,b]\rightarrow M$ be a future-directed causal curve. Let $a=t_0<t_1<\ldots <t_N=b$ be a partition of $[a,b]$, 
then 
\begin{equation*}
 L_g(\gamma) = \sum_{i=0}^{N-1} L_g(\gamma\rvert_{[t_i,t_{i+1}]}) \leq \sum_{i=0}^{N-1}\tau(\gamma(t_i),\gamma(t_{i+1}))\,. 
\end{equation*}
Taking now the infimum over all partitions of $[a,b]$ we obtain $L_g(\gamma)\leq L_\tau(\gamma)$.

Let $t\in[a,b]$ such that $\dot\gamma(t)$ exists and is future-directed causal and set $e:=\mathrm{exp}_{\gamma(t)}\colon 
\tilde{U}\overset{\cong}{\longrightarrow} U$, where $U$ is a convex neighborhood 
of $\gamma(t)$ such that $e$ is a diffeomorphism from $\tilde{U}:=e^{-1}(U)$ onto $U$. Let $V\subseteq U$ be a neighborhood 
of $\gamma(t)$ such that every causal curve with endpoints in $V$ is contained in $U$ and let $\delta>0$ be such that 
$\gamma([t,t+\delta])\subseteq V$. Then we obtain for $0<h<\delta$
\begin{equation*}
 \frac{1}{h} \tau(\gamma(t),\gamma(t+h)) = \|\frac{1}{h} e^{-1}(\gamma(t+h))\|\,,
\end{equation*}
where $\|.\| = \sqrt{|g(.,.)|}$. Thus, taking the limit $h\searrow 0$ we get
\begin{equation}\label{eq-smo-len-exp}
 \lim_{h\searrow 0} \frac{\tau(\gamma(t),\gamma(t+h))}{h} = \Big \| \frac{d}{d h}\Big|_0\, 
e^{-1}(\gamma(t+h)) \Big\| = \| (T_0 \mathrm{exp}_{\gamma(t)})^{-1}(\dot\gamma(t))\| = \|\dot\gamma(t)\|\,,
\end{equation}
where we used that $T_0 \exp_{\gamma(t)} = \mathrm{id}$ and that $e^{-1}\circ\gamma$ is differentiable at $t$, and so the 
one-sided derivative agrees with the derivative. Furthermore,
\begin{equation*}
 \frac{\tau(\gamma(t),\gamma(t+h))}{h} \geq  \frac{1}{h}L_\tau(\gamma\rvert_{[t,t+h]}) \geq  
\frac{1}{h} L_g(\gamma\rvert_{[t,t+h]}) = \frac{1}{h}\int_t^{t+h}\|\dot\gamma(s)\| \mathrm{d}s\,,
\end{equation*}
where we used the above ($L_\tau\geq L$). Now the left hand side goes to $\|\dot\gamma(t)\|$ as $h\searrow 0$ by 
\eqref{eq-smo-len-exp} and obviously so does the right hand side. Consequently,  
$\lim_{h\searrow 0}\frac{1}{h}L_\tau(\gamma\rvert_{[t,t+h]})=\|\dot\gamma(t)\|$ as well. Finally, we obtain for any segment 
of $\gamma$ that is contained in such a convex neighborhood $U$, say $\gamma([t_0,t_1])\subseteq U$ and almost all $t$ in 
$[t_0,t_1]$ that 
\begin{equation}\label{eq-smo-len-der}
 D_+(t\mapsto L_\tau(\gamma\rvert_{[t_0,t]})) = \|\dot\gamma(t)\|\,,
\end{equation}
where $D_+$ denotes the right-sided derivative.

We now establish that $\phi\colon[t_0,t_0+\delta]\rightarrow [0,\infty)$ given by $\phi(t):= L_\tau(\gamma\rvert_{[t_0,t]})$ 
is absolutely continuous. 
Let $([a_i,b_i])_{i=1}^N$ 
be a collection of non-overlapping intervals in $[t_0,t_0+\delta]$ with $\sum_{i=1}^N (b_i-a_i) < \alpha$, $\alpha$ to be 
given later. Then we calculate
\begin{equation*}
 \sum_{i=1}^N |\phi(b_i)-\phi(a_i)| = \sum_{i=1}^N L_\tau(\gamma\rvert_{[a_i,b_i]})\leq \sum_{i=1}^N 
\tau(\gamma(a_i),\gamma(b_i))=:(*)\,,
\end{equation*}
where we used that $\phi$ is monotonically increasing (cf.\ Lemma \ref{lem-phi}) and that $L_\tau$ is additive by Lemma 
\ref{lem-L-tau-add}. In the convex neighborhood $U$ we know that for $p< q$ with $p,q \in V$, the maximal causal curve 
joining $p$ and $q$ is contained in $U$ and its length is given by $\|\Delta(p,q)\|$. Here $\Delta:=E^{-1}\colon U\times U 
\rightarrow TM$ is a diffeomorphism onto its image and $E(v) = (\pi(v),\mathrm{exp}(v))$ for $v\in TM$ in the
domain of $E$, cf.\ \cite[Lemma 5.9]{ONe:83}. This implies 
that 
\begin{equation*}
 (*)\leq \sum_{i=1}^N \|\Delta(\gamma(a_i),\gamma(b_i))\|\leq C \sum_{i=1}^N \|\Delta(\gamma(a_i),\gamma(b_i))\|_2=:(**)\,,
\end{equation*}
for some constant $C$ (depending only on $g$ and $U$), where $\|.\|_2$ denotes the Euclidean norm in these 
coordinates. Since $\Delta(p,.)$ is smooth, it is locally Lipschitz continuous, and since $\Delta(p,p)=0$ for all $p\in 
U$ we get $\|\Delta(\gamma(a_i),\gamma(b_i))\|_2 = \|\Delta(\gamma(a_i),\gamma(b_i)) - 
\Delta(\gamma(a_i),\gamma(a_i))\|_2\leq C' \|\gamma(a_i)-\gamma(b_i)\|_2$ $\leq \tilde{C}(b_i-a_i)$. In the last inequality we 
used the Lipschitz continuity of $\gamma$. Finally, we get for $\eps>0$ and $0< \alpha<\frac{\eps}{C \tilde{C}}$ that
\begin{equation*}
 (**)\leq C \tilde{C} \sum_{i=1}^N (b_i-a_i) < C \tilde{C} \alpha < \eps\,,
\end{equation*}
establishing the absolute continuity of $\phi$.  It follows that there exists a subset of full measure in $[t_0,t_0+\delta]$ 
on which $\phi$ is differentiable and where its derivative is given by \eqref{eq-smo-len-der}. 
This enables us to apply the Fundamental Theorem of Calculus to obtain
\begin{equation*}
L_\tau(\gamma\rvert_{[t_0,t_1]}) = \int_{t_0}^{t_1} \frac{d}{d s}L_\tau(\gamma\rvert_{[t_0,s]})\mathrm{d}s = 
\int_{t_0}^{t_1} \|\dot\gamma(s)\|\mathrm{d}s  = L_g(\gamma\rvert_{[t_0,t_1]})\,.
\end{equation*}
Covering $\gamma([a,b])$ with finitely many such neighborhoods $V$ and using the additivity of both $L_\tau$ (Lemma 
\ref{lem-L-tau-add}) and $L$ yields $L_\tau(\gamma) = L_g(\gamma)$.
\end{pr}

\subsection{Maximal causal curves}

\begin{defi}
 Let \Xll be a Lorentzian pre-length space. A future-directed causal curve $\gamma\colon [a,b]\rightarrow X$ is 
\emph{maximal} if $L_\tau(\gamma) = \tau(\gamma(a),\gamma(b))$, and analogously for past-directed causal curves.
\end{defi}
A note on terminology is in order here: according to the above definition, a maximal curve 
$\gamma$ is a time-separation realizing curve. Any such $\gamma$ is also maximal in the following 
sense: let $\sigma$ be another causal curve connecting $p=\gamma(a)$ and $q=\gamma(b)$ with
$L_\tau(\sigma)\ge L_\tau(\gamma)$. Then $L_\tau(\sigma) = L_\tau(\gamma)$: in fact,
by the definition of $\tau$-length we have $L_\tau(\sigma)\le \tau(p,q)=L_\tau(\gamma)$.
For Lorentzian length spaces (see Section \ref{sec-lorentzian-length-spaces}) both notions of maximality 
in fact coincide.

\begin{prop}\label{prop-max-prop}
 Let \Xll be a Lorentzian pre-length space.
 \begin{enumerate}[label=(\roman*)]
  \item A null curve is always maximal on any compact interval.
  \item A maximal curve is maximal on any subinterval. \label{prop-max-sub-int}
  \item If a maximal curve is timelike then it is rectifiable. \label{prop-max-tl-rec}
 \end{enumerate}
\end{prop}
\begin{pr}
 Let $\gamma\colon[a,b]\rightarrow X$ be a future-directed causal curve.
 \begin{enumerate}[label=(\roman*)]
  \item Let $\gamma$ be null. Then for all $a\leq t_1 < t_2 \leq b$ we have $\tau(\gamma(t_1),\gamma(t_2))= 0 $, which 
implies $L_\tau(\gamma) = 0$. Thus $L_\tau(\gamma) = 0 = \tau(\gamma(a),\gamma(b))$ and $\gamma$ is maximal.
  \item Let $\gamma$ be maximal and $a\leq c < d \leq b$ a subinterval. Assume that $\gamma$ is not maximal on $[c,d]$, 
i.e., $L_\tau(\gamma\rvert_{[c,d]}) < \tau(\gamma(c),\gamma(d))$. Then by Lemma \ref{lem-L-tau-add} and the reverse triangle 
inequality we get
\begin{align*}
\tau(\gamma(a),\gamma(b)) &= L_\tau(\gamma) = L_\tau(\gamma\rvert_{[a,c]}) + L_\tau(\gamma\rvert_{[c,d]}) + 
L_\tau(\gamma\rvert_{[d,b]})\\
&< L_\tau(\gamma\rvert_{[a,c]}) +\tau(\gamma(c),\gamma(d)) + L_\tau(\gamma\rvert_{[d,b]})\\
&\leq \tau(\gamma(a),\gamma(c)) + 
\tau(\gamma(c),\gamma(d)) + \tau(\gamma(d),\gamma(b))\\
&\leq \tau(\gamma(a),\gamma(b))\,.
\end{align*}
This is a contradiction, thus establishing that $\gamma$ is maximal on $[c,d]$.
 \item Let $\gamma$ be timelike. Then for all $a\leq t_1 < t_2 \leq b$ we have $0<\tau(\gamma(t_1),\gamma(t_2))= 
L_\tau(\gamma\rvert_{[t_1,t_2]})$. Thus $\gamma$ is rectifiable.
\end{enumerate}
\end{pr}

\subsection{Causality conditions}

\begin{defi}
 A causal space $(X,\ll,\leq)$ is called
\begin{enumerate}[label=(\roman*)]
 \item \emph{chronological} if the relation $\ll$ is irreflexive, i.e., $x\nll x$ for all $x\in X$, and
 \item \emph{causal} if the relation $\leq$ is a partial order, i.e., $x\leq y$ and $y\leq x$ implies that $x=y$ for all 
$x,y\in X$.
\end{enumerate}
A \LpLS \Xll is called

\begin{enumerate}[label=(\roman*)]
\setcounter{enumi}{2}
 \item \emph{non-totally imprisoning} if for every compact set $K\Subset X$ there is a $C>0$ such that the $d$-arclength 
of all causal curves contained in $K$ is bounded by $C$,
 \item \emph{strongly causal} if the Alexandrov topology $\mathcal{A}$ agrees with the metric topology $\mathcal{D}$ 
(and hence also with the chronological topology $\mathcal{I}$), and
 \item \emph{globally hyperbolic} if \Xll is non-totally imprisoning and for every $x,y\in X$ the set $J^+(x)\cap J^-(y)$ is 
compact in $X$.
\end{enumerate}
\end{defi}

\begin{rem}
 Causality does not imply chronology in general as can be seen from the simple example: $X:=\{*\}$, $\ll:=\leq := 
\{(*,*)\}$. Clearly, both $\ll$ and $\leq$ are transitive and reflexive, hence $(X,\ll,\leq)$ is not chronological but it is 
causal.
\end{rem}

\begin{defi}
A causal space $(X,\ll,\leq)$ is called \emph{interpolative} if for all $x,y\in X$ with $x\ll y$ there is a 
$z\in X$ such that $x\ll z\ll y$ and $x\neq z\neq y$.
\end{defi}

\begin{lem}\label{lem-cc}
 Let \Xll be a Lorentzian pre-length space. Then
\begin{enumerate}[label=(\roman*)]
 \item if \Xll is causal and interpolative it is chronological,\label{lem-cc-cai-ch}
 \item if \Xll is chronological then the time separation function $\tau$ is zero on the diagonal, i.e., $\tau(x,x)=0$ for 
all $x\in X$, and\label{lem-cc-ch-tau}
 \item if \Xll is strongly causal, then for all $x\in X$, for every neighborhood $U$ of $x$, there is a neighborhood 
$V\subseteq U$ of $x$ such that for every causal curve $\gamma\colon[a,b]\rightarrow X$ with $\gamma(a),\gamma(b)\in V$ one 
has $\gamma([a,b])\subseteq U$ (i.e., the usual definition of strong causality for spacetimes).\label{lem-cc-str-cau}
\end{enumerate}
\end{lem}
\begin{pr}
 \begin{enumerate}[label=(\roman*)]
  \item Assume that there is an $x\in X$ such that $x\ll x$. Then since \Xll is interpolative there is a $z\in X$ such that 
$x\ll z \ll x$ and $x\neq z$. This implies $x\leq z \leq x$ and since $\leq$ is a partial order $x=z$ --- a contradiction.
 \item This follows from $\tau(x,x)>0 \Leftrightarrow x\ll x$ for all $x\in X$.
 \item Let $x\in X$ and let $U$ be a neighborhood of $x$. Then since $\mathcal{A}=\mathcal{D}$ there is a 
$V\in \mathcal{A}$ such that $x\in V$ and $V\subseteq U$. We may assume that $V=(I^+(x_1)\cap I^-(y_1)) \cap 
\ldots \cap (I^+(x_n)\cap I^-(y_n))$ for some $x_1,y_1,\ldots,x_n,y_n\in X$. Now let $\gamma\colon[a,b]\rightarrow X$ be a 
(without loss of generality) future-directed causal curve with $\gamma(a),\gamma(b)\in V$. We claim that $\Gamma\subseteq V$, 
thus $\gamma([a,b])\subseteq U$. Let $t\in[a,b]$, then $\gamma(t)\in J^+(\gamma(a))$. For each $1\leq i \leq n$ we have $x_i\ll 
\gamma(a)\leq \gamma(t)$ and so $x_i\ll\gamma(t)$ by push-up (Lemma \ref{lem-pup}). This is equivalent to $\gamma(t)\in 
I^+(x_i)$. Analogously, one shows that $\gamma(t)\in I^-(y_i)$.
 \end{enumerate}
\end{pr}

It is not clear at this moment if strong causality is equivalent to the condition of Lemma 
\ref{lem-cc} \ref{lem-cc-str-cau}, i.e., the non-existence of almost closed causal curves, as is the case for smooth 
spacetimes, see \cite[Thm.\ 3.27]{MS:08} without further assumptions or more structure on the \LpLSn. However, 
for \LLSn s the crucial additional ingredient will be \emph{localizability} and for these spaces the conditions will be 
equivalent, see Theorem \ref{thm-cau-lad-lls} \ref{thm-cau-lad-lls-str-cau-no-acc}.

\section{Lorentzian length spaces}\label{sec-lorentzian-length-spaces}

\subsection{Causal connectedness}

\begin{defi}\label{def-cau-pc}
 A \LpLS \Xll is called \emph{causal\-ly path connected} if for all $x,y\in X$ with $x\ll y$ there is a future-directed 
timelike curve from $x$ to $y$ and for $x<y$ there is a future-directed causal curve from $x$ to $y$.
\end{defi} 

\begin{lem}\label{lem-cpc-i}
 A causally path connected \LpLS is interpolative.
\end{lem}
\begin{pr}
 Let $x,y\in X$ with $x\ll y$, then there is a future-directed timelike curve $\gamma\colon[a,b]\rightarrow X$ from $x$ to 
$y$. Since $\gamma$ is not constant, there is a $t\in [a,b]$ with $x=\gamma(a)\neq \gamma(t)=:z$, and because $\gamma$ is 
timelike we have $x\ll z$. If $x=y$ we are done, and if $x\neq y$ there is a $\delta>0$ such that $\bar{B}^d_\delta(x)\cap 
\bar{B}^d_\delta(y)=\emptyset$. Then if $\gamma([a,b])\subseteq \bar{B}^d_\delta(x)\cup \bar{B}^d_\delta(y)$ it would
follow that $\gamma([a,b])$ can be written as a disjoint union of the non-empty closed sets $\bar{B}^d_\delta(x)\cap\gamma([a,b])$ and 
$\bar{B}^d_\delta(y)\cap\gamma([a,b])$, contradicting connectedness. Thus there is a $t'\in[a,b]$ such that $x\neq 
z=\gamma(t') \neq y$ and by assumption $x\ll z \ll y$.
\end{pr}

\begin{lem}\label{lem-cpc-cc}
 Let \Xll be a causally path connected \LpLSn.
\begin{enumerate}[label=(\roman*)]
 \item \Xll is chronological if and only if there are no closed timelike curves in $X$.\label{lem-cpc-cc-chr}
 \item \Xll is causal if and only if there are no closed causal curves in $X$.\label{lem-cpc-cc-cau}
\end{enumerate}
\end{lem}
\begin{pr}
Let \Xll be a causally path connected \LpLSn, which for brevity we just denote by $X$.
\begin{enumerate}[label=(\roman*)]
\item
\begin{itemize}
\item [($\Rightarrow$):] Let $X$ be chronological and $\gamma$ a closed timelike curve. Then for all $x$ in the 
image of $\gamma$ we 
have $x\ll x$ --- a contradiction.
\item [($\Leftarrow$):] Let $X$ be such that there are no closed timelike curves. Let $x\in X$ with $x\ll x$, then by 
the causal path-connectedness there is a future-directed timelike curve from $x$ to $x$ --- a contradiction.
\end{itemize}
\item
\begin{itemize}
\item [($\Rightarrow$):] Let $X$ be causal and $\gamma$ a closed causal curve. Since $\gamma$ is not constant, there are 
points $x,y$ on $\gamma$ with $x\neq y$ and by assumption $x\leq y\leq x$ --- a contradiction.
\item [($\Leftarrow$):] Let $X$ be such that there are no closed causal curves. Let $x,y\in X$ with $x< y< x$, then by 
the causal path-connectedness there is a future-directed causal curve from $x$ to $y$ to $x$ --- a contradiction.
\end{itemize}
\end{enumerate}
\end{pr}

\subsection{Limit curves}

\begin{defi}
 Let \Xll be a \LpLS and let $x\in X$. A neighborhood $U$ of $x$ is called \emph{causally closed} if 
$\leq$ is closed in $\bar{U}\times\bar{U}$, i.e.,  if $p_n,q_n\in U$ 
with $p_n\leq q_n$ for all $n\in\N$ and $p_n\to p\in\bar{U}$, $q_n\to q\in\bar{U}$, then $p\leq q$. A \LpLS \Xll is called 
\emph{locally causally closed} if every point has a causally closed neighborhood.
\end{defi}

\begin{prop}\label{prop-cau-loc-cc}
 Strongly causal spacetimes with continuous metrics are locally causally closed.
\end{prop}
\begin{pr}
 Let $(M,g)$ be a strongly causal spacetime with $g$ continuous. Then $(M,g)$ is non-totally imprisoning, cf. 
\cite[p.\ 1437]{Sae:16}. Let $p\in M$ and $U$ an open, relatively compact neighborhood of $p$, then by \cite[Lemma 2.7]{Sae:16} 
there is a $C>0$ such that $L^h(\gamma)\leq C$ for all causal curves $\gamma$ with image contained in $U$, where 
$h$ is a complete Riemannian metric on $M$. Strong causality implies the existence of a neighborhood $V$ of $p$, 
$V\subseteq U$ such that for all causal curves $\lambda\colon[a,b]\rightarrow M$ with $\lambda(a),\lambda(b)\in V$ one has 
$\lambda([a,b])\subseteq U$. Now let $(x_n)_n$, $(y_n)_n$ be sequences in $V$ with $x_n\leq y_n$ for all $n\in\N$ and 
$x_n\to x\in \bar{V}$, $y_n\to y\in\bar{V}$. Thus there is a sequence $(\gamma_n)_n$ of future-directed causal curves 
$\gamma_n\colon[0,1]\rightarrow M$ with $\gamma_n(0)=x_n\in V$, $\gamma_n(1)=y_n\in V$ for all $n\in\N$. Hence 
$\gamma_n([0,1])\subseteq U$ and so $L^h(\gamma_n)\leq C$ for all $n\in\N$. Finally, the limit curve theorem 
(\cite[Thm.\ 1.5]{Sae:16}) establishes the existence of a future-directed causal curve from $x$ to $y$, if $x\not=y$, thus 
$x<y$.
\end{pr}

\begin{lem}\label{lem-lim-cc}
Let \Xll be a locally causally closed \LpLS and let $(\gamma_n)_n$ be a sequence of future-directed causal curves 
$\gamma_n\colon[a,b]\rightarrow X$ converging pointwise to a non-constant Lipschitz curve 
$\gamma\colon [a,b]$ $\rightarrow X$. Then $\gamma$ is a future-directed causal curve.
\end{lem}
\begin{pr}
 For every $t\in[a,b]$ there is an open, causally closed neighborhood $U_t$ of $\gamma(t)$. Let $a\leq t_1<t_2\leq b$ such 
that $\gamma(t_1),\gamma(t_2)\in U_t$. Then there is an $n_0\in\N$ such that for all $n\geq n_0$ we have 
$\gamma_n(t_1),\gamma_n(t_2)\in U_t$. Since $\gamma_n(t_1)\leq \gamma_n(t_2)$ for all $n\in\N$ and by assumption 
$\gamma_n(t_i)\to\gamma(t_i)$ for $i=1,2$, we conclude that $\gamma(t_1)\leq\gamma(t_2)$. This gives an open cover of the 
compact set $\gamma([a,b])$, from which we may extract a finite sub-cover $U_1,\ldots,U_N$. Additionally, this 
finite cover has a Lebesgue number $\delta>0$. Let $L>0$ be the Lipschitz constant of $\gamma$, then if 
$|t_1-t_2|\leq \frac{\delta}{L}$, one has $\gamma(t_1),\gamma(t_2)\in U_i$ for some $i\in\{1,\ldots,N\}$. Now let $a\leq 
t_1<t_2\leq b$, and let $t_1 =: s_0 < s_1 < \ldots < s_{k-1} < s_k:= t_2$ such that $|s_{j+1}-s_j|\leq \frac{\delta}{L}$ for 
all $j\in\{0,\ldots,k-1\}$. Thus, by construction $\gamma(s_j),\gamma(s_{j+1})\in U_{l_j}$ for 
all $j\in\{0,\ldots,k-1\}$ and corresponding $l_j\in\{1,\ldots,N\}$ and so $\gamma(t_1) = \gamma(s_0)\leq \gamma(s_1)\leq 
\ldots \leq\gamma(s_k)=\gamma(t_2)$, hence $\gamma(t_1)\leq\gamma(t_2)$.
\end{pr}

\begin{thm}(Limit curve theorem)\label{thm-lim-curve}
 Let \Xll be a locally causally closed \LpLSn. Let $(\gamma_n)_n$ be a sequence of future-directed causal curves 
$\gamma_n\colon[a,b]\rightarrow X$ that are uniformly Lipschitz continuous, i.e., there is an $L>0$ such that 
$\Lip(\gamma_n)\leq L$ for all $n\in\N$. Suppose that there exists a compact set that contains every $\gamma_n([a,b])$ or that
$d$ is proper (i.e., all closed and bounded sets are compact) and that the curves $(\gamma_n)_n$ accumulate at some point, 
i.e., there is a $t_0\in[a,b]$ such that $\gamma_n(t_0)\to x_0\in X$. Then there exists a subsequence $(\gamma_{n_k})_k$ of 
$(\gamma_n)_n$ and a Lipschitz continuous curve $\gamma\colon[a,b]\rightarrow X$ such that $\gamma_{n_k}\to\gamma$ 
uniformly. If $\gamma$ is non-constant, then $\gamma$ is a future-directed causal curve. In particular, if $\gamma_n(a)=p$, 
$\gamma_n(b)=q$ for all $n\in\N$, with $p\neq q$, then $\gamma$ is a future-directed causal curve connecting $p$ and $q$.
\end{thm}
\begin{pr}
 The sequence $(\gamma_n)_n$ is equicontinuous and either the $\gamma_n$s are contained in a compact set by assumption 
or $(X,d)$ is proper. In the latter case we have for $t\in[a,b]$ and all $n\in\N$
\begin{align*}
 d(x_0,\gamma_n(t))&\leq d(x_0,\gamma_n(t_0)) + d(\gamma_n(t_0),\gamma_n(t))\\
&\leq C + \Lip(\gamma_n)|t-t_0| \leq C + L (b-a)\,,
\end{align*}
where $C>0$ is some constant determined by the convergence of $\gamma_n(t_0)$ to $x_0$. Thus in both cases
$(\gamma_n(t))_n$ is relatively compact for all $t\in[a,b]$ and so we can apply the Arzel\`a-Ascoli theorem 
(e.g.\ \cite[Thm.\ 1.4.9]{Pap:14}) to get a uniformly converging subsequence $(\gamma_{n_k})_k$. The uniform limit 
$\gamma:=\lim_{k\to\infty}\gamma_{n_k}$ is Lipschitz continuous, with $\Lip(\gamma)\leq L$ and thus Lemma 
\ref{lem-lim-cc} shows that $\gamma$ is a future-directed causal curve. 
\end{pr}

\begin{lem}(A sufficient condition that the limit curve is not constant, cf.\ \cite[Thm.\ 3.1]{Min:08a})\label{lem-suff-cond-non-const}
 Let $(\gamma_n)_n$ be a sequence of (continuous) curves defined on $[a,b]$ that converge uniformly to a curve 
$\gamma\colon[a,b]\rightarrow X$. If there is a $t\in[a,b]$ and a neighborhood $U$ of $\gamma(t)$ such that only finitely 
many $\gamma_n$ are contained in $U$, then $\gamma$ is not constant.
\end{lem}
\begin{pr}
Let $\eps>0$ be such that $B^d_\eps(\gamma(t))\subseteq U$. The assumptions yield that there is an $n_0\in\N$ such that for 
all $n\geq n_0$ there is an $s_n\in[a,b]$ with $\gamma_n(s_n)\not\in U$, hence $0<\eps \le d(\gamma(t),\gamma_n(s_n))$. 
Without loss of generality we can assume that $s_n\to s^*\in[a,b]$ and thus $0<\eps\leq d(\gamma(t),\gamma(s^*))$.
\end{pr}

\begin{rem}
 We do not require $(X,d)$ to be a proper metric space in the definition of a \LpLSn, since in the case where we will
apply such results (as the Limit curve theorem above) in the development of the theory, it will be to (relatively) compact 
subsets.
\end{rem}

We now introduce inextendible causal curves.

\begin{defi}
 Let $-\infty\leq a < b\leq \infty$ and let $\gamma\colon[a,b)\rightarrow X$ be a future (or past)-directed causal (or 
timelike) curve. It is called \emph{extendible} if there exists a future (past)-directed causal (timelike) curve 
$\tilde{\gamma}\colon[a,b]\rightarrow X$ such that $\tilde{\gamma}\rvert_{[a,b)}=\gamma$. The curve $\gamma$ is called 
\emph{inextendible} if it is not extendible. Analogously for the other endpoint of the interval.
\end{defi}

\begin{rem}
 An extendible causal curve is Lipschitz continuous on its (open) domain of definition.
\end{rem}

\begin{lem}\label{lem-inext-cc-infty}
Let \Xll be a locally causally closed \LpLSn, let $-\infty < a<b\leq\infty$ and let $\gamma\colon[a,b)\rightarrow X$ be a 
(without loss of generality) future-directed causal curve parametrized with respect to $d$-arclength. If $(X,d)$ is a 
proper metric space or the image of $\gamma$ is contained in a compact set, then $\gamma$ is inextendible if and only if 
$b=\infty$. In this case $L^d(\gamma)=\infty$. Moreover, $\gamma$ is inextendible if and only if $\lim_{t\nearrow 
b}\gamma(t)$ does not exist.
\end{lem}
\begin{pr}
We first show the equivalence of $\gamma$ inextendible and $b=\infty$.
\begin{itemize}
  \item[($\Leftarrow$):] Since $\gamma$ is parametrized with respect to $d$-arclength, we have $d(\gamma(a),\gamma(t))$ $= 
t-a$ for all $t\in[a,b)$. Thus if $\gamma$ were extendible we would have $b=d(\gamma(a),\gamma(b)) + a<\infty$ --- a 
contradiction. 
  \item[($\Rightarrow$):] Assume that $b<\infty$. In both cases we have that $\gamma$ is contained in a compact set. 
Either by assumption or if $(X,d)$ is proper, then $\gamma([a,b))\subseteq B^d_{b-a}(\gamma(a))$, which is relatively compact. 
%This can be seen as follows: for $t\in[a,b)$, we have $d(\gamma(a),\gamma(t)) = t-a<b-a$.
Thus there exists a sequence 
$(t_n)_n$ with $t_n\nearrow b$ and $\lim_{n\to\infty}\gamma(t_i)=:p$. 
This is the only limit point of $\gamma$ as the 
parameter approaches $b$. Assume that there is another sequence $(s_n)_n$ such that $s_n\nearrow b$ and 
$\lim_{n\to\infty}\gamma(s_n)=:q\neq p$. Then we have $\delta:= d(p,q)>0$ and thus there is an $n_0\in\N$ such that for all 
$n\geq n_0$ we have $(b-t_n)<\frac{\delta}{4}$, $(b-s_n)<\frac{\delta}{4}$, $d(p,\gamma(t_n))<\frac{\delta}{4}$ and 
$d(q,\gamma(s_n))<\frac{\delta}{4}$. The curve $\gamma$ is $1$-Lipschitz continuous and so we obtain
\begin{align*}
 \delta &= d(p,q) \leq d(p,\gamma(t_n)) + d(\gamma(t_n),\gamma(s_n)) + d(\gamma(s_n),q)\\
        &< \frac{\delta}{4} + |t_n - s_n| + \frac{\delta}{4} \leq \frac{\delta}{2} + (b - t_n) + (b - s_n) < \delta\,,
\end{align*}
a contradiction. At this point we show that we can extend $\gamma$ via $\tilde{\gamma}$ given by 
$\tilde{\gamma}\rvert_{[a,b)}:=\gamma$, $\tilde{\gamma}(b):=p$. Clearly, $\tilde{\gamma}$ is $1$-Lipschitz continuous, and 
so it remains to show that $\tilde{\gamma}(t)\leq p$ for all $t\in[a,b)$. Let $U$ be a causally closed neighborhood of $p$. 
Then there is a $t^*\in[a,b)$ such that for all $t\in(t^*,b)$, $\tilde{\gamma}(t)=\gamma(t)\in U$. Fix $t\in(t^*,b)$ and let 
$(t_n)_n$ be a sequence in $(a,b)$ with $t_n\nearrow b$ and $t\leq t_0$. This yields that $\gamma(t)\leq \gamma(t_n)$ for 
all $n\in\N$ since $\gamma$ is causal. By construction $\gamma(t_n)\to p$ and hence by causal closedness of $U$ we obtain 
$\gamma(t)\leq p$. This shows that $\gamma(t)\leq p$ for all $t\in(t^*,b)$. Now fix $t\in[a,t^*]$ and let $t'\in(t^*,b)$. 
Then $\gamma(t)\leq \gamma(t')\leq p$ and by transitivity $\gamma(t)\leq p$, as required.
 \end{itemize}
The latter implication shows that if $\lim_{t\nearrow b}\gamma(t)$ exists then $\gamma$ is extendible. Conversely, if 
$\gamma$ is extendible the limit obviously exists.
\end{pr}

\begin{defi}
 A \LpLS \Xll is called \emph{$d$-compatible} if for every $x\in X$ there exists a neighborhood $U$ of $x$ and a constant $C>0$ 
such that $L^d(\gamma)\leq C$ for all causal curves $\gamma$ contained in $U$.
\end{defi}

\begin{thm}(Limit curve theorem for inextendible curves)\label{thm-lim-curve-inext}
 Let \Xll be a locally causally closed and $d$-compatible \LpLSn. Let $(\gamma_n)_n$ be a sequence of future-directed causal 
curves $\gamma_n\colon[0,L_n]\rightarrow X$ that are parametrized with respect to $d$-arclength with 
$L_n:=L^d(\gamma_n)\to\infty$. If there exists a compact set that contains every $\gamma_n([0,L_n])$ or if $d$ is proper and 
$\gamma_n(0)\to x$ for some $x\in X$, then there exists a subsequence $(\gamma_{n_k})_k$ of $(\gamma_n)_n$ and a future 
directed causal curve $\gamma\colon[0,\infty)\rightarrow X$ such that $\gamma_{n_k}\to\gamma$ locally uniformly. Moreover, 
$\gamma$ is inextendible.
\end{thm}
\begin{pr}
Extend each $\gamma_n$ constantly to $[0,\infty)$ and denote it again by $\gamma_n$. Then the sequence $(\gamma_n)_n$ is 
equicontinuous since every $\gamma_n$ is $1$-Lipschitz continuous and as in the proof of Theorem \ref{thm-lim-curve} we have
$\gamma_n([0,\infty))\subseteq K$, for some compact set $K\Subset X$. Again, the Arzel\`a-Ascoli theorem 
%(e.g.\ \cite[Theorem 1.4.9, p.\ 33f.]{Pap:14}) 
gives a locally uniformly converging subsequence $(\gamma_{n_k})_k$. The limit curve 
$\gamma:=\lim_{k\to\infty}\gamma_{n_k}\colon[0,\infty)\rightarrow X$ is $1$-Lipschitz continuous. To see that 
$\gamma$ is causal, observe that for every $t>0$ there is a $k_0\in\N$ such that $L_{n_k}\geq t$ for all $k\geq k_0$ and we 
have $\gamma_{n_k}\rvert_{[0,t]}\to \gamma\rvert_{[0,t]}$ uniformly. To apply Lemma \ref{lem-lim-cc} we need to show that 
$\gamma\rvert_{[0,t]}$ is not constant (at least for $t>0$ sufficiently large). Let $U$ be a neighborhood of $\gamma(0)$ 
such that there is a $C>0$ that bounds the $d$-arclength of all causal curves in $U$ and let $t>C$. Then since 
$L_{n_k}\geq t>C$, the $\gamma_{n_k}$s cannot be contained in $U$, hence Lemma \ref{lem-suff-cond-non-const} yields that 
$\gamma\rvert_{[0,t]}$ is not constant. Thus $\gamma\rvert_{[0,t]}$ is future-directed causal by Lemma \ref{lem-lim-cc}, 
since $\gamma_{n_k}\rvert_{[0,t]}$ is causal as a segment of the original causal curve $\gamma_{n_k}$ defined on 
$[0,L_{n_k}]$ for $k\geq k_0$. Now let $0\leq t_1 < t_2<\infty$, then, by the above, there is a $t\geq t_2$ such that
$\gamma\rvert_{[0,t]}$ is causal and hence $\gamma(t_1)\leq\gamma(t_2)$. It remains to show the inextendibility of 
$\gamma$. Assume, to the contrary, that $\gamma$ is extendible and set $\lim_{t\nearrow\infty}\gamma(t)=:p$. Let $V$ be a 
neighborhood of $p$ such that there is a $C>0$ that bounds the $d$-arclength of all causal curves in $V$. There is 
a $t^*\in[0,\infty)$ such that $\gamma([t^*,\infty))\subseteq V$ and hence for $T>t^*$ with $T-t^*>C$ there is a $k\in\N$ 
such that $\gamma_{n_k}([t^*,T])\subseteq V$ by the uniform convergence on $[t^*,T]$ and $L_{n_k}>T$. This is a 
contradiction as $L^d(\gamma_{n_k}\rvert_{[t^*,T]}) = T-t^* >C$.
\end{pr}

\begin{cor}\label{cor-nti-inext}
 Let \Xll be a locally causally closed and $d$-com\-pat\-i\-ble \LpLSn. Then $X$ is non-totally imprisoning if and only if 
no compact set in $X$ contains an inextendible causal curve.
\end{cor}
\begin{pr}

 \item[($\Rightarrow$):] Assume that there is a compact set $K\Subset X$ and an inextendible causal curve $\gamma$ contained 
in $K$. By Lemma \ref{lem-inext-cc-infty} we have $L^d(\gamma)=\infty$ --- a contradiction to $X$ being non-totally 
imprisoning.
 
 \item[($\Leftarrow$):] Assume that there is a compact set $K\Subset X$ and a sequence of (without loss of generality) 
future-directed causal curves $\gamma_n\colon I_n\rightarrow X$ contained in $K$ with $L^d(\gamma_n)\to\infty$. 
Parametrizing them with respect to $d$-arclength gives a sequence $\lambda_n\colon[0,L_n]\rightarrow X$, with 
$L_n:=L^d(\gamma_n)=L^d(\lambda_n)$. Now Theorem \ref{thm-lim-curve-inext} yields a limit curve of this sequence that is an 
inextendible causal curve contained in $K$ --- a contradiction.
\end{pr}

\subsection{Localizability}

We now try to capture the idea that locally the geometry and causality of a (smooth) Lorentzian manifold is better 
behaved than globally. The following definition generalizes to our current setting a number of essential properties 
inherent to convex neighborhoods in smooth Lorentzian manifolds. Also in metric length 
spaces the corresponding notion would be that of a convex neighborhood (in the sense of 
\cite[Def.\ 3.6.5]{BBI:01}).

\begin{defi}\label{def-loc-LpLS}
 A \LpLS \Xll is called \emph{localizable} if $\forall x\in X$ there is an open neighborhood $\Omega_x$ of $x$ in $X$
 with the following properties:
 \begin{enumerate}[label=(\roman*)]
  \item \label{def-loc-LpLS-cau-comp} There is a $C>0$ such that $L^d(\gamma)\leq C$ for all causal curves $\gamma$ 
contained in $\Omega_x$ (hence $X$ is 
$d$-compatible).
  \item \label{def-loc-LpLS-om-con} There is a continuous map $\omega_x\colon \Omega_x \times \Omega_x\rightarrow 
[0,\infty)$ such that\\
$(\Omega_x, d\rvert_{\Omega_x\times\Omega_x}, \ll\rvert_{\Omega_x\times \Omega_x}, \leq\rvert_{\Omega_x\times\Omega_x}, 
\omega_x)$ is a \LpLS with the following non-triviality condition: for every $y\in\Omega_x$ we have 
$I^\pm(y)\cap\Omega_x\neq\emptyset$.
  \item \label{def-loc-LpLS-max-cc} For all $p,q\in \Omega_x$ with $p<q$ there is a future-directed causal curve 
$\gamma_{p,q}$ from $p$ to $q$ that is 
maximal in $\Omega_x$ and satisfies
\begin{equation}
 L_\tau(\gamma_{p,q}) = \omega_x(p,q) \leq \tau(p,q)\,.
\end{equation}
(That the curve $\gamma_{p,q}$ is maximal in $\Omega_x$ means that for every other future-directed causal curve $\lambda$ 
connecting $p$ and $q$ with image contained in $\Omega_x$ we have that $L_\tau(\gamma_{p,q}) \geq L_\tau(\lambda)$.) 
 \end{enumerate}
 We call such a neighborhood $\Omega_x$ a \emph{localizing neighborhood of $x$}.
If, in addition, the neighborhoods $\Omega_x$ can be chosen such that
\begin{enumerate}
 \item[(iv)]\label{def-loc-LpLS-4} Whenever $p,q\in\Omega_x$ satisfy $p\ll q$ then $\gamma_{p,q}$ is timelike and strictly longer than any future-directed 
causal curve in $\Omega_x$ from $p$ to $q$ that contains a null segment,
\end{enumerate}
then \Xll is called {\em regularly localizable}. Finally, if every point $x\in X$ has a neighborhood basis
of open sets $\Omega_x$ satisfying (i)--(iii), respectively (i)--(iv), then \Xll is called {\em strongly localizable} respectively
{\em SR-localizable}. 
\end{defi}
\begin{prop}\label{prop-L-tau-up-sc}
 Let \Xll be a strongly causal and localizable \LpLSn. Then $L_\tau$ is upper semicontinuous, i.e., if $(\gamma_n)_n$ is a 
sequence of future-directed causal curves (defined on $[a,b]$) converging uniformly to a future-directed causal curve 
$\gamma\colon[a,b]\rightarrow X$, then
\begin{equation}
 L_\tau(\gamma)\geq \limsup_n L_\tau(\gamma_n)\,.
\end{equation}
\end{prop}
\begin{pr}
By strong causality (and Lemma \ref{lem-cc} \ref{lem-cc-str-cau}) every point $x\in X$ has an open neighborhood $U_x\subseteq 
\Omega_x$ such that any causal curve with endpoints in $U_x$ is contained in $\Omega_x$.

Let $\eps>0$, then there is a partition $(t_i)_{i=0}^N$ of $[a,b]$ such that
 \begin{equation}
 \sum_{i=0}^{N-1} \tau(\gamma(t_i),\gamma(t_{i+1})) < L_\tau(\gamma) + \frac{\eps}{2}\,.
 \end{equation}
By making the partition finer (and by the reverse triangle inequality) we can assume that $\gamma(t_i),\gamma(t_{i+1})\in 
U_{x_i}\subseteq \Omega_{x_i}$ for some $x_i\in \gamma([a,b])$, $i=0,\ldots,N-1$.

Thus
\begin{equation}\label{deltadef}
  \sum_{i=0}^{N-1} \tau(\gamma(t_i),\gamma(t_{i+1})) \geq  \sum_{i=0}^{N-1} \omega_{x_i}(\gamma(t_i),\gamma(t_{i+1})) =: 
\Delta.
\end{equation}
Now we choose $n_0\in\N$ such that for all $n\geq n_0$ we have $\gamma_n(t_i),\gamma_n(t_{i+1})\in U_{x_i}$ and 
$|\omega_{x_i}(\gamma(t_i),\gamma(t_{i+1})) - \omega_{x_i}(\gamma_n(t_i),\gamma_n(t_{i+1}))|< \frac{\eps}{2 N}$ for 
$i=0,\ldots,N-1$. By construction $\gamma_n([t_i,t_{i+1}]) \subseteq \Omega_{x_i}$. Localizability then implies that for every 
$i=0,\ldots,N-1$ there is a future-directed causal curve $\lambda_i$ from $\gamma_n(t_i)$ to 
$\gamma_n(t_{i+1})$ that is maximal in $\Omega_{x_i}$ with $\omega_{x_i}(\gamma_n(t_i),\gamma_n(t_{i+1})) = L_\tau(\lambda_i)$ and so
\begin{align*}
 \Delta &\geq  \sum_{i=0}^{N-1} \omega_{x_i}(\gamma_n(t_i),\gamma_n(t_{i+1})) - \frac{\eps}{2} =  \sum_{i=0}^{N-1} 
L_\tau(\lambda_i) - \frac{\eps}{2}\\
&\geq  \sum_{i=0}^{N-1} L_\tau(\gamma_n\rvert_{[t_i,t_{i+1}]}) - \frac{\eps}{2} = 
L_\tau(\gamma_n) - \frac{\eps}{2}\,,
\end{align*}
where in the last step we used the additivity of the $\tau$-length proved in Lemma \ref{lem-L-tau-add}. Together with 
\eqref{deltadef} this yields $L_\tau(\gamma)\geq L_\tau(\gamma_n) - \eps$ for every $n\geq n_0$.
\end{pr}

\begin{thm}\label{thm-max-cc-char}
In a regularly localizable Lorentzian pre-length space, maximal causal curves have a causal character, i.e., if for a (future-directed) maximal causal 
curve $\gamma\colon[a,b]\rightarrow X$ there are $a\leq t_1 < t_2\leq b$ with $\gamma(t_1)\ll\gamma(t_2)$, then $\gamma$ 
is timelike. Otherwise it is null.
\end{thm}
\begin{pr}
First we establish that it suffices to show the claim for $t_1 = a$ and $t_2=b$. Indeed, let $\gamma\colon[a,b]\rightarrow 
X$ be a future-directed maximal causal curve and assume that there are $a<t_1<t_2<b$ with $x:=\gamma(t_1)\ll\gamma(t_2)=:y$. 
Thus $\gamma(a)\leq x\ll y \leq \gamma(b)$ and hence by push-up (Lemma \ref{lem-pup}) we conclude that 
$\gamma(a)\ll\gamma(b)$.

We begin the main part of the proof by showing that the claim follows if there exist points $\gamma(t_1)=x\ll 
y=\gamma(t_2)$ such that $\gamma([t_1,t_2])$ lies in a regularly localizing neighborhood $\Omega$ as in Definition 
\ref{def-loc-LpLS}. In fact, since $\gamma$ is maximizing on $[a,b]$, it also is on $[t_1,t_2]$ (Proposition 
\ref{prop-max-prop} \ref{prop-max-sub-int}), and we claim that $\gamma\rvert_{[t_1,t_2]}$ is timelike. 
Otherwise there would exist $t_1\leq s_1 < s_2\leq t_2$ such that $r_1:=\gamma(s_1) \not\ll\gamma(s_2)=:r_2$, implying 
that $\tau(r_1,r_2)=0$ and thus by maximality $L_\tau(\gamma\rvert_{[s_1,s_2]})=0$, i.e., $\gamma\rvert_{[s_1,s_2]}$ is null.
By regular localizability (in $\Omega$) we know that there is a future-directed timelike curve $\gamma_{x,y}$ in $\Omega$ from 
$x$ to $y$ that is strictly longer than $\gamma\rvert_{[t_1,t_2]}$, since the latter contains the null segment 
$\gamma\rvert_{[s_1,s_2]}$. This is a contradiction to the maximality of $\gamma$. We now cover $\gamma([t_1,b])$ by 
finitely many regularly localizing neighborhoods $\Omega_1=\Omega, \dots, \Omega_N$ as in Definition \ref{def-loc-LpLS} and 
pick $t_i$ ($i=3,\dots,N+2$), $t_{N+2}=b$, such that $\gamma(t_i) \in \Omega_{i-2}\cap \Omega_{i-1}$ for 
$i=3,\dots,N+1$ and $\gamma([t_i,t_{i+1}])\subseteq \Omega_{i-1}$ for $i=2,\dots N+1$. It then follows as above that, 
since $\gamma(t_1)\ll \gamma(t_2)\le \gamma(t_3)$, and hence $\gamma(t_1)\ll\gamma(t_3)$, $\gamma$ must be 
timelike on $[t_1,t_3]$. Then picking some $t'<t_3$ such that $\gamma([t',t_3])\subseteq \Omega_2$ we find ourselves in the 
same situation as before, only with $[t',t_3]$ replacing $[t_1,t_2]$. Consequently, we can iterate the procedure and obtain 
that $\gamma$ is timelike on $[t_1,b]$. Since we may symmetrically argue into the past, $\gamma$ must in fact be timelike on 
all of $[a,b]$.

It remains to show that points $x\ll y$ as above always exist on $\gamma$. Since $0<\tau(p,q)=L_\tau(\gamma)$, by Lemma 
\ref{def-cc-len} it follows that $\gamma\rvert_{[a,m]}$ or $\gamma\rvert_{[m,b]}$ has strictly positive $\tau$-length, where 
$m=\frac{1}{2}(b-a)$. Iterating this bisection it follows that for any $\delta>0$ there exist $t_1<t_2$ in $[a,b]$ such that 
$|t_1-t_2|<\delta$ and $\tau(\gamma(t_1),\gamma(t_2)) = L_\tau(\gamma\rvert_{[t_1,t_2]})>0$, and so $\gamma(t_1)\ll 
\gamma(t_2)$. We now cover $\gamma([a,b])$ by finitely many regularly localizing neighborhoods as in Definition 
\ref{def-loc-LpLS} and let $\eps$ be a Lebesgue number of this cover. Since $\gamma$ is uniformly continuous, by choosing 
$\delta$ small enough we can guarantee that $d(\gamma(t_1),\gamma(t_2))<\eps$, and so both points lie in one of the 
neighborhoods from the cover.
\end{pr} 
\begin{ex}\label{ex-funnel}
\begin{enumerate}[label=(\roman*)] \ 
\item \label{ex-funnel-1} {\em Causal/timelike funnels} 

In Minkowski space $\R^n_1$, let $\lambda$ be a future-directed causal curve connecting two points
$p$ and $q$. Let $X$ be the union of $J^-(p)$, $J^+(q)$ (or, alternatively, of $I^-(p)$, $I^+(p)$) and the image of $\lambda$.
\begin{figure}[h!]
\begin{center}
\includegraphics[width=80mm, height= 60mm]{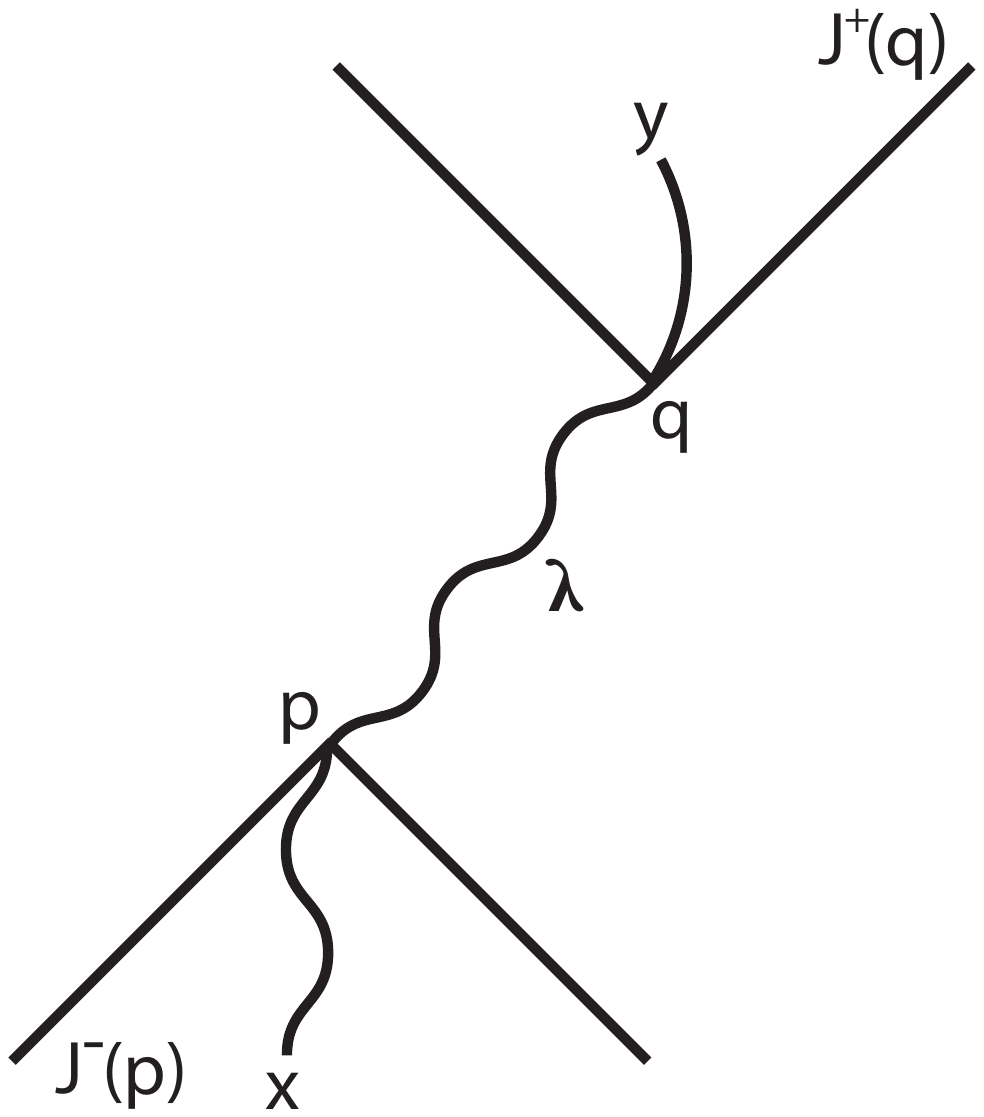}
\end{center}
\caption{A causal funnel.}
\end{figure}
For $x,y\in X$, let $x\le y$ if $x$ can be connected to $y$ within $X$ by a curve that is future-directed causal in $\R^n_1$, and
let $x\ll y$ if this curve can be chosen to contain a timelike segment. Define $\tau(x,y)$ to be the supremum over all
lengths of such curves connecting $x$ and $y$ if such curves exist, and $0$ otherwise. Also, let $d$ be the
restriction of the standard metric on $\R^n$.
Then it is easily verified that \Xll is a Lorentzian pre-length space.
If $\lambda$ is null (hence also null in the sense of Definition \ref{def-cau-cur}), 
$x\ll p$, and $q\le y$, then the maximal curve from $x$ to $y$ necessarily changes its causal character.
\item As can be seen from Corollary \ref{cg_branch} below, even in spacetimes with continuous metrics it can happen that 
maximal causal curves change their causal character. 
\end{enumerate}
\end{ex} 
We next turn to a fundamental property of smooth spacetimes, namely the push-up principle (cf.\ \cite{Chr:11}): 
causal curves that connect timelike related points and contain a null segment
can be deformed into timelike curves with the same endpoints and strictly greater length. This principle can 
be extended to the current setting as follows:
\begin{thm}\label{thm-length-increasing-push-up}
Let \Xll be a regularly localizable Lorentzian pre-length space, 
and let $\gamma: [a,b]\to X$ be a future-directed causal curve with $L_\tau(\gamma)>0$. %$\gamma(a)\ll \gamma(b)$. 
If $\gamma|_{[c,d]}$ is null on some (non-trivial) sub-interval $[c,d]$ of $[a,b]$, 
then there exists a strictly longer future-directed timelike 
curve $\sigma$ from $\gamma(a)$ to $\gamma(b)$. If $X$ is even (SR)-localizable, then $\sigma$
can be chosen to lie in any given neighborhood of $\gamma([a,b])$. 
\end{thm} 
\begin{pr} Without loss of generality we may suppose that $a<c$ and $d=b$, the other cases can be reduced to this one or
proved analogously. Let $t_1:= \inf \{t \in [a,b] : L_\tau(\gamma|_{[t,b]}) = 0\}$, then $a<t_1$: Suppose,
to the contrary, that $t_1=a$ and let $\Omega_{\gamma(a)}$ be a regularly localizing neighborhood of $\gamma(a)$.
Then since $L_\tau(\gamma|_{[s,t]})=\omega_{\gamma(a)}(\gamma(s),\gamma(t))$ depends continuously on $s$ and $t$
for $s$, $t$ small, it would follow that $L_\tau(\gamma)=0$, contradicting our assumption.

Now let $\Omega$ be a regularly localizing neighborhood of $\gamma(t_1)$. 
Then we can pick $t_0<t_1<t_2$ sufficiently
close to $t_1$ to secure $\gamma(t_0), \gamma(t_2) \in \Omega$. Also, $\gamma(t_0)\ll \gamma(t_2)$
and $\gamma([t_1,t_2])$ contains a null segment. Thus we can connect $\gamma(t_0)$ to $\gamma(t_2)$
by a future-directed timelike curve $\sigma$ in $\Omega$ that is strictly longer than $\gamma|_{[t_0,t_2]}$. 

Similarly to the proof of Theorem \ref{thm-max-cc-char} we cover $\gamma([t_0,b])$ by 
finitely many regularly localizing neighborhoods $\Omega_1=\Omega_2:=\Omega,\Omega_3, \dots, \Omega_N$ 
as in Definition \ref{def-loc-LpLS} and 
pick $t_3< \dots < t_{N}=b$ in $(t_2,b]$ such that $\gamma([t_i,t_{i+1}])\subseteq \Omega_{i+1}$ for $i=0,\dots N-1$.
Next, we choose a point $p$ on $\sigma$ that lies in $\Omega_3$ and is timelike related to $\gamma(t_2)$
and concatenate $\sigma$ from $p$ onward to a maximal curve from $p$ to $\gamma(t_3)$ within $\Omega_3$. 
Iterating this procedure, we obtain a timelike curve from $\gamma(t_0)$ to $\gamma(b)$ that is strictly
longer than $\gamma|_{[t_0,b]}$. Analogously, we can argue for $\gamma|_{[a,t_0]}$ to construct the 
claimed curve. Finally, if $X$ is (SR)-localizable then the regularly localizing neighborhoods and thereby
the timelike curve constructed above can be chosen to lie within any prescribed neighborhood of the image of $\gamma$.
\end{pr}
Recalling Lemma \ref{lem-rect-is-timelike}, we obtain the following generalization of 
\cite[Cor.\ 2.4.16]{Chr:11}:
\begin{cor}\label{car-push-up}
Let \Xll be a regularly localizable Lorentzian pre-length space, 
and let $\gamma: [a,b]\to X$ be a future-directed causal curve such that for some $a\le c<d\le b$,
$\gamma|_{[c,d]}$ is rectifiable. Then there exists a timelike future-directed curve from $\gamma(a)$
to $\gamma(b)$. If $X$ is even (SR)-localizable, then this curve
can be chosen to lie in any given neighborhood of $\gamma([a,b])$.
\end{cor}

\subsection{Lorentzian length spaces} 
Finally, we have the concepts at hand to define the following notion of \emph{intrinsic} time separation function.

\begin{defi}\label{bigtau}
Let \Xll be a locally causally closed, causally path connected and localizable \LpLS and let $x,y\in X$. Then set
\begin{equation*}
 \mathcal{T}(x,y):= \sup\{L_\tau(\gamma):\gamma \text{ future-directed causal from }x \text{ to } y\}\,, 
\end{equation*}
if the set of future-directed causal curves from $x$ to $y$ is not empty. Otherwise set $\mathcal{T}(x,y):=0$. We call $X$ 
a \emph{Lorentzian length space} if $\mathcal{T} = \tau$. If, in addition, $X$ is regularly
localizable, then it is called a \emph{regular} Lorentzian length space.
\end{defi}

\begin{rem}\label{lls_rem}\

\begin{enumerate}[label=(\roman*)]
\item The above definition is a close analogue of the notion of length spaces in metric geometry: 
A metric space $(X,d)$ is a length space if for any points $x,y\in X$, $d(x,y)$ equals the infimum over the
length of all paths connecting them, where the length of a path is defined as the supremum of the lengths of inscribed polygons (cf.\ \cite{BBI:01,Pap:14}).
\item Since a \LLS is causally path connected, the set of all future-directed causal curves connecting two causally related points is never empty.
\item \label{lls_rem-3}
In any Lorentzian pre-length space, ${\mathcal T}(x,y)\le \tau(x,y)$ for all $x, y\in X$:  
This is obvious if $\mathcal{T}(x,y)=0$. If, on the
other hand, $\mathcal{T}(x,y)>0$, then for any $\eps>0$ there exists a future-directed causal curve 
$\gamma$ from $x$ to $y$ with $\mathcal{T}(x,y)< L_\tau(\gamma)+\eps \le \tau(x,y)+\eps$.
\end{enumerate}

\end{rem}

\begin{ex} \label{ex-lls}
\begin{enumerate}[label=(\roman*)]\ 
 \item Let $(M,d^h,\ll,\leq,\tau)$ be the \LpLS induced by a smooth and strongly causal spacetime $(M,g)$, see Example 
\ref{ex-lpls-smo}. Then by Proposition \ref{prop-smo-stc-len} we know that $L_\tau = L_g$, the usual Lorentzian length 
functional. Thus the definition of $\Tau$ is the same as for the time separation function $\tau$ of $(M,g)$, cf.\ 
\cite[Def.\ 14.15]{ONe:83}. Using the exponential map and convex neighborhoods it is not hard to see 
that $(M,d^h,\ll,\leq,\tau)$ is also locally causally closed and regularly localizable. Moreover, causal path connectedness holds
due to the definition of the causal relations. Consequently, $(M,d^h,\ll,\leq,\tau)$ is a 
regular and (SR)-localizable \LLSn.
\item As in Example \ref{ex-funnel} \ref{ex-funnel-1}, let $X$ be a timelike funnel. If the connecting curve $\lambda$ is
timelike, then \Xll is causally path connected and causally closed. It then readily follows that \Xll is
a strongly localizable Lorentzian length space.
\item 
In Section \ref{sec-app} we will give further examples of Lorentzian length spaces, and in particular we will show 
that spacetimes of low regularity can be viewed as Lorentzian length spaces,
although not necessarily as {\em regular} Lorentzian length spaces. This connection is our motivation
for the terminology introduced after Definition \ref{def-loc-LpLS} (iv).
\end{enumerate}

\end{ex}

\begin{lem}\label{lem-lls-rec}
 In a \LLS two timelike related points can always be connected via a causal curve of positive $\tau$-length.
\end{lem}
\begin{pr}
Let \Xll be a \LLS and let $x,y\in X$ with $x\ll y$. Then $0<\tau(x,y)=\Tau(x,y)$. Moreover, for every $\eps>0$ there is a 
future-directed causal curve $\gamma$ from $x$ to $y$ such that $L_\tau(\gamma) > \Tau(x,y) - \eps$. By choosing 
$\eps=\frac{\Tau(x,y)}{2}>0$ it follows that $L_\tau(\gamma)>0$. 
\end{pr}

\subsection{The causal ladder for Lorentzian length spaces}

\begin{thm}\label{thm-cau-lad-lls}
 For \LLSn s 
\begin{enumerate}[label=(\roman*)]
 \item causality implies chronology,
 \item non-total imprisonment implies causality,\label{thm-cau-lad-lls-nti-cau}
 \item strong causality implies non-total imprisonment,
 \item strong causality is equivalent to the non-existence of almost closed causal curves (i.e., the converse to Lemma 
\ref{lem-cc} \ref{lem-cc-str-cau} holds for \LLSn s), and\label{thm-cau-lad-lls-str-cau-no-acc}
 \item global hyperbolicity implies strong causality.
\end{enumerate}
\end{thm}
\begin{pr}
Let \Xll be a \LLSn, which for brevity we just denote by $X$. Lemma \ref{lem-cpc-i} shows that $X$ is 
interpolative.
 \begin{enumerate}[label=(\roman*)]
  \item Let $X$ be causal, then Lemma \ref{lem-cc} \ref{lem-cc-cai-ch} establishes that $X$ is chronological.
  \item Let $X$ be non-totally imprisoning. Assume that $X$ is not causal, then by Lemma 
\ref{lem-cpc-cc} \ref{lem-cpc-cc-cau} we know that there is a closed causal curve $\gamma\colon[a,b]\rightarrow X$. Since 
$\gamma$ is not constant we have that $L^d(\gamma)>0$. By going infinitely often around this loop we get a causal curve 
$\tilde{\gamma}$ such that $L^d(\tilde{\gamma})=\infty$ and whose image is contained in the compact set $\gamma([a,b])$ --- a 
contradiction to $X$ being non-totally imprisoning.

\item Let $X$ be strongly causal and assume that $X$ is totally imprisoning. By Corollary \ref{cor-nti-inext} 
this means that there is a compact set $K\Subset X$ and an inextendible (future-directed) causal curve 
$\gamma\colon[0,\infty)\rightarrow X$ contained in $K$. Moreover, by Lemma \ref{lem-inext-cc-infty} we know that 
$\lim_{t\nearrow \infty}\gamma(t)$ does not exist. However, for any sequence that convergences to $\infty$ there is a 
convergent subsequence $(t_n)_n$ with $t_n\nearrow \infty$ and $\lim_{n\to\infty}\gamma(t_n)=:p$, since $\gamma$ is 
contained in the compact set $K$. Now, since $\lim_{t\nearrow \infty}\gamma(t)$ does not exist there is another 
sequence $(s_n)_n$ with $s_n\nearrow \infty$ with $\lim_{n\to\infty}\gamma(s_n)=:q\neq p$. Let $U$ be a neighborhood of $p$ 
that does not contain $q$. By strong causality and Lemma \ref{lem-cc} \ref{lem-cc-str-cau} there exists a neighborhood $V$ 
of $p$ with $V\subseteq U$ and such that any causal curve with endpoints in $V$ is contained in $U$. There is an $n_0\in \N$ 
such that $\gamma(t_n)\in V$ for all $n\geq n_0$. By mixing the sequences $(t_n)_n$ and $(s_n)_n$ to get a strictly 
monotonically increasing sequence $(r_n)_n$ one can find $n_1 < n_2 < n_3$ such that $\gamma(r_{n_1}),\gamma(r_{n_3})\in V$ 
and $\gamma(r_{n_2})\notin U$. This is a contradiction since $\gamma\rvert_{[t_{r_1},t_{r_3}]}$ is a causal curve with 
endpoints in $V$ that leaves $U$.

\item Let $X$ be such that for all $x\in X$, for every neighborhood $U$ of $x$, there is a neighborhood 
$V\subseteq U$ of $x$ such that for every causal curve $\gamma\colon[a,b]\rightarrow X$ with $\gamma(a),\gamma(b)\in V$ one 
has $\gamma([a,b])\subseteq U$. Assume to the contrary that $X$ is not strongly causal, i.e., there is a $p\in X$ and a $\delta>0$ 
such that for all $A\in\mathcal{A}$ (the subbase for the Alexandrov topology, cf.\ Subsection \ref{subsec-top}) with $p\in 
A$ one has $A\not\subseteq B^d_\delta(p)$. Now the assumptions yield that there is a $d$-neighborhood $V$ of $p$, 
$V\subseteq B^d_\delta(p)$ such that all causal curves with endpoints in $V$ are contained in $B^d_\delta(p)$. Let $\Omega$ 
be a localizing neighborhood for $p$, then $I^\pm(p)\cap \Omega\neq\emptyset$, thus by causal path-connectedness there is a 
timelike curve $\gamma$ through $p$. Now choose $p^-,p^+\in\gamma([a,b])\cap V$ with $p^-\ll p\ll p^+$. Then, $p\in I^+(p^-)\cap 
I^-(p^+)\in \mathcal{A}$ but on the other hand $I^+(p^-)\cap I^-(p^+)\subseteq B^d_\delta(p)$ --- a contradiction.

\item Let $X$ be globally hyperbolic and assume that $X$ is not strongly causal, i.e., there is a point $x\in X$ and a 
neighborhood $U$ of $x$ such that for all neighborhoods $V$ of $x$ with $V\subseteq U$ there is a causal curve with 
endpoints in $V$ that leaves $U$. As above there is a timelike curve $\lambda$ through $x$, hence we can choose $p,q\in 
U$ on $\lambda$ with $p\ll x \ll q$. Moreover, since $I^+(p)\cap I^-(q)$ is open there is a $\delta_0>0$ such that for all 
$0<\delta\leq \delta_0$ we have $B^d_\delta(x)\subseteq I^+(p)\cap I^-(q) \subseteq J^+(p)\cap J^-(q)$, which is compact by 
assumption. Let $n_0\in\N$ with $\frac{1}{n_0}< \delta_0$, then for all $n\geq n_0$ there is a future-directed causal curve
$\gamma_n\colon  [a_n,b_n] \rightarrow X$ with $\gamma_n(a),\gamma_n(b)\in B^d_{1/n}(x)$, $\gamma_n([a_n,b_n])\subseteq 
J(p,q)$ and $\gamma_n([a_n,b_n])\not\subseteq U$. Thus we can apply the limit curve theorem \ref{thm-lim-curve} to obtain a closed 
causal curve (which is not constant since it leaves $U$) --- a contradiction to non-total imprisonment via point 
\ref{thm-cau-lad-lls-nti-cau} above.
 \end{enumerate}
\end{pr}

\subsection{Geodesic length spaces}

\begin{defi}
 A \LpLS \Xll is called \emph{geodesic} if for all $x,y\in X$ with $x<y$ there is a future-directed causal curve 
$\gamma$ from $x$ to $y$ with $\tau(x,y)=L_\tau(\gamma)$ (hence maximizing).
\end{defi}

\begin{thm}\label{thm-gh-LLS-tau-fin-cont}
 Let \Xll be a globally hyperbolic \LLSn, then $\tau$ is finite and continuous.
\end{thm}
\begin{pr}
We first show that $\tau$ is finite. Suppose, to the contrary, that $\tau(p,q)=\infty$
for two points $p, q\in X$. Hence there exists a sequence $\gamma_k$ of future directed
causal curves from $p$ to $q$ with $L_\tau(\gamma_k)\to \infty$. Moreover, all of them are
contained in the compact set $J(p,q)$ and so by the limit curve theorem \ref{thm-lim-curve}
and the upper semicontinuity of $L_\tau$ (Proposition \ref{prop-L-tau-up-sc}) 
the corresponding limit curve $\gamma$ has infinite $\tau$-length. Since $X$ is non-totally
imprisoning by Theorem \ref{thm-cau-lad-lls}, $\gamma$ is defined on a compact interval,
say $[0,b]$. Now fix a subdivision $0=t_0<\ldots<t_N=b$ of $[0,b]$ such that each
$\gamma|_{[t_i,t_{i+1}]}$ is contained in a localizing neighborhood $U_i$. By strong causality,
the local time separation function on $U_i$ (which is finite) can be assumed to agree 
with $\tau$. But then each $L_\tau(\gamma|_{[t_i,t_{i+1}]})\le \tau(\gamma(t_i),\gamma(t_{i+1}))
<\infty$, a contradiction.

Next, to establish that $\tau$ is continuous, we assume to the contrary that 
$\tau$ is not upper semicontinuous at $(p,q)\in X\times X$. 
Thus there exist some $\delta>0$ and sequences $p_n \to p$, $q_n\to q$ 
such that
\begin{equation}
 \tau(p_n,q_n)\geq \tau(p,q) + \delta\,,
\end{equation}
for all $n\in\N$. Since $\tau(p,q)\geq 0$ we have that $\tau(p_n,q_n)>0$ and hence $p_n\ll q_n$ for all $n\in\N$. 
Furthermore, for $n\geq 1$ there is a future-directed causal curve $\gamma_n$ from $p_n$ to $q_n$ with $L_\tau(\gamma_n) > 
\Tau(p_n,q_n)- \frac{1}{n} = \tau(p_n,q_n)-\frac{1}{n}$. Note that, by strong causality, this shows that $p\neq q$.
 Again by strong causality and localizability, there are $p_-,q_+\in X$ such that $p\in I^+(p_-)$ and $q\in I^-(q_+)$. So 
there is an $n_0\in\N$ such that $p_n,q_n\in J(p_-,q_+)$ for all $n\geq n_0$. By global hyperbolicity $J(p_-,q_+)$ is compact 
and the image of $\gamma_n$ is contained in $J(p_-,q_+)$ for all $n\geq n_0$ and so by the Limit curve theorem \ref{thm-lim-curve} and 
$p\neq q$ we get that there is a subsequence $(\gamma_{n_k})_k$ of $(\gamma_n)_n$ that converges uniformly to a 
future-directed causal curve $\gamma$ from $p$ to $q$. Moreover, by construction and the upper semicontinuity of $L_\tau$ 
(Proposition \ref{prop-L-tau-up-sc}) this yields that
\begin{align*}
 \Tau(p,q) &\geq L_\tau(\gamma)\geq \limsup_k L_\tau(\gamma_{n_k})\\
&\geq  \limsup_k \Big(\tau(p_{n_k},q_{n_k}) - \frac{1}{n_k}\Big) 
\geq \tau(p,q) + \delta = \Tau(p,q) + \delta > \Tau(p,q)\,,
\end{align*}
a contradiction. Thus $\tau$ is continuous.
\end{pr}

\begin{rem}
 Finiteness of $\tau$ precludes the pathological situation where a maximal curve could have infinite length.
\end{rem}

Finally, we obtain the following generalization of the Avez-Seifert theorem:

\begin{thm}\label{thm-gh-LLS-geo}
 Any globally hyperbolic \LLS \Xll is geodesic.
\end{thm}
\begin{pr}
 By Theorem \ref{thm-gh-LLS-tau-fin-cont} we know that $\tau$ is finite and continuous. Let $x,y\in X$ with $x< y$, then 
$\tau(x,y)<\infty$ and we get a sequence $(\gamma_n)_n$ of future-directed causal curves $\gamma_n\colon[a,b]\rightarrow X$ 
from $x$ to $y$ such that $L_\tau(\gamma_n)\to \tau(x,y)$. These curves are all contained in the compact set $J(x,y)$ and 
so by the Limit curve theorem \ref{thm-lim-curve} we get a limit curve $\gamma$ from $x$ to $y$ with $L_\tau(\gamma) = 
\limsup_n L_\tau(\gamma_n) = \lim_n L_\tau(\gamma_n) = \tau(x,y)$. Thus $\gamma$ is a maximal future-directed causal curve 
from $x$ to $y$.
\end{pr}

\subsection{Parametrization by arclength}

We will now establish that a rectifiable curve (which is timelike) can be parametrized with respect to $\tau$-arclength. The 
only drawback is that this parametrization need not be Lipschitz continuous. Thus the resulting curve will not be a causal 
curve in the sense of Definition \ref{def-cau-cur}. To handle this issue we introduce the following notion.

\begin{defi}
Let \Xll be a \LpLS and let $\gamma\colon[a,b]\rightarrow X$ be a future-directed causal curve. A \emph{weak 
parametrization} of $\gamma$ is a curve of the form $\gamma\circ\phi\colon[c,d]\rightarrow X$, where 
$\phi\colon[c,d]\rightarrow [a,b]$ is continuous and strictly monotonically increasing.
\end{defi}
Note that if $\phi$ is Lipschitz continuous, then $\gamma\circ\phi$ is Lipschitz continuous, hence a causal curve. 
Moreover, the $\tau$-length of a weak parametrization can be defined as in Definition \ref{def-cc-len}, and Lemma 
\ref{lem-L-tau-inv} shows that the $\tau$-length is invariant under such a reparametrization, too.

\begin{lem}
Let \Xll be a \LpLSn, let $\gamma\colon[a,b]$ $\rightarrow X$ be a future-directed causal curve and let
$\lambda:=\gamma\circ\phi\colon[c,d] \rightarrow X$ be a weak parametrization of $\gamma$. Then $\lambda$ has the same causal 
character as $\gamma$.
\end{lem}
\begin{pr}
 We show the case when $\gamma$ is causal, the timelike case is completely analogous. Let $c\leq t_1 < t_2 \leq d$, 
then $a\leq \phi(t_1)<\phi(t_2)\leq b$ and so $\lambda(t_1)=\gamma(\phi(t_1))\leq \gamma(\phi(t_2)) = \lambda(t_2)$.
\end{pr}

\begin{lem}\label{lem-phi}
 Let \Xll be a \LpLS and let $\gamma\colon[a,b]\rightarrow X$ be a future-directed causal curve with 
$L:=L_\tau(\gamma)<\infty$. Then the map $\phi\colon[a,b]\rightarrow [0,L]$, $t\mapsto L_\tau(\gamma\rvert_{[a,t]})$ is 
monotonically increasing. Moreover, if the time separation function $\tau$ is continuous and satisfies $\tau(x,x)=0$ for all 
$x\in X$ then $\phi$ is continuous.
\end{lem}
\begin{pr}
First we show that $\phi$ is monotonically increasing. Let $a\leq s < t \leq b$ and let $a=t_0<t_1<\ldots<t_N=t$ be a 
partition of $[a,t]$. If there is a $k\in\{1,\ldots,N\}$ such that $t_k = s$, then $(t_i)_{i=0}^k$ is a 
partition of $[a,s]$ and thus $\phi(s) \leq \sum_{i=0}^{k-1}\tau(\gamma(t_i),\gamma(t_{i+1})) \leq 
\sum_{i=0}^{N-1}\tau(\gamma(t_i),\gamma(t_{i+1}))$ %(because $\tau\geq 0$). 
On the other hand, if there is no such $k$, define 
$j:=\max\{1\leq i \leq N: t_i< s\}$. Then $(t_i)_{i=0}^j \cup \{s\}$ is a partition of $[a,s]$. This yields
$\phi(s) \leq \sum_{i=0}^{j-1}\tau(\gamma(t_i),\gamma(t_{i+1})) + \tau(\gamma(t_j),\gamma(s)) \leq 
\sum_{i=0}^{N-1}\tau(\gamma(t_i),\gamma(t_{i+1}))$, where we again used that $\tau\geq 0$ and the reverse triangle 
inequality. Taking the infimum over all partitions of $[a,t]$ gives $\phi(s)\leq  L_\tau(\gamma\rvert_{[a,t]})=\phi(t)$.

To show continuity of $\phi$ at any $t\in [a,b]$, we make use of the continuity of the maps $y\mapsto \tau(x,y)$ and 
$y\mapsto \tau(y,x)$ for $x\in X$ fixed. Let $t\in[a,b]$ and $\eps>0$, then there is a neighborhood $U$ of $\gamma(t)$ in 
$X$ such that for all $y\in U$
\begin{equation}\label{eq-lem-phi}
 \tau(\gamma(t),y) < \eps \text{ and } \tau(y,\gamma(t))<\eps\,,
\end{equation}
since $\tau\geq 0$ and $\tau(\gamma(t),\gamma(t))=0$. By the continuity of $\gamma$, there is a $\delta>0$ such that 
$\gamma((t-\delta,t+\delta))\subseteq U$. For $s\in (t-\delta,t]$ we have by Lemma \ref{lem-L-tau-add} and \eqref{eq-lem-phi}
\begin{align*}
 |\phi(t)-\phi(s)|= L_\tau(\gamma\rvert_{[s,t]}) \leq 
\tau(\gamma(s),\gamma(t))<\eps\,.
\end{align*}
Analogously for $s\in[t,t+\delta)$ we have $|\phi(t)-\phi(s)| = L_\tau(\gamma\rvert_{[t,s]}) \leq 
\tau(\gamma(t),\gamma(s))\\ < \eps$.
\end{pr}

\begin{prop}\label{prop-rec-reparam}
 Let \Xll be a \LpLS with $\tau$ continuous and $\tau(x,x)=0$ for all $x\in X$. Let 
$\gamma\colon[a,b]\rightarrow X$ be a future-directed rectifiable curve with $L:=L_\tau(\gamma)<\infty$. Then 
there exists a weak parametrization $\tilde{\gamma}$ of $\gamma$ such that $\tilde{\gamma}$ is parametrized with respect to 
$\tau$-length, i.e., $\tilde{\gamma}\colon[0,L]\rightarrow X$ with $L_\tau(\tilde{\gamma}\rvert_{[0,s]}) = s$ for all 
$s\in[0,L]$.
\end{prop}
\begin{pr}
 Define $\phi\colon[a,b]\rightarrow [0,L]$, $t\mapsto L_\tau(\gamma\rvert_{[a,t]})$ as in Lemma \ref{lem-phi}. Then $\phi$ 
is strictly monotonically increasing and continuous and thus gives rise to a weak parametrization $\tilde{\gamma}:= \gamma 
\circ \phi^{-1}\colon [0,L]\rightarrow [a,b]$. Note that Lemma \ref{lem-L-tau-inv} applies also to weak parametrizations and 
hence we conclude that $L_\tau(\tilde{\gamma}\rvert_{[0,s]}) = L_\tau(\gamma\rvert_{[a,\phi^{-1}(s)]}) = \phi(\phi^{-1}(s)) = 
s$.
\end{pr}

\begin{cor}\label{cor-max-tl-par}
Let \Xll be a \LpLS with $\tau$ continuous and $\tau(x,x)=0$ for all $x\in X$. Then a maximal timelike curve $\gamma$ 
with finite $\tau$-length has a weak parametrization $\lambda$ such that $\tau(\lambda(s_1),\lambda(s_2)) = s_2 - s_1$ for 
all $s_1 < s_2$ in the corresponding interval.
\end{cor}
\begin{pr}
Let $\gamma$ be timelike and maximal. Then by Proposition \ref{prop-max-prop} \ref{prop-max-tl-rec} $\gamma$ is 
rectifiable and hence by Proposition \ref{prop-rec-reparam} there is a weak parametrization $\lambda=\gamma 
\circ \phi^{-1}$ on $[0,L_\tau(\gamma)]$ such that $L_\tau(\lambda\rvert_{[0,s]}) = s$. Moreover, as noted above, Lemma 
\ref{lem-L-tau-inv} applies also to weak parametrizations, hence we have that $L_\tau(\gamma) = L_\tau(\lambda)$ and thus for 
$0\leq s_1 < s_2 \leq L_\tau(\gamma)$ we get 
\begin{align*}
 s_2 - s_1 &= L_\tau(\lambda\rvert_{[s_1,s_2]}) = 
L_\tau(\gamma\rvert_{[\phi^{-1}(s_1),\phi^{-1}(s_2)]})\\ &= \tau(\gamma(\phi^{-1}(s_1)),\gamma(\phi^{-1}(s_2))) = 
\tau(\lambda(s_1),\lambda(s_2))\,.
\end{align*}
\end{pr}

\section{Curvature bounds via triangle comparison}\label{sec-curvature}
In close analogy to the theory of CAT$(k)$- and Alexandrov spaces, in this section we introduce spaces whose 
curvature is bounded above or below, in terms of triangle comparison with respect to Lorentzian model spaces
of constant curvature. The comparison conditions will be formulated with respect to the time separation function $\tau$.

\subsection{Timelike geodesic triangles}
We begin by considering \emph{timelike geodesic triangles} in a \LLS and compare them to timelike geodesic triangles in a model space 
of constant curvature.

\begin{defi} \label{def_tl_geo_tri}
 A \emph{timelike geodesic triangle} in a \LpLS \Xll is a triple $(x,y,z)\in X^3$ with $x\ll y\ll z$ such that 
$\tau(x,z)<\infty$ and such that the sides are realized by future-directed causal curves:  
there are future-directed causal curves $\alpha$ from $x$ to $y$, $\beta$ from
$y$ to $z$, and $\gamma$ from $x$ to $z$ such that $L_\tau(\alpha) = \tau(x,y)$, $L_\tau(\beta) = 
\tau(y,z)$ and $L_\tau(\gamma) = \tau(x,z)$.
\end{defi}

The reason for merely requiring the realizing curves in the previous definition to be causal instead
of timelike is that in general it may happen that maximizing curves change their causal 
character. For a concrete example see Corollary \ref{cg_branch} below.
By Theorem \ref{thm-max-cc-char}, however, such a situation cannot occur in regularly
localizable pre-length spaces.

\begin{lem}\label{lem-LpLS-geo-tri}
 Let \Xll be a \LpLS and let $(x,y,z)$ be a geodesic triangle. Let $\alpha, \beta, \gamma$ be future-directed causal curves from $x$ to $y$, $y$ to $z$, and $x$ to $z$, respectively,
such that $L_\tau(\alpha) = \tau(x,y)=:a$, $L_\tau(\beta) = \tau(y,z)=:b$ and $L_\tau(\gamma) = 
\tau(x,z)=:c$. Then $a<\infty$, $b<\infty$ and $\alpha, \beta, \gamma$ are maximal. 
\end{lem}
\begin{pr}
 By the reverse triangle inequality we get $a = \tau(x,y) < \tau(x,y) + \tau(y,z) \leq \tau(x,z)<\infty$. 
Analogously one shows that $b<\infty$. Let $\alpha$ be defined on the interval $[t_0,t_1]$, then $L_\tau(\alpha) = \tau(x,y) 
=\tau(\alpha(t_0),\alpha(t_1))$, hence $\alpha$ is maximal. Similarly $\beta$ and $\gamma$ are maximal. 
\end{pr}

\begin{rem}\label{triangle_rem}  In the situation of the previous Lemma, if $\tau$ is continuous and 
$\tau(x,x)=0$ for all $x$, then for any $0< s < \tau(x,y)$ there is a point $q$
on the image of $\alpha$ with $\tau(x,q)=s$ (and analogously for $\beta$ and $\gamma$). 
In fact, let $\alpha: [a,b]\to X$. Then since $t\mapsto \tau(x,\alpha(t))$ is
continuous, it attains any value between $0=\tau(x,\alpha(a))$ and $\tau(x,y)=\tau(x,\alpha(b))$.

Moreover, if $\alpha$, $\beta$, $\gamma$ are timelike then Corollary \ref{cor-max-tl-par} 
allows one to obtain weak parametrizations of $\alpha, \beta, \gamma$ with respect to 
$\tau$-length.
\end{rem}

\begin{lem}
 Let \Xll be a globally hyperbolic \LLSn, then any triple of points $(x,y,z)\in X^3$ with $x\ll y\ll z$ is a geodesic 
timelike triangle, whose sides, if timelike, can be weakly parametrized with respect to $\tau$-length. 
\end{lem}
\begin{pr}
 By Theorem \ref{thm-gh-LLS-tau-fin-cont} $\tau$ is finite and continuous, implying in particular $\tau(x,z)<\infty$. By 
Theorem \ref{thm-cau-lad-lls} we know that $X$ is chronological and thus $\tau$ is zero on the diagonal by Lemma
\ref{lem-cc} \ref{lem-cc-ch-tau}. Furthermore, by Theorem \ref{thm-gh-LLS-geo} $X$ is geodesic, hence there are maximal 
causal curves realizing the sides of the triangle, which, if timelike, can be weakly parametrized with respect to $\tau$-length
by Corollary \ref{cor-max-tl-par}.
\end{pr}

\subsection{Model spaces of constant curvature}

Curvature bounds for Lorentzian length spaces will be based on triangle comparison
in relation to model spaces of constant curvature. In the present section we introduce these
model spaces, following \cite{Har:82,AB:08}.

\begin{defi} Let $K\in \R$. By $M_K$ we denote the simply connected two-dimensional Lorentz space form of constant curvature $K$.
\end{defi}

Following the notation of \cite[Ch.\ 8]{ONe:83}, we have
\[
M_K = \left\{ \begin{array}{ll}
\tilde S^2_1(r) & K=\frac{1}{r^2}\\
\R^2_1 & K=0\\
\tilde H^2_1(r) & K= -\frac{1}{r^2}.
\end{array}
\right.
\]
Here,  $\tilde S^2_1(r)$ is the simply connected covering manifold of the two-dimensional Lorentzian pseudosphere
$S^2_1(r)$, $\R^2_1$ is two-dimensional Minkowski space, and $\tilde H^2_1(r)$ is the simply connected covering manifold
of the two-dimensional Lorentzian pseudohyperbolic space.

Concerning the existence of comparison triangles in the model spaces, we may directly utilize 
the {\em Realizability Lemma} \cite[Lemma 2.1]{AB:08} to obtain conditions on a triple $(a,b,c)$
to be realized as the side lengths of a timelike triangle in a model space $M_K$. Below, we
will set $\frac{\pi}{\sqrt{K}}:=\infty$ if $K\le 0$.
\begin{lem}\label{lem-tl-tri-mk} (Realizability) Let $K\in \R$. Let $(a,b,c)\in \R_+^3$ with $c\ge a+b$. 
If $c=a+b$, then let $c<\frac{\pi}{\sqrt{K}}$. 
If $c>a+b$ and
$K<0$, then assume $c<\frac{\pi}{\sqrt{-K}}$.
Then there exists a timelike geodesic triangle in $M_K$ with side lengths $a$, $b$, $c$. 
\end{lem}
\begin{pr} To deduce this from \cite[Lemma 2.1]{AB:08}, note that in \cite{AB:08} lengths are always signed. Since we 
consider only timelike geodesic triangles (and unsigned lengths), $(a,b,c)$ corresponds to $(-a,-b,-c)$ in \cite{AB:08}. The result can
then immediately be read off from \cite[Lemma 2.1, points 2.\ and 3.]{AB:08}.
\end{pr}
Following \cite{AB:08}, a triple $(a,b,c)$ as in the assumptions of Lemma \ref{lem-tl-tri-mk} will be said to {\em satisfy
timelike size bounds} for $K$.

\subsection{Timelike curvature bounds}

To concisely formulate our notion of timelike curvature bounds in Lorentz\-ian pre-length spaces, we
introduce the following terminology: Let $(x,y,z)$ be a timelike geodesic triangle in a Lorentz\-ian
pre-length space as in Definition \ref{def_tl_geo_tri}, realized by maximal causal curves $\alpha, \beta, \gamma$, 
and suppose that $(\bar{x},\bar{y},\bar{z})$
is a timelike geodesic triangle in a model space $M_K$ with identical side lengths realized by (necessarily)
timelike geodesics $\bar\alpha$, $\bar{\beta}$, $\bar{\gamma}$. Denote the time separation function
in $M_K$ by $\bar{\tau}$. We say that a point $q$ on $\alpha$ \emph{corresponds} to a point $\bar q$ on $\bar{\alpha}$
if $\tau(x,q)=\bar{\tau}(\bar x, \bar q)$, and analogously for $\beta$ and $\gamma$. By Remark \ref{triangle_rem}, under the assumptions made in the following Definition,
any intermediate value of $\tau$ along $\alpha, \beta, \gamma$ is actually attained.

\begin{defi}\label{def:cbb}
A Lorentzian pre-length space \Xll has timelike curvature bounded below (above) by $K\in\R$ if every point in $X$ 
possesses a neighborhood $U$ such that:
\begin{enumerate}[label=(\roman*)]
\item $\tau|_{U\times U}$ is finite and continuous.
\item \label{def:cbb-max-cc} Whenever $x$, $y \in U$ with $x\ll y$, there exists a future-directed causal curve $\alpha$ in 
$U$ with $L_\tau(\alpha) = \tau(x,y)$.
\item Let $(x,y,z)$ be a timelike geodesic triangle in $U$, realized by maximal causal curves $\alpha, \beta, \gamma$
whose side lengths satisfy timelike size bounds for $K$, and let $(\bar{x},\bar{y},\bar{z})$ be a comparison triangle of 
$(x,y,z)$ 
in $M_K$ as given by Lemma \ref{lem-tl-tri-mk}, realized by timelike geodesics $\bar\alpha$, $\bar{\beta}$, $\bar{\gamma}$.
Then whenever $p$, $q$ are points on the sides of $(x,y,z)$ and $\bar p$, $\bar q$ are corresponding
points of the sides of $(\bar{x},\bar{y},\bar{z})$, 
we have $\tau(p,q)\le \bar{\tau}(\bar p, \bar q)$ $($respectively $\tau(p,q)\ge \bar{\tau}(\bar p, \bar q))$.
\end{enumerate}

Such a neighborhood $U$ is called \emph{comparison neighborhood with respect to $M_K$}.
\end{defi}

\begin{figure}[h!]
  \begin{center}
\includegraphics[width=80mm, height= 60mm]{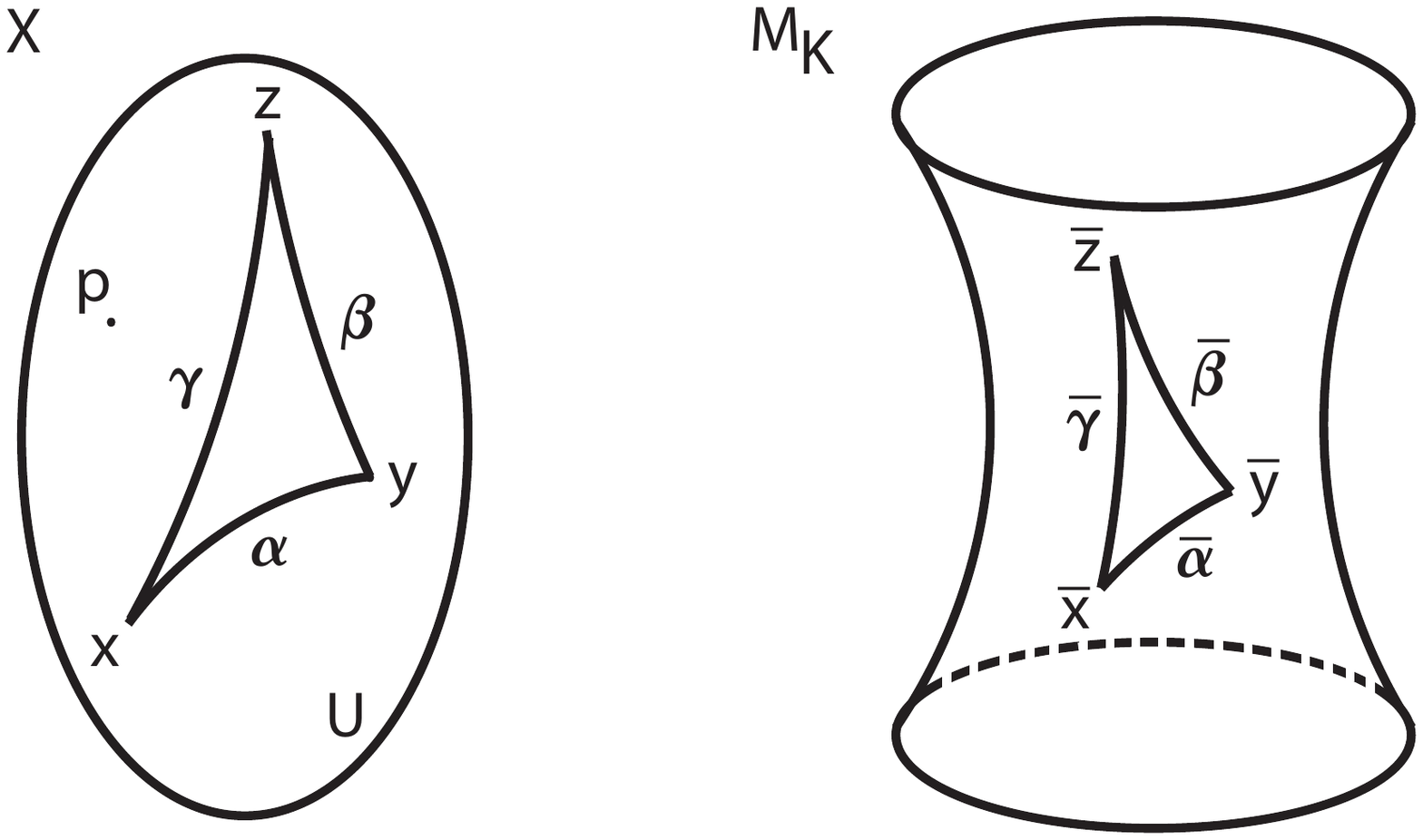}
\end{center}
\caption{Timelike triangle in $X$ and comparison triangle in $M_K$.}
\end{figure}
\vskip-2em

\begin{rem} \ 
\begin{enumerate}[label=(\roman*)]
\item The above definition is as close as possible to the corresponding definition of curvature bounds
in its metric analogue, the theory of CAT$(k)$- respectively Alexandrov spaces (\cite[Def.\ 4.1.9, Def.\ 9.1.1]{BBI:01}).
\item Condition 1 of Definition \ref{def:cbb} in particular secures that $\tau(x,x)=0$ for every $x\in X$ via 
Proposition \ref{prop-tau-0-inf}, so Remark \ref{triangle_rem} applies.
\end{enumerate}
\end{rem}

\begin{ex}\label{AB_rem} In \cite{AB:08}, sectional curvature bounds for general semi-Riemannian manifolds were introduced.
A smooth Lorentzian manifold $M$ is defined to satisfy a lower sectional curvature bound $R\ge K$ if spacelike
sectional curvatures are $\ge K$ and timelike sectional curvatures are $\le K$ (and $R\le K$ with ``timelike''
and ``spacelike'' reversed). It then follows from \cite[Prop.\ 5.2]{AB:08} that $R\ge K$ (respectively $R\le K$) in this sense implies that 
$M$ has timelike curvature bounded below (respectively above) by $K$ in the sense of Definition \ref{def:cbb}. Hence a smooth strongly causal Lorentzian manifold with $R\ge K$ in the sense of \cite{AB:08},
while having timelike sectional curvature bounded \emph{above} by $K$ has timelike sectional 
curvature bounded \emph{below} by $K$ in the sense of Definition \ref{def:cbb}, and analogously
for $R\le K$.
\end{ex}

\subsection{Branching of maximal curves}
\begin{defi}\label{branch_def} (Definition of a branching point)
 Let \Xll be a \LpLS and let $\gamma\colon[a,b]\rightarrow X$ be a maximal causal curve. A point $x:=\gamma(t)$ with $t\in (a,b)$ is 
called \emph{branching point} of $\gamma$ if there exist maximal causal curves $\alpha,\beta\colon[a,c]\rightarrow X$ with $c\ge b$ 
such that $\alpha\rvert_{[a,t]} = \beta\rvert_{[a,t]} = \gamma\rvert_{[a,t]}$ and $\alpha([t,c])\cap \beta([t,c]) = \{x\}$.
If $\alpha, \beta, \gamma$ are timelike then $x$ is called a timelike branching point.
\end{defi}
\begin{ex}\label{ex-funnel-branch} \ 
\begin{enumerate}[label=(\roman*)]
\item \label{ex-funnel-branch1} In a causal/timelike funnel (see Examples \ref{ex-funnel} and \ref{ex-lls}), 
any maximal causal curve from $J^-(p)$ to $J^+(q)$ (resp.\ from $I^-(p)$ to $I^+(p)$) has $q$ as a branching point.
\item For an example of branching in the setting of spacetimes with continuous Lorentzian metrics, see Corollary \ref{cg_branch} below. 
\end{enumerate}
\end{ex} 
In preparation for the following result, call a Lorentzian pre-length space timelike locally uniquely geodesic (l.u.g.)
if every point $x\in X$ has a neighborhood such that, if $p\ll q$ and $p,q\in U$ then there exists a unique
maximal future-directed causal curve from $p$ to $q$ in $U$. Hence compared to \ref{def:cbb-max-cc} of \ref{def:cbb}
one additionally requires uniqueness. Already for low regularity Lorentzian metrics, timelike l.u.g.\ and non-branching
are independent properties. In fact, the classical paper \cite{HW:51} contains examples of $\mathcal{C}^1$-Riemannian metrics that are 
locally uniquely geodesic but display branching, as well as non-branching metrics that fail to be locally uniquely geodesic.
These examples can be translated into the Lorentzian setting, cf.\ \cite{SS:18}. 

Contrary to the case of metric spaces, in the Lorentzian setting the fact that
the time separation function satisfies the reverse triangle inequality precludes
a direct way of generating non-degenerate triangles (i.e., such that 
the strict triangle inequality holds for $\tau$ on their vertices), as required in the standard
proof of non-branching under lower curvature bounds (cf.\ \cite[Lemma 2.4]{Shi:93}).
Conditions (i) and (ii) of the following theorem are sufficient to exclude degeneracy
of comparison triangles.

\begin{thm}\label{nonbranch_th} 
 Let \Xll be a strongly causal Lo\-ren\-tzian pre-length space with timelike curvature bounded below by some $K\in\R$ such that  either
\begin{enumerate}[label=(\roman*)]
 \item Any point in $X$ has a relatively compact, causally closed
neighborhood $\Omega$ such that for any $p\ll q$ in $\Omega$ there is a maximal future-directed
timelike curve from $p$ to $q$ in $\Omega$ that is strictly longer than any future-directed causal curve
from $p$ to $q$ in $\Omega$ that contains a null segment, and there is some $C>0$ bounding the $d$-length of any causal curve
in $\Omega$ (cf.\ (i) and (iv) from Definition \ref{def-loc-LpLS}),  or
\item $X$ is timelike locally uniquely geodesic.
\end{enumerate}
Then maximal timelike curves in $X$ do not have timelike branching points.
\end{thm}
\begin{pr}
 Assume there is a (without loss of generality future-directed) maximal timelike curve $\lambda\colon[a,b]\rightarrow X$ 
that has a timelike branching point $x=\lambda(t_0)$ ($t_0\in(a,b)$). Then there are future-directed timelike maximal curves 
$\alpha,\beta\colon[a,c]\rightarrow X$ with $c>b$ such that $\alpha(t_0) = x = \beta(t_0)$, $\alpha\rvert_{[a,t_0]} = 
\beta\rvert_{[a,t_0]} = \lambda\rvert_{[a,t_0]}$ and $\alpha([t_0,c])\cap \beta([t_0,c]) = \{x\}$. Let $U$ be an open 
comparison neighborhood of $x$ with respect to $M_K$ and let $\Omega$ be a neighborhood of
$x$ as provided by either (i) or (ii). Let $V\subseteq U\cap \Omega$ be an open neighborhood of $x$ such that all 
causal curves with endpoints in $V$ are contained in $U\cap \Omega$ (cf.\ Lemma \ref{lem-cc} \ref{lem-cc-str-cau}).

Our first aim is to show that under any of the assumptions (i) or (ii) we can construct a non-degenerate timelike triangle $(p,q,r)$ in $U\cap\Omega$, i.e.\ with $\tau(p,q)+\tau(q,r) < \tau(p,r)$.

Let $t'\in[a,t_0)$ such that $p:=\lambda(t')\in V$. Choose $s\in(t_0,c]$ such that $\beta([t_0,s])\subseteq V$ and such that there exists some $\tilde s\in (t_0,c]$ with $\tau(p,\beta(s))=\tau(p,\alpha(\tilde s))$, and set $r:=\beta(s)$ and $r':=\alpha(\tilde s)$.

Assuming (ii), note that $I^-(r)\cap V$ is an open neighborhood of $x$, thus $\alpha^{-1}(I^-(r)\cap V)$ is an open neighborhood of $t_0$ in $[a,c]$, so there is an $s'\in [a,c]$, $s'>t_0$ with 
$q:=\alpha(s')\in I^-(r)\cap V$. By our choice of $V$, we have that $\alpha([t_0,s'])\subseteq U\cap \Omega$.  
Consequently, we obtain $p\ll q\ll \beta(s)=r$ and 
there is a unique future-directed maximal causal curve $\gamma$ of positive length
from $q$ to $r$ in $U\cap\Omega$ by Definition \ref{def:cbb},\ref{def:cbb-max-cc}. This 
gives a timelike geodesic triangle $(p,q,r)$ contained in $U\cap\Omega$. Moreover,
by the timelike l.u.g.-property of $\Omega$, the sidelengths of $(p,q,r)$ 
satisfy the strict triangle inequality.

\begin{figure}[h!]
\begin{center}
\includegraphics[width=72mm, height= 48mm]{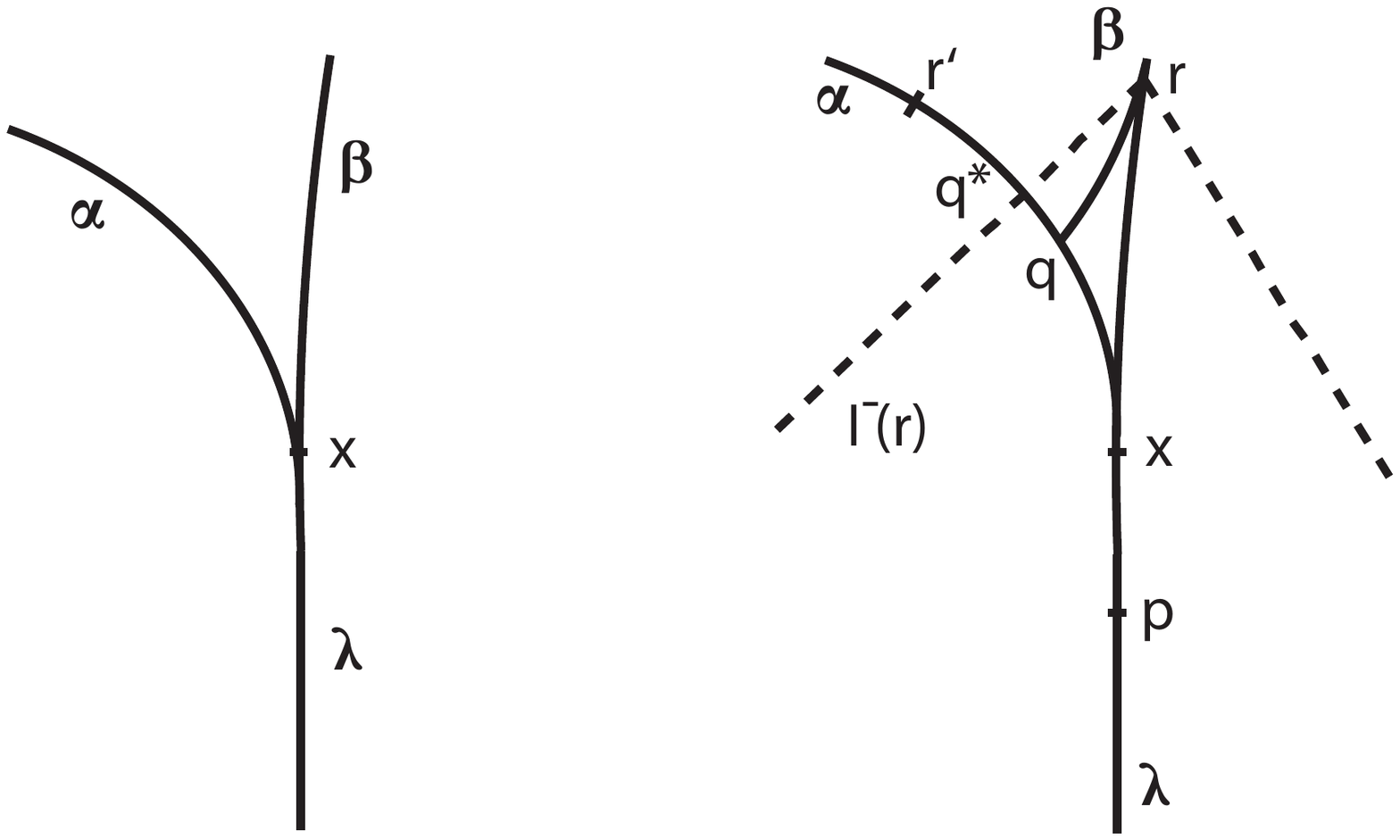}
\end{center}
\caption{Timelike branching.}
\end{figure}

Alternatively assuming (i), we first show $\tau(r',r)=0$. If $\tau(r',r)>0$ we get $\tau(p,r')<\tau(p,r') + \tau(r',r) \leq \tau(p,r) = \tau(p,r')$ --- a contradiction. Defining $s^*:=\sup\{t\in[t',c]: \tau(\alpha(t),r)>0\}$, we obtain $s^*\leq \tilde s$ and thus $q^*:=\alpha(s^*)\leq r'$ and $\tau(q^*,r)=0$. We distinguish two cases: First, let $q^*\ll r'$. Then, $\tau(p,q^*) < \tau(p,r') = \tau(p,r)$ and we pick $\eps>0$ such that $\tau(p,q^*) + \eps < \tau(p,r)$. By continuity, we can choose $q\ll q^*$ on $\alpha$ with $\tau(q,r)<\eps$, and hence
$\tau(p,q) + \tau(q,r) < \tau(p,q^*) + \eps < \tau(p,r)$. This gives the non-degenerate timelike triangle. In the second case, where $q^*=r'$ we derive a contradiction as follows. Picking a sequence $t_n\nearrow s^*$, by (i) there exists a sequence $\sigma_n$ of future-directed timelike curves in $U\cap \Omega$ from
$\alpha(t_n)$ to $r$. By Theorem \ref{thm-lim-curve}, $\sigma_n$ possesses a limit curve $\sigma$ that is future directed causal and
connects $r'$ to $r$. As $\tau(r',r)=0$, $\sigma$
has to be null. Then, the concatenation of $\alpha|_{[t',s*]}$ and $\gamma$ is a future-directed causal curve from $p$ to $r$ that contains a null segment and so by (i), 
there exists a strictly longer timelike curve $\chi$ from $p$ to $r$. Hence
\[
\tau(p,r') + \tau(r',r) = L_\tau(\alpha|_{[t',s*]}) + L_\tau(\sigma)  < L_\tau(\chi) \le
\tau(p,r) = \tau(p,r')\,,
\]
which is impossible.

Thus under any of the assumptions (i) or (ii) of the theorem we arrive at a non-degenerate timelike triangle, and it is clear that the points in the above constructions can be chosen in such a way that the side lengths of this triangle satisfy timelike curvature bounds.

Let $\bar\Delta:=(\bar{p},\bar{q},\bar{r})$ be a comparison triangle in $M_K$. Denote the sides of  $\bar\Delta$ by 
$\bar{\alpha},\bar{\beta},\bar{\gamma}$ and let $\bar{x}_1$ be a point on $\bar{\alpha}$ and $\bar{x}_2 \neq \bar{r}$ 
a point on $\bar{\beta}$ such that $\bar{\tau}(\bar{p},\bar{x}_1) = \bar{\tau}(\bar{p},\bar{x}_2) = \tau(p,x)>0$. Note that 
since $\bar\Delta$ is a non-degenerate triangle we have that $\bar{x}_1\neq \bar{x}_2$. Moreover, we have 
$\bar{\tau}(\bar{x}_1,\bar{r})<\bar{\tau}(\bar{x}_2,\bar{r})$ since otherwise the broken future-directed 
timelike geodesic going from $\bar{p}$ to $\bar{x}_1$ to $\bar{r}$ would be at least as long as the unbroken 
future-directed timelike geodesic $\bar{\beta}$ from $\bar{p}$ to $\bar{r}$. Finally, since the timelike curvature is bounded 
from below by $K$ we obtain that
\begin{equation*}
 \tau(x,r) = \bar{\tau}(\bar{x}_2,\bar{r}) > \bar{\tau}(\bar{x}_1,\bar{r}) \ge \tau(x,r)\,,
\end{equation*}
a contradiction.
\end{pr}
As an immediate consequence we obtain:
\begin{cor}
Let \Xll be a strongly causal Lorentzian length space with timelike curvature bounded below by some $K\in\R$ that is either regular and locally compact or timelike locally uniquely geodesic. Then maximal timelike curves in $X$ do not have timelike branching points.
\end{cor}

\subsection{Causal curvature bounds}\label{sec-causal-curv-bounds}
As already indicated in Remark \ref{AB_rem}, Alexander and Bishop in \cite{AB:08} introduced 
Alexandrov curvature bounds for smooth semi-Riemannian manifolds. In fact, they achieve a
complete characterization of sectional curvature bounds in terms of triangle comparison
(\cite[Thm.\ 1.1]{AB:08}). In their approach, the sides of any given geodesic triangle
(in a sufficiently small normal neighborhood) are parametrized on the interval $[0,1]$,
and this affine parameter is used for comparing with corresponding triangles in the
model spaces. Moreover, timelike geodesics are assigned negative lengths, and spacelike
geodesics positive lengths. 

When trying to generalize this approach to Lorentzian pre-length spaces, there are 
two basic problems. First, there is no built-in notion of spacelike curve
in our setting. In addition, while in the timelike case a viable substitute for a geodesic
in Lorentzian manifolds is given by the notion of maximal causal curve between timelike related
points (as implemented in Definition \ref{def:cbb}), even for maximal null curves there
is no natural parametrization. Indeed, the affine parametrizations on $[0,1]$ employed in 
the comparison results of \cite{AB:08} rely on the fact that geodesics satisfy a
second order ODE, which already for Lorentzian metrics on spacetimes of regularity 
below $C^1$ has no classical counterpart, hence is also unavailable in our setting.

Nevertheless, a restricted notion of causal curvature bounds can also be implemented
for Lorentzian pre-length spaces. In fact, in addition to timelike geodesic triangles in the sense
of Definition \ref{def_tl_geo_tri} we may consider triangles $(x,y,z)$ that satisfy 
$x\ll y \le z$ or $x \le y \ll z$ such that 
$\tau(x,z)<\infty$ and such that the sides (if non-trivial) are realized by future-directed causal curves.
Such triangles will be called {\em admissible causal geodesic triangles}. Those sides
of any such triangle whose vertices are timelike related are called the timelike
sides of the triangle. Since one of the sides of the triangle may have vanishing $\tau$-length,
when realizing such sides by maximal curves in $X$, or the comparison triangle by causal geodesics
in the model spaces, we now allow for the degenerate cases where the realizing curve
(either in $X$ or in the model space) is in fact constant.
By \cite[Lemma 2.1]{AB:08}, the realizability Lemma \ref{lem-tl-tri-mk} only needs minor modifications
to also cover the current setup: Let $(a,b,c)\in \R_{\ge 0}$, with $c\ge a+b$ and at most one entry equal to $0$.
Then the same bounds as is Lemma \ref{lem-tl-tri-mk} guarantee the existence of causal comparison
triangles in the model spaces.
With these conventions, the analogue of Definition \ref{def:cbb} reads:

\begin{defi}\label{def:cbb_causal}
A Lorentzian pre-length space \Xll has causal curvature bounded below (above) by $K\in\R$ if every point in $X$ 
possesses a neighborhood $U$ such that:
\begin{enumerate}[label=(\roman*)] 
\item $\tau|_{U\times U}$ is finite and continuous.
\item Whenever $x,y \in U$ with $x < y$, there exists a causal curve $\alpha$ in $U$ with $L_\tau(\alpha) = \tau(x,y)$.
\item Let $(x,y,z)$ be an admissible causal geodesic triangle in $U$, realized by maximal causal curves 
(or a constant curve, respectively) $\alpha, \beta, \gamma$
whose side lengths satisfy timelike size bounds for $K$, and let $(\bar{x},\bar{y},\bar{z})$ be a comparison triangle of 
$(x,y,z)$ in $M_K$ realized by causal geodesics (or a constant curve, respectively) $\bar\alpha$, $\bar{\beta}$, $\bar{\gamma}$.
Then whenever $p$, $q$ are points on the timelike sides of $(x,y,z)$ and $\bar p$, $\bar q$ are corresponding
points of the timelike sides of $(\bar{x},\bar{y},\bar{z})$,
we have $\tau(p,q)\le \bar{\tau}(\bar p, \bar q)$ $($respectively $\tau(p,q)\ge \bar{\tau}(\bar p, \bar q))$.
\end{enumerate}
Again such a neighborhood $U$ is called a \emph{comparison neighborhood with respect to $M_K$}.
\end{defi} 

Since as explained above there is no natural parametrization for the null side of an admissible
causal geodesic triangle, there is also no natural notion of corresponding points for these
sides and the null side of the comparison triangle in the model space. Thus the restriction
to the timelike sides in the above definition.

Despite this limitation, causal curvature bounds make it possible to establish properties of maximal 
curves and of length-increasing push-up that are 
closely analogous to those of regularly localizable Lorentzian length spaces (cf.\
Theorems \ref{thm-max-cc-char} and \ref{thm-length-increasing-push-up}). In the formulation
of the following result we will call a Lorentzian pre-length space $X$ \emph{locally timelike geodesically connected} if every point in $X$ has a neighborhood $U$ such that for any
points $x\ll y$ in $U$ there exists a future directed timelike maximal curve from $x$
to $y$ that is contained in $U$.
%To show these, we start with the following observation:

\begin{prop}\label{push_up_curvature} Let \Xll be a strongly causal Lorentzian pre-length space with causal curvature bounded above. Moreover, suppose that $X$ is locally causally closed and 
locally timelike geodesically connected. 
Let $\gamma: [a,b]\to X$ be a future-directed causal curve with $\gamma(a)\ll \gamma(b)$ 
and suppose there exists some (non-trivial) sub-interval $[c,d]$ of $[a,b]$ 
such that $\gamma|_{[c,d]}$ is null. Then $\gamma$ is not maximal.
\end{prop}
\begin{pr} 
It suffices to consider the case where $a<c$ and $d=b$. So suppose, to the contrary, that $\gamma$ were maximal and 
let $t_0:= \inf \{t \in [a,b] : L_\tau(\gamma|_{[t,b]}) = 0\}$. By Proposition \ref{prop-max-prop} \ref{prop-max-sub-int}
and the lower semi-continuity of $\tau$, $L_\tau(\gamma|_{[t_0,b]})=\tau(\gamma(t_0),\gamma(b))=0$. 
Also, $a<t_0 (\le c)$ since otherwise $\tau(\gamma(a),\gamma(b))$ would
vanish, contrary to our assumption. Moreover, for any $a\le s <t_0$, $\tau(\gamma(s),\gamma(t_0))=
L_\tau(\gamma|_{[s,t_0]})>0$ by definition of $t_0$, and $\gamma(s)\ll \gamma(b)$.

\begin{figure}[h!]
  \begin{center}
\includegraphics[width=80mm, height= 60mm]{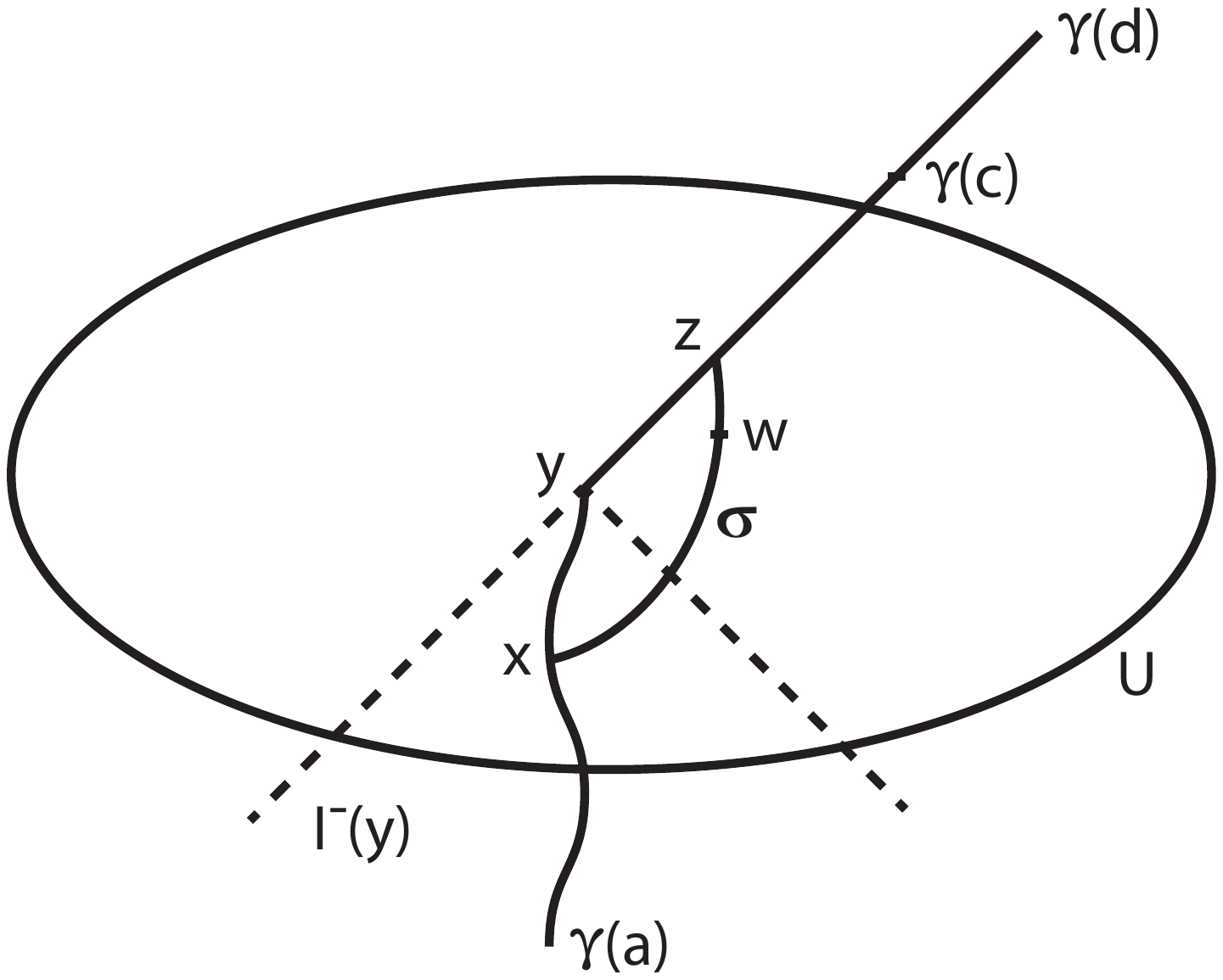}
\end{center}
\caption{Push-up via upper causal curvature bound.}
\end{figure}

Let $V_1$ be a comparison neighborhood of $y:=\gamma(t_0)$ and $V_2, V_3$ be neighborhoods of $y$ that are causally closed and timelike geodesically connected, respectively.
Since $X$ is strongly causal, there exists an open causally convex neighborhood 
$U\subseteq V_1\cap V_2\cap V_3$ of $y$. Then $U$ is still a comparison neighborhood and at the same time is both causally closed and timelike geodesically connected.

Pick some $x:=\gamma(s)$ and $z:=\gamma(t)$ for $s<t_0<t$ and $s,t$ sufficiently close to $t_0$
so that both points are contained in $U$.
Again by Proposition \ref{prop-max-prop} \ref{prop-max-sub-int}, 
$\tau(x,z) = L_\tau(\gamma_{[s,t]})=L_\tau(\gamma_{[s,t_0]})=\tau(x,y)>0$, and $\tau(y,z) = L_\tau(\gamma_{[t_0,t]})=0$.
Also, we may choose $s$ and $t$ such that the side lengths satisfy timelike size bounds for the upper bound $K$ on the curvature.
Then the corresponding comparison triangle $(\bar x, \bar y, \bar z)$ in the model space $M_K$ is degenerate,
with the timelike sides coinciding and the null side collapsing to a single point $\bar y = \bar z$. 

Since $U$ is timelike geodesically connected there exists a future directed maximal timelike
curve $\sigma$ from $x$ to $z$ in $U$ (recall that $\tau(x,z)>0$, so $x\ll z$). Let $w_k$
($k\in \N$) be a sequence of points on $\sigma$ with $w_k\ll z$ and $w_k\to z$. Then for 
$k$ sufficiently large we must have $w_k \not\in I^-(y)$. Indeed, otherwise since $U$ is
causally closed we would obtain from  $w_k\ll y < z$ that $z \le y < z$, contradicting
strong causality ($y$ and $z$ cannot be separated by chronological diamonds).

Picking $k\in \N$ such that $w_k\not\in I^-(y)$ we arrive at a contradiction because 
our assumption on the
upper curvature bound implies $0< \bar \tau(\overline{w_k},\bar y) \le \tau(w_k,y)$, i.e.,
$w_k\ll y$.
\end{pr} 
\begin{rem}\label{rem-curvbound-vs-localizable}
Combining Proposition \ref{push_up_curvature} with the general assumptions on comparison neighborhoods, we obtain the following: 
In a strongly causal Lorentzian pre-length space with causal curvature bounded above 
that additionally  is locally causally closed and locally timelike geodesically connected,
every point possesses a 
neighborhood in which property (iv) from Definition \ref{def-loc-LpLS} is satisfied. This allows us to adapt the proofs
of Theorems \ref{thm-max-cc-char} and \ref{thm-length-increasing-push-up} almost verbatim to obtain
the following two results. 
\end{rem}
\begin{thm}\label{thm-max-cc-char_curvature}
Let $X$ be a strongly causal Lorentzian pre-length space with causal curvature bounded above,
that is also locally causally closed and 
locally timelike geodesically connected. Then
maximal causal curves have a causal character, i.e., if for a (future-directed) maximal causal 
curve $\gamma\colon[a,b]\rightarrow X$ there are $a\leq t_1 < t_2\leq b$ with $\gamma(t_1)\ll\gamma(t_2)$, then $\gamma$ 
is timelike. Otherwise it is null.
\end{thm}
Furthermore, also in the current setting the principle of length-increasing push-up of causal curves is valid:
\begin{thm}\label{thm-length-increasing-push-up_curvature}
Let \Xll be as in Theorem \ref{thm-max-cc-char_curvature}
and let $\gamma: [a,b]\to X$ be a future-directed causal curve. 
\begin{enumerate}[label=(\roman*)]
\item\label{thm-length-increasing-push-up_curvature-1} If $L_\tau(\gamma)>0$ and if $\gamma|_{[c,d]}$ is null on some (non-trivial) 
sub-interval $[c,d]$ of $[a,b]$, then there exists a strictly longer future-directed timelike 
curve $\sigma$ from $\gamma(a)$ to $\gamma(b)$.
\item If there exist $a\le c<d\le b$ such that
$\gamma|_{[c,d]}$ is rectifiable, then there exists timelike future-directed curve from $\gamma(a)$
to $\gamma(b)$.
\end{enumerate} 
\end{thm} 
This result provides an interesting perspective on length-increasing push-up,
namely that it is a necessary consequence of any upper bound on synthetic causal curvature,
while in the smooth setting it is usually traced back to the Gauss Lemma.
As already pointed out in Remark \ref{rem-curvbound-vs-localizable}, 
there is also a close connection to regular localizability in the sense of Definition \ref{def-loc-LpLS},
cf.\ also Example \ref{ex_reg_lls} \ref{ex_reg_lls-2} below.

Another consequence of the observation in Remark \ref{rem-curvbound-vs-localizable} is the following corollary of Theorem \ref{nonbranch_th}.
\begin{cor}
Let \Xll be a strongly causal $d$-compatible Lo\-ren\-tzian pre-length space such 
that any point has a relatively compact neighborhood which is causally closed and
timelike geodesically connected. If $X$ has timelike curvature bounded below and causal curvature bounded above, then maximal timelike curves in $X$ do not have timelike branching points.
\end{cor}

\subsection{Curvature singularities}
The synthetic approach to curvature bounds developed in the previous sections in particular 
allows one to define curvature singularities in Lorentzian pre-length spaces, and thereby in 
particular in spacetimes of low regularity, where a classical description in terms of
the Riemann curvature tensor is not viable, or even in settings where there is no spacetime metric
available at all. In the present section we introduce the
necessary notions. 
\begin{defi}\label{def-curv-singularity} A Lorentzian pre-length space \Xll has timelike (respectively causal) 
curvature unbounded below/above
if some point in $X$ possesses a neighborhood $U$ such that 1.\ and 2.\ from Definition \ref{def:cbb}
(respectively Definition \ref{def:cbb_causal}) are satisfied, but such that the corresponding part of 
property 3 from these definitions fails to hold for any $K\in \R$. In this case we say that
$X$ has a curvature singularity. 
\end{defi}

Thus we assume that locally there always exist timelike triangles, but that
the comparison conditions fail somewhere in $X$. 

\begin{ex} Consider a causal or timelike funnel $X$ as in Example \ref{ex-funnel} with $\lambda$ timelike. Since $X$ is clearly 
timelike uniquely geodesic,  Example \ref{ex-funnel-branch} \ref{ex-funnel-branch1} 
and Theorem \ref{nonbranch_th} imply that $X$ has timelike curvature unbounded below.
Moreover, if $\lambda$ is null then one easily constructs maximal curves violating Theorem
\ref{thm-length-increasing-push-up_curvature} \ref{thm-length-increasing-push-up_curvature-1}, 
so in this case the causal curvature of $X$ is unbounded above. 
\end{ex}

Curvature singularities are of central 
importance in General Relativity. As an application of the notions introduced above, we therefore demonstrate that the 
central singularity of the interior Schwarzschild metric can be detected by timelike 
triangle comparison. 
\begin{ex}
The curvature singularity of the interior of a Schwarzschild black hole.

Consider the interior Schwarzschild metric (cf., e.g., \cite{ONe:83,Wal:84,DLC:08})
\begin{equation}\label{eq-schwarzschild-metric}
ds^2 = -\left(1-\frac{2M}{r} \right)dt^2 +\left(1-\frac{2M}{r}  \right)^{-1} dr^2 + r^2(d\theta^2 + \sin^2\theta d\phi^2),
\end{equation}
where $M>0$, $t\in \R$, $r\in (0,2M)$ and $\theta$, $\phi$ parametrize the two-sphere $S^2$. 
The metric \eqref{eq-schwarzschild-metric} has the form of a warped product, with fiber $S^2$
and leafs isometric to $\R\times (0,2M)$ with metric
\begin{equation}\label{eq-schwarzschild-leaf}
g = -\left(\frac{2M}{r} - 1 \right)^{-1} dr^2 +\left(\frac{2M}{r} -1 \right)dt^2.
\end{equation}
We write $g$ in this form to emphasize that, in the Schwarzschild interior, the coordinate
$r$ is timelike, and $t$ is spacelike.
As is well-known, for $r\to 0+$, this metric incurs a curvature singularity. In fact, 
its sectional curvature (i.e., its Gauss curvature), is given by $K=\frac{2M}{r^3}$ (\cite[Lemma 13.3]{ONe:83}),
hence goes to $+\infty$ as one approaches the spacelike hypersurface $r=0$.
 
Although \cite[Thm.\ 1.1]{AB:08} provides a characterization of sectional curvature in terms of triangle
comparison, we cannot directly from that result conclude that the spacetime \eqref{eq-schwarzschild-leaf}
displays a timelike curvature singularity in the sense of Definition \ref{def-curv-singularity}. In fact, 
divergence of timelike sectional curvature (as is the case here) only implies that triangle
comparison for triangles of arbitrary causal character will fail, and does not necessarily
entail that {\em timelike} triangles will be the culprits for this behavior.
Instead, to verify the conditions of Definition \ref{def-curv-singularity}
we will explicitly study a family of timelike geodesic triangles approaching $r=0$.

For brevity, put $h(r):=1-\frac{2M}{r}$. Then denoting by $\tau$ proper time along a timelike geodesic,
by \cite[Prop.\ 13.11]{ONe:83} in the case $L=0$, for the constant of motion $E$ we have
\begin{align}
E &= h(r)\frac{dt}{d\tau}  \\
E^2 &= \left(\frac{dr}{d\tau}\right)^2 + h(r).
\end{align}
By \cite[Eq.\ (81)]{DLC:08}, for $E=1$ there are two families of pregeodesics $\gamma_\pm(r)=(r,t_\pm(r)+ \mathrm{const})$, where
\begin{equation}\label{eq-tpm}
t_\pm(r) = \pm \frac{2}{3}(6M+r)\sqrt{\frac{r}{2M}} \mp 4M \mathrm{artanh}\left(\sqrt{\frac{r}{2M}}\right).
\end{equation} 
Here, the corresponding proper time is given by (see \cite[Eq.\ (76)]{DLC:08})
\begin{equation}\label{eq-prop-time-E1}
\tau(r) = \pm \sqrt{\frac{2r}{M}} \frac{r}{3} + \mathrm{const}.
\end{equation}
Also, for $E=0$ there is a pregeodesic of the form $\gamma_0(r)=(r,0)$, with 
\begin{equation}\label{eq-prop-time-E0}
\frac{dr}{d\tau}=\sqrt{\frac{2M}{r}-1}.
\end{equation}
Now fix the constant in $\gamma_-$ to be $-2C$ for some $C>0$ to be specified later, and denote, for $k\in \N$, 
by $\gamma^{(k)}_+$ the pregeodesic
$\gamma_+$ with constant $\frac{C}{k}$. Let $x$ be the intersection of $\gamma_0$ and $\gamma_-$,
$y_k$ that of $\gamma_-$ and $\gamma^{(k)}_+$, and $z_k$ that of $\gamma^{(k)}_+$ and $\gamma_0$
(see Figure \ref{fig_geo}).
As the time-orientation in the Schwarzschild interior is directed towards $r=0$, $(x,y_k,z_k)$ is a timelike
triangle with the corresponding pregeodesics from above as realizing sides. 

\begin{figure}[h!] 
  \begin{center}
\includegraphics[width=120mm, height= 60mm]{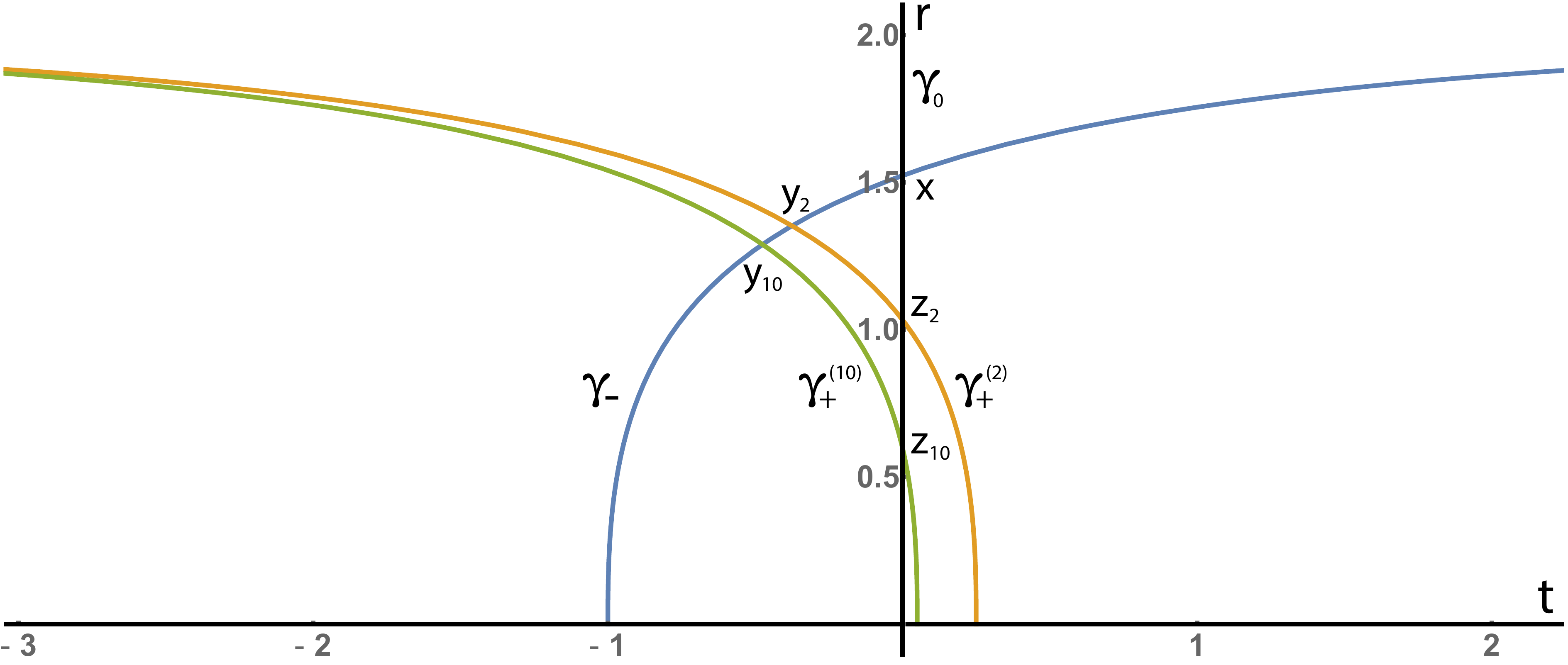} 
\end{center}
\caption{Infalling timelike geodesics for $M=1$ and $C=0.5$.}\label{fig_geo}
\end{figure}

We will use the timelike triangles $(x,y_k,z_k)$ to demonstrate that the interior Schwarzschild solution has
a curvature singularity in the sense of Definition \ref{def-curv-singularity}, more precisely that the
timelike curvature is unbounded below. 

Suppose, to the contrary, that the timelike curvature were bounded below by some $K\in \R$. It follows
from \eqref{eq-prop-time-E1} and \eqref{eq-prop-time-E0} that each of the above geodesics reaches
$r=0$ in finite proper time and that these finite times go to zero as $C\to 0+$. Thus we can choose $C>0$ so small that each triangle satisfies timelike size bounds for $K$.

Consider now the scalar product of the unit tangent vectors of $\gamma_0$ and $\gamma^{(k)}_+$ at
$z_k=:(0,r_k)$. A straightforward calculation shows that this is given by
\[
 -\frac{1}{\sqrt{1-(t_+')^2(r_k)\left(\frac{2M}{r_k}-1\right)^2}} \to -1 \quad (k\to \infty),
\]
confirming that the triangles become degenerate, with the hyperbolic angle at $z_k$ collapsing
to $0$ as $k\to \infty$.

One can now apply \cite[Prop.\ 5.1]{AB:08} to conclude that the comparison triangles
in the model space must display the same behavior. In fact, as in \cite{AB:08} 
denote by $\gamma_{pq}$ (respectively $\bar\gamma_{\bar p \bar q}$) the maximal timelike
geodesic connecting $p$ to $q\gg p$ in the spacetime (respectively $\bar p$ to $\bar q\gg \bar p$ in $M_K$),
parametrized on $[0,1]$. 
Moreover, let $E_q(r):=\langle \gamma_{qr}'(0),\gamma_{qr}'(0)\rangle$. Then by \cite[Eq.\ (5.3)]{AB:08}, for a timelike geodesic triangle $(p,q,r)$,
\begin{equation}\label{eq-AB-angle}
(E_q\circ \gamma_{pr})'(0) = -2 \langle \gamma_{pq}'(0),\gamma_{pr}'(0)\rangle,
\end{equation}
and analogously for $(\bar p,\bar q,\bar r)$ in $M_K$. The key observation now is that to calculate the left hand side of this equation, information is only required about $(E_q\circ \gamma_{pr})(t)$
for $t$ arbitrarily close to $0$ and that, since $p\in I^-(q)$, so is $\gamma_{pr}(t)$ for $t$ small. 
Hence knowledge of timelike distances suffices to determine $(E_q\circ \gamma_{pr})(t)$. 
More precisely, $-(E_q\circ \gamma_{pr})(t)$ equals the square of the time separation from 
$\gamma_{pr}(t)$ to $q$. Consequently, our assumption on timelike triangle comparison implies 
that $(E_q\circ \gamma_{pr})'(0)\ge (\bar E_{\bar q}\circ \bar\gamma_{\bar p \bar r})'(0)$, hence
by \eqref{eq-AB-angle} that $\langle \gamma_{pq}'(0),\gamma_{pr}'(0)\rangle \le \langle \bar\gamma_{\bar p \bar q}'(0),\gamma_{\bar p \bar r}'(0)\rangle$. Since the sidelengths of $(p,q,r)$ and $(\bar p,\bar q,\bar r)$
coincide by definition, the same relation must hold between the scalar products of the 
unit tangent vectors. 

Applying this to the timelike triangles $(x,y_k,z_k)$ and $(\bar x,\bar y_k,\bar z_k)$, 
it follows that the assumption
of the lower bound $K$ on the timelike curvature in the sense of Definition \ref{def:cbb}
implies that the hyperbolic angle at $\bar z_k$ must also go to $0$ as $k\to \infty$. 
But in the model space
$M_K$, $(\bar x,\bar y_k,\bar z_k)$ converges to a nondegenerate timelike triangle
whose side lengths equal the (non-trivial) limits of the sidelengths of $(x,y_k,z_k)$ (hence satisfy
the strict reverse triangle inequality), a contradiction.
\end{ex}

\section{Classes of examples}\label{sec-app}
\subsection{Continuous Lorentzian metrics} \label{subsec-cont_lor}
As a main application of the theory developed so far, in this subsection we are going to show that
any smooth manifold endowed with a continuous, strongly causal, and causally plain (as defined in \cite{CG:12}, see below) 
Lorentzian metric provides a natural example of a Lorentzian length space.

Let $M$ be a smooth manifold and let $g$ be a continuous Lorentzian metric on $M$ (cf.\ Example \ref{ex-cau-spa-con}).
By a causal respectively timelike curve in $M$ we mean a locally Lipschitz curve whose tangent vector is causal respectively 
timelike almost everywhere. It would also be possible to start from absolutely continuous curves, but since
causal absolutely continuous curves always possess a reparametrization that is Lipschitz (cf.\ the discussion
in \cite[Sec.\ 2.1, Rem.\ 2]{Min:17}), the above convention is not a restriction.

The time separation function $\tau: M\times M \to [0,\infty]$ is defined in the standard way, i.e., 
$\tau(x,y) = \sup\{L_g(\gamma) : \gamma \text{ future-directed causal from } x \text{ to }y \}$,
if $x\le y$ and $\tau(x,y)=0$ otherwise. Here, by $L_g(\gamma)$ we denote the $g$-length
of a causal curve $\gamma:[a,b]\to M$, i.e., $L_g(\gamma) = \int_a^b \sqrt{-g(\dot \gamma,\dot{\gamma})}\,dt$.
Then \ref{def-LpLS} \eqref{eq-rev-tri}, i.e., the reverse triangle inequality for $\tau$,
just as in the smooth case, follows directly from the definition.
Also, we fix any complete Riemannian metric $h$ on $M$ and denote by $d^h$ the metric introduced by $h$.

\begin{rem}\label{max_curve_rem} Any $L_g$-maximal curve $\gamma$ is also $L_\tau$-maximal, and $L_g(\gamma)=L_\tau(\gamma)$.
In fact, suppose that $\gamma$ is future-directed causal
from $p$ to $q$ with $L_g(\gamma)=\tau(p,q)$. Since $\tau(p,q)\ge L_\tau(\sigma)$ for any
future-directed causal curve $\sigma$ from $p$ to $q$, $\gamma$ is also $L_\tau$-maximal.
Moreover, $L_g(\gamma)\ge L_\tau(\gamma)$
and since the converse inequality always holds (see Lemma \ref{lem-str-caus-tau-Tau} below), $L_g(\gamma)=L_\tau(\gamma)$. 
\end{rem}

\begin{ex}\label{bubble-ex} It was shown in \cite{CG:12} that without further assumptions, causality theory of continuous Lorentzian
metrics displays a number of unexpected (and unwanted) new phenomena. Consider the following metric on $\R^2$ 
(\cite[Ex.\ 1.11]{CG:12}):
\begin{equation}\label{cg_metric}
g = -(du+(1-|u|^\lambda) dx)^2 + dx^2,
\end{equation}
where $\lambda\in (0,1)$. Then $g\in C^{0,\lambda}(\R^2)$ and $g$ is smooth everywhere except on the $x$-axis.
However, the light cone $J(p)\setminus I(p)$ of any point $p$ on the $x$-axis has non-zero measure (and is 
covered by the null curves emanating
from $p$). For points $q$ in the interior of this region (the so-called {\em bubbling region}), push-up fails.
Also, although $\tau(p,q)>0$, there does not exist any timelike curve connecting $p$ to $q$, so $p\not\ll q$. 
  In addition, as already noted in 
\cite{GLS:18}, $\tau$ is not lower semicontinuous for this spacetime: Let $p$ be the origin and 
let $q\in \partial I^+(p)$. Then $\tau(p,q)>0$, but taking $p_n:=(\frac{1}{n},0)$, $\tau(p_n,q)=0$
for every $n$. Consequently (fixing any background Riemannian metric $h$ on $\R^2$),
$(\R^2,d^h,\ll,\le,\tau)$ is not a Lorentzian pre-length space. 
\end{ex}

In order to exhibit additional exotic causality properties of continuous Lorentzian metrics, let
us study the metric $g$ from \eqref{cg_metric} in greater depth. For concreteness, set $\lambda:=\frac{1}{2}$
and let $M:=(-1,1)\times\R$. 
Thus the metric is
\begin{equation*}
 g_{(u,x)} = -du^2 + 2 (\sqrt{|u|}-1)du\, dx + \sqrt{|u|}(2-\sqrt{|u|})dx^2\,,
\end{equation*}
and its inverse is given by
\begin{equation*}
 g_{(u,x)}^{-1}= \sqrt{|u|}(\sqrt{|u|}-2) du^2 + 2 (\sqrt{|u|}-1)du\, dx + dx^2\,.
\end{equation*}

We first collect some basic facts about the causality of $(M,g)$, choosing the time orientation by defining 
$\partial_u$ to be future-directed.

Let $\gamma=(\alpha,\beta)\colon[a,b]\rightarrow M$ be a causal curve, then
\begin{align*}
 0 &\geq -\dot{\alpha}(s)^2 + 2 (\sqrt{|u|}-1)\dot{\alpha}(s)\dot{\beta}(s) + 
\sqrt{|u|}\underbrace{(2-\sqrt{|u|})}_{>0}\dot{\beta}(s)^2\\
&\geq -\dot{\alpha}(s)^2 + 2 (\sqrt{|u|}-1)\dot{\alpha}(s)\dot{\beta}(s)\,,
\end{align*}
for all $s\in[a,b]$ such that $\dot{\gamma}(s)$ exists (i.e., for almost all $s\in[a,b]$).

Furthermore, the time orientation of $\gamma$ gives
\begin{align}
 0 &>  -\dot{\alpha}(s) + (\sqrt{|u|} - 1) \dot{\beta}(s) \text{ (if }\gamma\text{ is future-directed)}\,,\label{eq-fd}\\
 0 &<  -\dot{\alpha}(s) + (\sqrt{|u|} - 1) \dot{\beta}(s) \text{ (if }\gamma\text{ is past-directed)}\,,\label{eq-pd}
\end{align}
again for all $s\in[a,b]$ such that $\dot{\gamma}(s)$ exists.

\begin{lem}\label{lem-str-cau}
 The spacetime $(M,g)$ is strongly causal.
\end{lem}
\begin{pr}
 Define $f\colon M\rightarrow \R$ by $f(u,x):=u$, then the gradient of $f$ is given by $\grad(f)_{(u,x)} = 
\sqrt{|u|}(\sqrt{|u|}-2)\partial_u + (\sqrt{|u|}-1)\partial_x$. Thus
\begin{align*}
 g_{(u,x)}(\grad(f),& \grad(f)) = \\
& -(\sqrt{|u|}(\sqrt{|u|}-2))^2 + \sqrt{|u|}\underbrace{(\sqrt{|u|}-2)}_{<0}(\sqrt{|u|}-1)^2 
\leq 0\,,
\end{align*}
and $g_{(u,x)}(\grad(f),\grad(f))<0$ for $u\neq 0$. Moreover, $g_{(u,x)}(\grad(f),\partial_u) = 1$, so $\grad(f)$ is 
past-directed causal on $M$ and past-directed timelike on $(-1,0)\times\R$ and $(0,1)\times\R$, hence is a temporal 
function there. So strong causality holds on $(-1,0)\times\R$ and $(0,1)\times\R$. 

It remains to show strong causality at points $(0,x_0)$. Let $\gamma=(\alpha,\beta)\colon[a,b]$ $\rightarrow M$ be a 
future-directed causal curve with $\gamma(a)=(0,x_0)$ and $\gamma(b)=(u_1,x_1)$. Note that from $\grad(f)$ being 
past-directed causal follows that $0\leq\dot{\alpha}(s)$, whenever $\dot{\gamma}(s)$ exists. So $\alpha(s)\geq 0 $ for all 
$s\in[a,b]$ and if $u_1>0$ then $\gamma$ cannot return to a neighborhood of $(0,x_0)$ that does not contain $(u_1,x_1)$. 
Finally, if $\gamma$ does not enter $(0,1)\times\R$, then $\alpha=0$, and so by equation \eqref{eq-fd} we obtain 
$\dot{\beta}(s)>0$ almost everywhere. Consequently, $\gamma$ cannot return to a neighborhood not containing $(0,x_1)$ (note 
that necessarily $x_1>x_0$).
\end{pr}

Now fix a point $q:=(u_0,x_0)$ in the (upper right) bubble region, i.e., $0 < u_0 < \min(\frac{x_0}{4},1)$.
\begin{prop}\label{prop-ex-max}
 There exists a maximal causal curve from $0$ to $q$.
\end{prop}
\begin{pr} By Remark \ref{max_curve_rem}, it suffices to show the existence of a future-directed causal curve
from $0$ to $q$ whose $L_g$-length is maximal. To this end, 
we will show that $J^+(0)\cap J^-(q)$ is compact and then refer to results from \cite{Sae:16}. 
By the above we know that $J^+(0)\subseteq [0,1)\times\R$. The past 
lightcone emanating from $q$ is bounded by the null curves $\nu, \mu$ given below. Note that they are pregeodesics on 
$(0,1)\times\R$, since the metric is smooth there and so \cite[Thm.\ 4.13]{BEE:96} applies. The left bounding null 
curve $\nu\colon[0,x_0]\rightarrow M$ is given by
\begin{align*}
 \nu(x):= \begin{cases}
           (\frac{1}{4}(2\sqrt{u_0}-x)^2,\,x_0-x) \quad &(x\in[0,2\sqrt{u_0}])\,,\\
           (0,\,x_0-x)\quad &(x\in[2\sqrt{u_0},x_0])\,.
          \end{cases}
\end{align*}
It connects $q$ to $0$ and is past-directed null. Moreover, note that $\nu$ is parametrized as $x\mapsto (v(x),\,x_0-x)$ 
and so satisfies the equation $\dot{v}(x) = -\sqrt{|v(x)|}$ (cf.\ \cite[Eq.\ (1.20)]{CG:12}). The other 
past-directed null curve emanating from $q$ solves the equation 
\begin{equation*}\label{eq-mu-u}
 \dot{u}(x) = -2 + \sqrt{|u(x)|}\,,
\end{equation*}
if parametrized as $\mu(x)=(u(x),\,x)$. This corresponds to the $\eps=-1$ case in \cite[Eq.\ (1.20)]{CG:12}. As 
initial condition we impose $u(x_0)=u_0$. We can solve this equation as follows: Let $F\colon[-1,1]\rightarrow\R$ be given by
\begin{equation*}
 F(s):=\int_{u_0}^s \frac{1}{-2+\sqrt{|r|}}\,\mathrm{d} r\, + x_0\quad (s\in[-1,1])\,.
\end{equation*}
Then clearly $F'(s)\leq -\frac{1}{2}$, and hence $F$ is strictly monotonically decreasing and 
$F(u_0)=x_0$. Define $x':=F(0) = x_0 - \int^{u_0}_0 \frac{1}{-2+\sqrt{|r|}}\,\mathrm{d} r$, $a:=F(1)$ and $b:=F(-1)$. Then 
$a<x_0<x'<b$ and $F([-1,1])=[a,b]$, so the inverse of $F$ exists, $F^{-1}\colon[a,b]\rightarrow [-1,1]$. Finally, $u$ is 
given as $u:=F^{-1}|_{[x_0,x']}\colon[x_0,x']\rightarrow (-1,1)$ and satisfies $u(x_0)=u_0$ and $u(x')=0$. Thus 
$\mu\colon[x_0,x']\rightarrow M$, $x\mapsto(u(x),x)$ connects $q$ to $(0,x')$. Then as we saw above in Lemma 
\ref{lem-str-cau} any past-directed causal curve from $(0,x')$ to $0$ lies on the $x$-axis.

Now we see that
\begin{align*}
 &\qquad J^+(0)\cap J^-(q) =\\
 &\Big\{(u,x)\in M:\, x\in[0,x_0]\text{ and } 0\leq u \leq v(x)\text { or } x\in[x_0,x']\text{ and 
}0\leq u\leq u(x)\Big\}\,,
\end{align*}
which is a compact subset of $M$.

\begin{figure}[h]
\begin{center}
 \includegraphics[width=300pt, height=150pt]{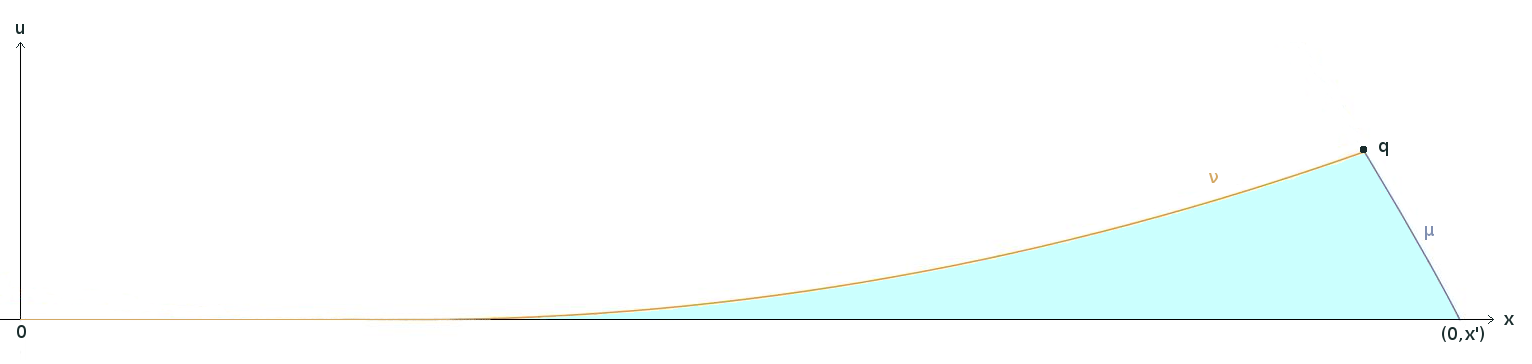}
 \caption{$J^+(0)\cap J^-(q)$ for $q=(\frac{1}{8},1)$}
\end{center}
\end{figure}

Finally we are able conclude that $\tilde{C}(0,q)$ is compact (with respect to the compact-open topology) as in the first 
part of the proof of \cite[Thm.\ 3.2]{Sae:16}. So, \cite[Prop.\ 6.4]{Sae:16} establishes the existence of a 
maximal causal curve connecting $0$ and $q$.
\end{pr}

\begin{cor}\label{cg_branch}
 Let $\gamma=(\alpha,\beta)\colon[0,1]\rightarrow M$ be a maximal curve from $0$ to $q$, as given by Proposition 
\ref{prop-ex-max}, then $\gamma$ has a branching point at which it changes its causal character (from null to
timelike).
\end{cor}
\begin{pr} We first note that since the future-directed causal curve given by going from $0$ along the $x$-axis to 
$(0,x_0)$ and then vertically to $(u_0,x_0)=q$ has length $u_0>0$, we know that $L_g(\gamma)\geq u_0 > 0$. Therefore $\gamma$ 
has to be somewhere timelike and in particular, it has to leave the $x$-axis after $x_0-2\sqrt{u_0}$ (where $\nu$ intersects 
the $x$-axis for the first time) and before $x'$. Off the $x$-axis, i.e. on $(0,1)\times\R$, $\gamma$ has to be pregeodesic, 
since it is maximizing everywhere and the metric is smooth in this region. Consequently, 
$\gamma$ has a causal character on $(0,1)\times\R$. First we observe that $\gamma$ cannot be null on $(0,1)\times\R$, since 
then $L_g(\gamma)$ would vanish, a contradiction. Thus, $\gamma$ has to be timelike on $(0,1)\times\R$. 
This demonstrates that $\gamma$ is a maximal causal curve whose causal character changes. Moreover, the point where
it leaves the $x$-axis is a branching point as defined in Definition \ref{branch_def} because the $x$-axis itself is
maximizing (and null) between any of its points.
\end{pr}
In light of the results of Section \ref{sec-curvature} it can be argued that the root cause of the
phenomena described above lies in the fact that the curvature of the metric $g$ is unbounded near $\{u=0\}$,
see \cite[Eq.\ (1.23)]{CG:12}.

The foregoing considerations demonstrate that Lorentzian metrics that are merely continuous do not 
provide a satisfactory causality theory, a problem clearly recognized already
by Chrusciel and Grant in \cite{CG:12}. 
As a sufficient condition for reasonable causality properties of continuous Lorentzian metrics, these authors 
introduced the notion of a {\em causally plain} Lorentzian metric. To define this concept, let us first
recall from \cite[Def.\ 1.3]{CG:12} that a locally Lipschitz curve $\gamma$ is called locally uniformly
timelike (l.u.t.) if there exists a smooth Lorentzian metric $\check g \prec g$ such that $\check 
g(\dot\gamma,\dot{\gamma})<0$ almost everywhere. For $U\subseteq M$ open and $p\in U$, by $\check I^\pm(p,U)$ we 
denote the set of all points that can be reached by a future directed, respectively past directed l.u.t.\ curve in $U$ starting in $p$. 
Moreover, a {\em cylindrical neighborhood}
of a point $p$ is a relatively compact chart domain containing $p$ such that in this chart $g$ equals
the Minkowski metric at $p$ and the slopes of the light cones of $g$ remain close to $1$ (see
\cite[Def.\ 1.8]{CG:12} for a precise definition). Finally, $(M,g)$ is called causally plain if every
$p\in M$ possesses a cylindrical neighborhood $U$ such that $\partial \check I^\pm(p,U)=\partial J^\pm(p,U)$.
Otherwise it is called bubbling. Causally plain spacetimes in particular satisfy the standard
push-up properties (\cite[1.22--1.24]{CG:12}), as well as $I^\pm(p)=\check I^\pm(p)$ for every $p\in M$.
Any spacetime with a Lipschitz continuous metric is causally plain
(\cite[Cor.\ 1.17]{CG:12}).

Our main aim in this section is to establish that $(M,d^h,\ll,\le,\tau)$ as defined above is a Lorentzian length space for any
continuous, strongly causal and causally plain metric $g$. 

\begin{lem}\label{con_met_tau} Let $(M,g)$ be causally plain and let $p,q\in M$. Then $p\ll q$ if and only if $\tau(p,q)>0$. 
\end{lem}
\begin{pr}
That $p\ll q$ implies $\tau(p,q)>0$ is immediate from the definition of $\tau$.

For the converse implication, suppose that
$\tau(p,q)>0$ and let $\gamma$ be a future-directed causal curve from $p$ to $q$ with
$L_g(\gamma)\ge \frac{1}{2}\tau(p,q)>0$. If we can show that $\gamma$ enters $I^+(p)$ then
$p\ll q$ will follow from push-up (\cite[Lemma 1.22]{CG:12}). Since the length of $\gamma$
is strictly positive, the restriction of $\gamma$ to a suitable subinterval will also have
positive length while at the same time being contained in a cylindrical neighborhood
as in the definition of causally plain spacetimes
around some point on $\gamma$. So without loss of generality we may suppose that 
$\gamma$ is contained in such a neighborhood  $U$ around $p$. Now let us suppose that
$\gamma$ never enters $I^+(p,U)$, i.e., that it remains within $J^+(p,U)\setminus I^+(p,U)$.
By our assumption on $U$, \cite[Prop.\ 1.10 (v)]{CG:12} implies that $\partial J^+(p,U) = \partial \check I(p,U)
= \partial I(p,U)$ is given (in the cylindrical chart over $U$) as the graph of a Lipschitz function $f$.
It follows that $\gamma$ lies entirely within the graph of $f$. But then the proof of
\cite[Prop.\ 1.21 (v)]{CG:12} shows that $\dot \gamma (t)$ (if it exists) cannot be $g$-timelike 
at any $t$. This contradicts the fact that $L_g(\gamma)>0$.
\end{pr}

\begin{prop} \label{prop-caus-plain-lsc}
Let $(M,g)$ be causally plain. Then the time-separation function 
$\tau: M\times M \to [0,\infty]$ is lower semicontinuous.
\end{prop}
\begin{pr} We basically follow \cite[Lemma 4.4]{BEE:96}.
Let $p,q\in M$. If $\tau(p,q)=0$, there is nothing to prove. Next, let $0<\tau(p,q)<\infty$ and,
given any $0<\epsilon< \tau(p,q)$, choose a
future-directed causal curve $\gamma: [0,1]\to M$ from $p$ to $q$ with $L_g(\gamma)\ge \tau(p,q)-\frac{\eps}{2}$.
Now pick $0<t_1<t_2<1$ such that $0<L_g(\gamma|_{[0,t_1]})<\frac{\eps}{4}$ and $0<L_g(\gamma|_{[t_2,1]})<\frac{\eps}{4}$. 
Set $p_1:=\gamma(t_1)$ and $q_1:=\gamma(t_2)$, as well as $U:=I^-(p_1)$ and $V:=I^+(q_1)$. Since $L_g(\gamma|_{[0,t_1]})>0$, 
$\tau(p,p_1)>0$, which by Lemma \ref{con_met_tau} implies that $p\in I^-(p_1)$. Analogously, $V:=I^+(q_1)$ is an open 
neighborhood of $q$. For any $(p',q')\in U\times V$ we obtain
\begin{equation*}
\begin{split}
\tau(p',q') &\ge \tau(p',p_1) + \tau(p_1,q_1) + \tau(q_1,q') \ge L_g(\gamma|_{[t_1,t_2]})\\
&= L_g(\gamma) - L_g(\gamma|_{[0,t_1]}) - L_g(\gamma|_{[t_2,1]}) \ge \tau(p,q) - \eps,
\end{split}
\end{equation*}
and thereby lower semicontinuity of $\tau$ at $(p,q)$.

Finally, if $\tau(p,q)=\infty$, there are causal curves of arbitrary length from $p$ to $q$. Then the 
previous argument shows that the same is true of any points $p'$, $q'$ in $U$ respectively $V$.
\end{pr}   

Collecting the previous results we obtain:
\begin{prop}\label{prop-caus-plain-pls}
Let $(M,g)$ be a spacetime with a continuous causally plain metric. 
Then $(M,d^h,\ll,\le,\tau)$ is a Lorentzian pre-length space.
\end{prop}

In order for a spacetime with a continuous causally plain metric to give rise
to a Lorentzian length space there is one further requirement, namely 
we have to make sure that the notions of causal curve and of $R$-causal curve
coincide. In fact, as our next result will show, this is guaranteed if the spacetime is strongly causal.
In its proof, we will make use of {\em time functions}, i.e., functions that increase along any 
future-directed causal curve. In any continuous spacetime there always exist smooth local time
functions (e.g., the time coordinate in any cylindrical neighborhood).
\begin{prop}\label{prop-strongly-causal}
For a continuous and strongly causal spacetime $(M,g)$ the notions of causal and R-causal curve coincide.
\end{prop}
\begin{pr} By \cite[Lemma 3.21]{MS:08} (which remains valid for continuous spacetimes), strong causality
implies that any point in $M$ possesses a neighborhood basis consisting of causally convex neighborhoods
(i.e., such that any causal curve with endpoints contained in the neighborhood remains entirely within it).
Suppose now that $\gamma:I \to M$ is an R-causal curve and let $t_0\in I$. By \cite[Lemma 15]{BS:18} it suffices to show
that for any smooth local time function $f$ in a neighborhood $U$ of $\gamma(t_0)$, $f\circ \gamma$ is non-decreasing near $t_0$.
To this end, let $V\subseteq U$ be causally convex and let $J\subseteq I$ be an open interval around $t_0$
such that $\gamma(J)\subseteq V$. Then by R-causality, for $t_1<t_2$ in $J$ there exists a future-directed
causal curve $\sigma$ in $M$ (and consequently in $V$) from $\gamma(t_1)$ to $\gamma(t_2)$. By definition, 
$f$ is non-decreasing along $\sigma$, so in particular $f(\gamma(t_1))\le f(\gamma(t_2))$.
\end{pr}

Adding the assumption of strong causality, and defining $\mathcal{T}$ 
as in Definition \ref{bigtau}, we have:

\begin{lem}\label{lem-str-caus-tau-Tau} Let $(M,g)$ be a causally plain spacetime with a continuous, strongly causal 
metric and let $p,q\in M$. Then $\tau(p,q)=\mathcal{T}(p,q)$. Moreover, for any causal curve
$\gamma$, $L_g(\gamma)\le L_\tau(\gamma)$.
\end{lem}
\begin{pr} Note first that the claimed inequality on the lengths of causal curves
follows from the first part of the proof of Proposition \ref{prop-smo-stc-len}.

Also, $\mathcal{T}(p,q)\le \tau(p,q)$ was already shown in Remark \ref{lls_rem} \ref{lls_rem-3}.

To show the converse inequality suppose, to the contrary, that there exists some $\epsilon>0$
such that $\mathcal{T}(p,q)+\eps < \tau(p,q)$. Then by the definitions of $\mathcal{T}$ and $\tau$ there exists
some future-directed causal curve $\sigma$ from $p$ to $q$ such that for any
future-directed causal curve $\gamma$ we have
\[
L_\tau(\gamma) < \tau(p,q) - \eps  < L_g(\sigma) \le L_\tau(\sigma).
\]
Setting $\gamma:=\sigma$ now gives a contradiction.
\end{pr}
\begin{rem}
 \begin{itemize}
  \item[(i)] This result indicates that the length coming from a time separation function $\tau$ that in turn is defined via a length functional, is better behaved than $L_\tau$ for a generic $\tau$. In fact, in the above Lemma \ref{lem-str-caus-tau-Tau} we show $\tau= \Tau$ without knowing whether $L_\tau$ equals $L_g$.
 \item[(ii)] The assumption of causal plainness in the above Lemma \ref{lem-str-caus-tau-Tau} is in fact not needed for the proof. However, it was added as otherwise $\tau$ would not be lower semicontinuous in general and hence one would not obtain a \LpLSn.
 \end{itemize}
\end{rem}

Based on this we can proceed to proving the main result of this section:

\begin{thm} Let $(M,g)$ be a spacetime with a continuous, strongly causal and causally plain metric. 
Then $(M,d^h,\ll,\le,\tau)$ is a strongly localizable Lorentzian length space.
\end{thm}
\begin{pr}
Due to Proposition \ref{prop-caus-plain-pls}, Proposition \ref{prop-cau-loc-cc}, Lemma \ref{lem-str-caus-tau-Tau} and the fact 
that spacetimes are causally path connected by definition, it remains to establish that $(M,d^h,\ll,\le,\tau)$ is 
strongly localizable (Definition \ref{def-loc-LpLS}). We first note that by Proposition \ref{prop-strongly-causal} $g$-causal 
curves are the same as causal curves in the sense of Definition \ref{def-cau-cur}, hence we may speak of causal curves 
without ambiguity.

Let $p\in M$, then there is a neighborhood $U$ of $p$ such that the $h$-arclength of all causal curves in $U$ is bounded by 
some constant $C>0$ by \cite[Lemma 2.6.5]{Chr:11} or \cite[Lemma 2.1]{GLS:18}. This gives 
\ref{def-loc-LpLS} \ref{def-loc-LpLS-cau-comp}.

At this point let $\hat{g}$ be a smooth Lorentzian metric on $M$, with $g\prec\hat{g}$ (see \cite{CG:12}). Then there 
exists an (arbitrarily small)
$\hat{g}$-globally hyperbolic neighborhood $(V, x^\mu)$ of $p$ that is causally convex in $U$ by \cite[Thm.\ 
2.14]{MS:08} (cf.\ also \cite[Thm.\ 2.2]{SS:18}). This means that in the $x^\mu$-coordinates one has that $x^0=0$ is a 
Cauchy hypersurface in $V$ with respect to $\hat{g}$ and that any $\hat{g}$-causal curve in $U$ with endpoints in $V$ is 
contained in $V$. Then $x^0=0$ is also a Cauchy hypersurface with respect to $g$ and, hence $(V,g\rvert_V)$ is globally 
hyperbolic by \cite[Thm.\ 5.7]{Sae:16} and thus maximal (in $V$) causal curves exist between any two (in $V$) causally 
related points by the Avez-Seifert result for continuous metrics \cite[Prop.\ 6.4]{Sae:16}. By strong causality we can 
without loss of generality assume that $V$ is actually causally convex in $M$. Clearly, $I^\pm(x)\cap V \neq \emptyset$ for 
every $x\in V$.

Now define $\omega\colon V\times V\rightarrow [0,\infty)$ by $\omega(x,y):=\tau(x,y)$ for $x,y\in V$. Note that any causal 
curve from $x\in V$ to $y\in V$ is contained in $V$, hence $\tau(x,y)<\infty$ as there exists a maximal causal curve. This 
curve is actually globally maximal since $V$ is causally convex in $M$. This also implies that $\omega$ is continuous. It is 
lower semicontinuous by Proposition \ref{prop-caus-plain-lsc}, so assume that it were not upper semicontinuous at $(x,y)\in 
V\times V$. Thus, there is a $\delta>0$ and sequences $x_n\to x$, $y_n\to y$ in $V$ such that
\begin{equation*}
 \tau(x_n,y_n)\geq \tau(x,y)+\delta\,,
\end{equation*}
which implies $\tau(x_n,y_n)>0$. Consequently, there is a future-directed causal curve $\alpha_n$ from $x_n$ to $y_n$ with 
$L_g(\alpha_n) > \tau(x_n,y_n)-\frac{1}{n}$. Since $V$ is globally hyperbolic there is a limit curve $\alpha$ from $x$ to $y$
of the $\alpha_n$s (by the limit curve theorem \cite{Min:08a}), that satisfies $L_g(\alpha)\geq \tau(x,y)+\delta$ --- a 
contradiction. This shows that $\omega$ is continuous and establishes \ref{def-loc-LpLS} \ref{def-loc-LpLS-om-con}.

By the above we have that for any $x,y\in V$ with $x<y$ there is a globally maximal causal curve $\gamma_{x,y}$ from $x$ to 
$y$. Thus $\gamma_{x,y}$ is also $\tau$-maximal by Remark \ref{max_curve_rem} and so 
$L_\tau(\gamma_{x,y}) = \tau(x,y) = \omega(x,y)$, which establishes \ref{def-loc-LpLS} \ref{def-loc-LpLS-max-cc}.

\end{pr}

\begin{ex}\label{ex_reg_lls} There are large classes of spacetimes that are in fact regular Lorentzian length spaces, 
namely:
\begin{enumerate}[label=(\roman*)]
\item Strongly causal Lorentzian metrics $g$ of regularity at least $C^{1,1}$. Indeed, for such metrics property 
(iv) from Definition \ref{def-loc-LpLS} is an immediate consequence of the Gauss Lemma (cf.\ \cite{Min:15}
or \cite{KSSV:14}).
\item \label{ex_reg_lls-2} Continuous causally plain and strongly causal Lorentzian metrics whose causal curvature is bounded 
above. In fact, in this case Definition \ref{def-loc-LpLS} (iv) is satisfied by Remark \ref{rem-curvbound-vs-localizable}.
\end{enumerate}
\end{ex}

\subsection{Closed cone structures}
Many results from smooth causality theory can be generalized to cone structures on smooth manifolds. The interest
in such generalizations originated in the problem of constructing smooth time functions in stably causal or 
globally hyperbolic spacetimes. Fathi and Siconolfi in \cite{FS:12} studied continuous cone structures and
employed methods from weak KAM theory to address this problem. In \cite{BS:18}, Bernard and Suhr considered
closed cone structures and developed a theory of Lyapunov functions for such cone structures, 
showing, among other results, the equivalence between stable causality and
the existence of temporal functions, or between global hyperbolicity and the existence of
steep temporal functions in this setting.
The deepest and most comprehensive study of causality theory for closed cone structures to date is the very recent
work \cite{Min:17} by E.\ Minguzzi. It provides a complete causality theory, establishing the full causal ladder
for such cone structures, and contains manifold applications, among others to time and temporal functions, singularity
theorems, embedding of Lorentzian manifolds into Minkowski spacetime, and noncommutative geometry. 

In this section we follow the approach in \cite{Min:17} and show that closed cone structures
provide a rich source of examples of Lorentzian pre-length and length spaces.
We begin by recalling some basic definitions.

A {\em sharp cone} in a vector space $V$ is a subset of $V\setminus \{0\}$ that is positively
homogeneous, closed in the trace topology of $V$ on $V\setminus \{0\}$, convex, and
does not contain any line through $0$. It is called proper if its interior is non-empty.
A {\em cone structure} on a smooth manifold $M$ is a map $x\mapsto C_x$ that assigns to each
$x\in M$ a sharp non-empty cone. The cone structure is called {\em closed} if 
it forms a closed subbundle of the slit tangent bundle of $M$. It is called {\em proper}
if it is closed and $\mathrm{int}(C)_x\not=\emptyset$ for each $x\in M$.
(Semi-)continuity of cone structures is formulated in terms
of the Hausdorff distance on local sphere bundles, see \cite[Sec.\ 2]{FS:12}, \cite[Sec.\ 2]{Min:17}.

An absolutely continuous curve $\gamma$ in $M$ is called causal for the cone structure $C$
if $\dot \gamma(t)\in C_{\gamma(t)}$ almost everywhere. Timelike curves are by definition
piecewise $C^1$-solutions of the differential inclusion $\dot \gamma(t)\in \mathrm{int}(C)_{\gamma(t)}$.
Based on these notions, one defines the chronological and causal relations $\ll$ and $<$ as usual,
whereby any locally Lipschitz cone structure induces a causal space in the sense of Definition \ref{def-causal-space} (cf.\ \cite[Thm.\ 8]{Min:17}).

The following notions were introduced in \cite[Sec.\ 2.13]{Min:17}: Given a closed cone structure $(M,C)$
and a concave, positively homogeneous function $\mathcal{F}: C\to [0,\infty)$, a cone structure on
$M^\times := M\times \R$ is defined by 
\[
C^\times_{(p,r)} := \{(y,z) : y\in C_p, |z|\le \mathcal{F}(y)\}.
\]
The corresponding cone structure $(M^\times,C^\times)$ is called a {\em Lorentz-Finsler space} $(M,\mathcal{F})$.
The latter is called closed, respectively proper, respectively locally Lipschitz if $(M^\times,C^\times)$ has these properties.
On a closed Lorentz-Finsler space the length of a causal curve $\gamma:[0,1]\to M$ is defined by
$L(\gamma):=\int_0^1 \mathcal{F}(\dot \gamma)\,dt$. The corresponding {\em Lorentz-Finsler distance} is 
defined as $\tau(p,q)= 0$ if $p\not\le q$, and $\tau(p,q):=\sup_\gamma L(\gamma)$ otherwise, where the supremum
is taken over all future-directed causal curves from $p$ to $q$. Finally, fix a complete Riemannian metric
$h$ on $M$. Then we have:
\begin{prop}
Let $(M,\mathcal{F})$ be a locally Lipschitz proper Lorentz-Finsler space such that $\mathcal{F}(\partial C)=0$.
Then $(M,d^h,\ll,\le,\tau)$ is a Lorentzian pre-length space.
\end{prop}
\begin{pr}
In fact, under these assumptions $\tau$ is lower semi-continuous by \cite[Thm.\ 52]{Min:17}. Furthermore,
$\tau(p,q)>0\Leftrightarrow p\ll q$ follows from \cite[Thm.\ 55]{Min:17}. The other properties
are immediate from the definitions.
\end{pr}
As in the case of continuous spacetimes, if we want to proceed to establishing the properties
of Lorentzian length spaces we first have to secure that the classes of R-causal and causal curves
coincide. In fact, this is true for any strongly causal closed cone structure, hence in particular
for any strongly causal proper Lorentz-Finsler space: since the existence of arbitrarily small 
causally convex neighborhoods in this case holds by definition, this follows exactly as in the proof
of Proposition \ref{prop-strongly-causal}. Moreover, defining $\mathcal{T}$ as in Definition \ref{bigtau},
the same proof as in Lemma \ref{lem-str-caus-tau-Tau} gives:
\begin{lem}\label{lem-lor-fins-tau-Tau}
Let $(M,\mathcal{F})$ be a locally Lipschitz proper Lorentz-Finsler space such that $\mathcal{F}(\partial C)=0$
and suppose that $(M,\mathcal{F})$ is strongly causal. Then for all $p, q \in M$, $\mathcal{T}(p,q)=\tau(p,q)$.
\end{lem}
Finally, we can prove:
\begin{thm}
Let $(M,\mathcal{F})$ be a locally Lipschitz proper Lorentz-Finsler space such that $\mathcal{F}(\partial C)=0$.
If $(M,\mathcal{F})$ is strongly causal, then $(M,d^h,\ll,\le,\tau)$ is a strongly localizable Lorentzian length space.
\end{thm}
\begin{pr}
By strong causality and \cite[Prop.\ 7]{Min:17}, any point $x$ in $M$ possesses a basis
of open neighborhoods that are globally hyperbolic and causally convex. We may therefore
pick such a neighborhood $\Omega_x$ such that $\Omega_x$ is causally closed.
Since $M$ is causally path connected by definition and taking into account Lemma \ref{lem-lor-fins-tau-Tau},
it only remains to establish strong localizability, i.e., properties (i)--(iii) from Definition
\ref{def-loc-LpLS} for a neighborhood basis. Now for any $\Omega_x$,
(i) follows from \cite[Prop.\ 7]{Min:17}, and to obtain (ii), by \cite[Thm.\ 52 and Thm.\ 58]{Min:17}, 
we may set $\omega_x:=\tau|_{\Omega_x\times \Omega_x}$. 
Finally, (iii) follows from the Avez-Seifert theorem \cite[Thm.\ 54]{Min:17} (or also from Theorem 
\ref{thm-gh-LLS-geo}). 
\end{pr}

Moreover, just as in Example \ref{ex_reg_lls} \ref{ex_reg_lls-2}, regularity and (SR)-localizability
can be achieved by assuming upper causal curvature bounds.

Finally, a natural open question is whether one can weaken the assumption of 
local Lipschitz continuity of 
the cone structure by an analogue of causal plainness, as in the case of continuous spacetimes.
\subsection{Outlook on further examples}\label{subsec-app-qg}
The framework developed in the previous sections makes it possible to handle situations where one might not have the 
structure of a manifold or Lorentz(-Finsler) metric. Even in these cases the theory of Lorentzian (pre-) length spaces allows 
one to define timelike and causal curvature (bounds) via triangle comparison. Thus, 
it provides a new perspective on curvature in such cases where there is no classical notion of curvature (Riemann 
tensor, Ricci and sectional curvature, etc.). This is applies, in particular, to certain
approaches to quantum gravity, as pointed out in e.g.\ \cite{MP:06} (see also the corresponding paragraph in the 
introduction, Section \ref{sec-intro}). In this non-rigorous outlook we briefly sketch two such approaches, namely \emph{causal 
Fermion systems} \cite{Fin:16, Fin:18} and the theory of causal 
sets \cite{BLMS:87}.

The underlying idea in both cases is that the structure of spacetime has to be modified on a microscopic scale to 
include quantum effects. This gives rise to non-smoothness of the underlying geometry, and only in the macroscopic picture the 
classical spacetime (i.e., a Lorentzian manifold) emerges. We briefly recall the relevant definitions, the causal structures 
and discuss the connections to Lorentzian (pre-)length spaces.
\bigskip

We start with the recent approach of causal Fermion systems.
Let $H$ be a separable complex Hilbert space and let $n\in\N$. Let $F\subseteq L(H)$ be the set of all self-adjoint operators 
on $H$ of finite rank that have at most $n$ positive and $n$ negative eigenvalues. Let $\rho$ be a measure defined on a 
$\sigma$-algebra of subsets of $F$, called the \emph{universal measure}. Then $(H,F,\rho)$ is called a \emph{causal Fermion 
system}. The \emph{spacetime} $M$ is defined as $M:=\supp(\rho)\subseteq F$, the support of the universal measure $\rho$, 
with the induced topology from $L(H)$. The causal structure arises as follows: For $x,y\in M$ the product $x y:=x\circ y$ has 
at most $2n$ eigenvalues. If all of them have the same absolute value, $x$ and $y$ are called \emph{spacelike} separated. If 
they do not all have the same absolute value and are real, then $x$ and $y$ are called \emph{timelike} separated. In all 
other cases $x$ and $y$ are called \emph{lightlike} separated. There is also a notion of time orientation as follows: For an 
operator $x\in M$ denote by $\pi_x$ the orthogonal projection on the subspace $x(H)$. Define the anti-symmetric function 
$C\colon M\times M\rightarrow \R$ by $C(x,y):= i\, \mathrm{tr}(y x \pi_y \pi_x - x y \pi_x \pi_y)$. One can therefore define that
\emph{$y$ lies to the future (past) of $x$} if $C(x,y)>0$ ($C(x,y)<0$, respectively). Now one can define \emph{timelike and 
causal chains} and the causal relations $x\ll y$ if there is (a future-directed) timelike chain from $x$ to $y$ and 
analogously $x<y$ for causal chains. At this point one is able to introduce a \emph{Lorentzian distance} on $M$ and obtain 
the structure of a \LpLSn, as described in \cite[Subsec.\ 5.1]{Fin:18}. Whether this gives the structure of a \LLS 
will be considered elsewhere.

\bigskip

Another approach to quantum gravity is the theory of \emph{causal sets}, which is closely related to \LpLSn s and Example 
\ref{ex-lpls-gra}.

Let $(X,\leq)$ be a partially ordered set that is locally finite, i.e., for every $x,y\in X$ the set $J(x,y):=\{z\in X: 
x\leq z\leq y\}$ is finite. Writing $x<y$ if $x\leq y$ and $x\neq y$ we define $I(x,y):=\{z\in X: 
x<z< y\}$. This minimal framework induces an analogous notion of geodesics or maximal curves as follows: $(x,y)\in 
X^2$ is called a \emph{link} if $x\leq y$ and $I(x,y)=\emptyset$. A \emph{chain} is a sequence of points $(x_i)_{i=1}^n$ 
with $x_i<x_{i+1}$ for $i=1,\ldots,n-1$, and moreover, a chain is a \emph{path} if every pair $(x_i,x_{i+1})$ is a link. The 
length of a chain $C=(x_i)_{i=1}^n$ is $l(C)=n$ and a \emph{geodesic} between $x$ and $y$ (for $x<y$) is a path from $x$ to 
$y$ whose length is maximal over all paths from $x$ to $y$.

To include causal sets in the framework of \LpLSn s, define $\ll:= <$, i.e., $x\ll y$ if and only if $x<y$. Then, $(X,\ll,\leq)$ 
is a causal space (in the sense of Definition \ref{def-causal-space}) and we define $\tau(x,y):=\sup\{l(C): C$ is a chain from 
$x$ to $y\}$ and $\tau(x,y)=0$ if there is no chain from $x$ to $y$. Putting any metric $d$ on $X$ that makes $\tau$ 
lower semicontinuous yields a \LpLS \Xlln. Since $J(x,y)$ is finite for every $x,y\in X$ there is no metric on $X$ 
that allows continuous parametrizations of $J(x,y)$ and thus it cannot be turned into a \LLSn. However, there is a close 
connection and it seems possible to discretize a \LLS as one can discretize a Lorentzian manifold to obtain a causal 
set.

\noindent
{\bf Acknowledgments.} 
The authors are grateful to Tobias Beran and Shantanu Dave for helpful discussions in various stages of this project and to Stephanie Alexander, 
Felix Finster, Stacey Harris, and Ettore Minguzzi  
for several insightful discussions and comments,
in particular during the meetings \emph{Non-regular spacetime
geometry}, Firenze, June 20-22, 2017, and \emph{Advances in General Relativity}, ESI
Vienna, August 28 - September 1, 2017.

This work was supported by research grants P26859 and P28770 of the Austrian Science Fund FWF.

\addcontentsline{toc}{section}{References}

%\newpage
\bibliographystyle{alphaabbr}
\bibliography{lls}

\end{document}